\documentclass[11pt, a4paper]{amsart}
\usepackage{fullpage} 
\usepackage{amsmath}
\usepackage{amsthm}
\usepackage{amssymb}
\usepackage[only,mapsfrom]{stmaryrd}
\usepackage{graphicx}
\usepackage{cancel} 
\usepackage{color} 
\usepackage[all]{xy} 
\usepackage{tikz-cd}
\usepackage{tikz}
\usetikzlibrary{arrows}
\usepackage{nicefrac}
\usepackage{enumitem}
\usepackage{hyperref}

\newtheorem{thm}{Theorem}[section]
\newtheorem{lem}[thm]{Lemma}
\newtheorem{prop}[thm]{Proposition}
\newtheorem{cor}[thm]{Corollary}

\newtheorem*{thm*}{Theorem}
\newtheorem*{conj*}{Conjecture}
\newtheorem*{cor*}{Corollary}
\newtheorem*{ques*}{Question}

\newtheorem*{namedthm}{\namedthmname}

\theoremstyle{definition}
\newtheorem{rem}[thm]{Remark}
\newtheorem{defn}[thm]{Definition} 
\newtheorem*{rem*}{Remark}

\DeclareMathOperator{\Lk}{Lk}
\DeclareMathOperator{\St}{St}
     
\newcommand{\BC}{{\mathbb {C}}} \newcommand{\BD}{{\mathbb {D}}}

 \newcommand{\BL}{{\mathbb {L}}}
 \newcommand{\BN}{{\mathbb {N}}}
 
 \newcommand{\BR}{{\mathbb {R}}}

 \newcommand{\BZ}{{\mathbb {Z}}}

\newcommand{\CI}{{\mathcal {I}}} 
 
\newcommand{\CM}{{\mathcal {M}}}

\newcommand{\id}{{{\rm id}}}

\def\Proof{\noindent{\bf Proof.}\quad}
\def\sProof{\noindent{\bf Sketch of Proof. }\quad}
\def\qed{\hfill$\square$\smallskip}

\newcommand{\elbraigecpx}{\mathcal{EB}}
\newcommand{\marc}{\mathcal{MA}}
\newcommand{\lmarc}{\mathcal{LMA}}
\newcommand{\cmarc}{\mathcal{CMA}}

\let\Bigast\Bigjoin

\DeclareMathOperator{\dlk}{{\Lk}{\downarrow}}
\DeclareMathOperator{\dst}{{\St}{\downarrow}}
\newcommand{\elpbraigecpx}{\mathcal{EPB}}
\newcommand{\elcbraigecpx}{\mathcal{ECB}}

\newcommand{\spraige}{\mathcal{S}}
\newcommand{\Poset}{\mathcal{P}}

\DeclareMathOperator{\Bot}{bot}

\title{Finiteness properties for relatives of braided Higman--Thompson groups}

\author{Rachel Skipper}   
\address{D\'{e}partement de Math\'{e}matiques et Applications, \'{E}cole normale sup\'{e}rieure, 45 Rue d'Ulm, 75005 Paris, France} 
\email{rachel.skipper@ens.fr}

\author{Xiaolei Wu}
\address{Shanghai Center for Mathematical Sciences, Jiangwan Campus, Fudan University, No.2005 Songhu Road, Shanghai, 200438, P.R. China}
\email{xiaoleiwu@fudan.edu.cn}



\date{March 2023}

\begin{document}

\begin{abstract}
We study the finiteness properties of the braided Higman--Thompson group $bV_{d,r}(H)$ with labels in $H\leq B_d$, and $bF_{d,r}(H)$ and $bT_{d,r}(H)$ with labels in $H\leq PB_d$ where $B_d$ is the braid group with $d$ strings and $PB_d$ is its pure braid subgroup. We show that for all $d\geq 2$ and $r\geq 1$, the group $bV_{d,r}(H)$ (resp. $bT_{d,r}(H)$ or $bF_{d,r}(H)$) is of type $F_n$ if and only if $H$ is.  Our result in particular confirms a recent conjecture of Aroca and Cumplido. 
 \end{abstract}

\maketitle

\section{Introduction}

The family of Thompson's groups and the many groups in the extended Thompson family have long been studied for their many interesting properties. Thompson's group $F$ is the first example of a type $F_\infty$, torsion-free group with infinite cohomological dimension \cite{BG84}, while Thompson's groups $T$ and $V$ provided the first examples of finitely presented simple groups. More recently, the braided and labeled braided Higman--Thompson groups have garnered attention in part due their connections with big mapping class groups. 

The braided version of Thompson's group $V$, which we refer to here as $bV$, was first introduced independently by Brin and Dehornoy \cite{Br07}, \cite{Deh06}. Brady, Burillo, Cleary, and Stein introduced braided $F$, or $bF$. The groups $bV$ and $bF$ were shown to be finitely presented in \cite{Br06} and \cite{BBCS08}, respectively, and this was extended to show that both of these groups are of type $F_\infty$ in \cite{BFM+16}. Braided $T$ was mentioned in \cite{BFM+16} and shown to be of type $F_\infty$ in \cite{Wi19}. The ribbon version of Thompson's group $V$ was first constructed by Thumann and proved to be type $F_\infty$ as well in \cite{Th17}. Moving to higher dimensions, Spahn showed the braided Brin--Thompson groups are of type $F_\infty$ \cite{Spa21}. Recently, Aroca and Cumplido \cite{AC20} broadened the definitions of braided groups in the extended Thompson family to what we will refer to as labeled braided Higman--Thompson groups $bV_{d,r}(H)$, which depend on a choice of a subgroup $H$ of the braid group $B_d$. For any subgroup $H$ of the pure braid group $PB_d$, it is natural to also consider $bF_{d,r}(H)$ and $bT_{d,r}(H)$ which we do here. Aroca and Cumplido prove that if $H\leq B_d$ is finitely generated, the groups $bV_{d,r}(H)$ are finitely generated for all $r\geq 1$ and $d\geq 2$.
They conjectured the following in \cite[Section 4.3]{AC20}.
\begin{conj*}
For all $d\geq 2$ and $r\geq 1$, the group $bV_{d,r}(H)$ is finitely presented when $H$ is finitely presented. 
\end{conj*}

Recall that a group $G$ is of \emph{type $F_n$} if there exists an aspherical CW-complex whose fundamental group is $G$ and whose $n$-skeleton is finite. Being of type $F_1$ is equivalent to the group being finitely generated and type $F_2$ is equivalent to the group being finitely presented.  A group is of \emph{type $F_\infty$} if it is of type $F_n$ for all $n\geq 1$. Our first theorem confirms Aroca and Cumplido's conjecture. In fact we completely determine the relationship between the finiteness properties of $H$ and that of the braided Higman--Thompson groups with labels in $H$. 

\begin{thm*}[\ref{thm-fin-bthomp}]
For any $d\geq 2$ and $r\geq 1$  and any subgroup $H$ of the braid group $B_d$ (resp. of the pure braid group $PB_d$), the group $bV_{d,r}(H)$ (resp.  $bT_{d,r}(H)$ or $bF_{d,r}(H)$) is of type $F_n$ if and only if $H$ is.
\end{thm*}

\begin{rem*}
The collection of subgroups of $B_d$ represents a class of groups with rich finiteness properties. In fact, Zaremsky showed in \cite{Za17} that there exists a subgroup of $PB_d$ which is of type $F_n$ but not $F_{n+1}$ for any $0\leq n\leq d-3$. In particular,  our theorem provides a new class of Thompson-like groups which is of type $F_n$ but not of type $F_{n+1}$ for each $n$. 
\end{rem*}

When $H$ is the trivial group, the groups $bV_{d,r}(H)$, $bF_{d,r}(H)$, and $bT_{d,r}(H)$ are the braided Higman--Thompson groups $bV_{d,r}$, $bF_{d,r},$ and $bT_{d,r}$, hence we have the following.

\begin{cor*}[\ref{cor-BVBTBFtypeFinfty}]
The braided Higman--Thompson groups $bV_{d,r}$, $bF_{d,r},$ and $bT_{d,r}$ are of type $F_{\infty}$.
\end{cor*}

\begin{rem*}
Genevois, Lonjou and Urech  in \cite{GLU20} introduced another braided version of the Higman--Thompson group $T_{d,r}$ and proved that the groups they study are of type $F_{\infty}$ as well. Their groups are different from the ones studied here as they naturally surject onto $T_{d,r}$  with kernel being an infinite braid group, while our group $bT_{d,r}$ surjects naturally onto $T_{d,r}$ with an infinite pure braid group as the kernel.
\end{rem*}

View the braid group $B_d$ as the mapping class group of the disk with $d$ marked points and let $C$ be the subgroup of $B_d$ generated by the half twist around the boundary. Then the corresponding group $bV_{d,r}(C)$ can be identified with the ribbon Higman--Thompson group $RV_{d,r}$. See Proposition \ref{lem-idf-rb-br} for a precise statement. Note also that when we take the label group to be the index $2$ subgroup of $C$ which is generated by a full Dehn twist around the boundary, we get the oriented ribbon Higman--Thompson groups $RV^+_{d,r},RF^+_{d,r},$ and $RT^+_{d,r}$.

\begin{cor*}[\ref{cor-rthompFinfty}]
The ribbon Higman--Thompson group $RV_{d,r}$ is of type $F_\infty$. Likewise, the oriented ribbon Higman--Thompson groups $RV^+_{d,r},RF^+_{d,r},$ and $RT^+_{d,r}$ are of type $F_\infty$.
\end{cor*}

\begin{rem*}
For a more thorough exploration of the ribbon Higman--Thompson groups $RV_{d,r}$ and oriented ribbon Higman--Thompson groups $RV^+_{d,r}$, we direct the reader to \cite{SW21A} where the authors showed these families of groups satisfy homological stability. 
\end{rem*}

There is a large amount of literature devoted to finding finiteness properties of groups in the extended family of Thompson's groups. Most often the groups are of type $F_\infty$, e.g \cite{belk16,Br87,BFM+16,farley15, FMWZ13, martinez-perez16,nucinkis18,sz17,Th17}, though not always, e.g., Belk--Forrest's basilica Thompson group $T_B$ \cite{belk15} is type $F_1$, but not $F_2$ \cite{witzel16} and the simple groups of type $F_{n}$ but not $F_{n+1}$ are given in \cite{swz19} and \cite{BZ20}.  
 
The ``only if" part of Theorem \ref{thm-fin-bthomp} is proved using a quasi-retract argument inspired by \cite[Section 4]{BZ20}. For the ``if" part, as in \cite{BFM+16}, our proof uses Brown's criterion. Ultimately, it reduces to proving that certain $d$-marked-point-disk complexes are highly connected. See Sections \ref{Sct:conn-dcpx} and \ref{Sct:conn-i-ldcpx} for the details. Given a surface $S$ with $m$ marked points, a $k$-simplex in the $d$-marked-point-disk complex  $\BD_{d}(S)$ is an isotopy class of a system of disjointly embedded disks $\langle D_0,D_1,\cdots, D_k\rangle$ such that each disk $D_i$ encloses precisely $d$ marked points in its interior. The face relation is given by the subset relation.  Note that except some singular cases, $\BD_d(S)$ can be viewed as a full subcomplex of the curve complex first defined by Harvey \cite{Hav81}. The connectivity properties of the curve complex played an important role in Harer's proof of homological stability for the mapping class groups \cite{Har85}.  We have the following.

\begin{thm*}[\ref{cor:conn-i-dcomplex}]
Let $S$ be a surface with $m$ marked points. Then for any $d\geq2$, the complex $\BD_{d}(S)$ is $(\lfloor\frac{m+1}{2d-1}\rfloor-2)$-connected.
\end{thm*}

\subsection*{Outline of paper}
In Section \ref{sec:connectivitytools}, we describe the connectivity tools that will be necessary for the remainder of the paper. In Section \ref{sec:braidedgroups}, we introduce the definition of the labeled braided Higman--Thompson groups using braided paired forest diagrams to define the elements. Next, in Section \ref{sec:finitenessprop}, we build the Stein space on which the labeled braided Higman--Thompson groups act and use it to prove the ``if" part of Theorem \ref{thm-fin-bthomp}  by applying a combination of Brown's criterion with Bestvina-Brady discrete Morse theory. In the same section, we then prove the ``only if" part by a quasi-retract argument. 

\subsection*{Notation and convention.} All surfaces in this paper are assumed to be connected and orientable unless otherwise stated. Given a simplicial complex $X$ and a cell $\sigma \in X$, we denote the link of $\sigma$ in $X$ by $\Lk_X(\sigma)$ (resp. the star of $\sigma$ by $\St_X (\sigma)$). When the situation is clear, we quite often omit $X$ and simply denote the link by $\Lk (\sigma)$ and the star by $\St (\sigma)$.   We also use the convention that  $(-1)$-connected means non-empty and that every space is $(-2)$-connected. In particular, the empty set is $(-2)$-connected. Finally, we adopt the convention that elements in groups are multiplied from left to right.

\subsection*{Acknowledgements.}  
The first part of this project was done while the first author was a visitor in the Unit\'{e} de math\'{e}matiques pures et appliqu\'{e}es at the ENS de Lyon and during a visit to the University of Bonn. She thanks them for their hospitality. She was also supported by the GIF, grant I-198-304.1-2015, ``Geometric exponents of random walks and intermediate growth groups" and NSF DMS--2005297 ``Group Actions on Trees and Boundaries of Trees". This project also has received funding from the European Research Council (ERC) under the European Union’s Horizon 2020 research and innovation program (grant agreement No.725773).

Part of this work was done when the second author was a member of the  Hausdorff Center of Mathematics. At the time, he was  supported by Wolfgang L\"uck's ERC Advanced Grant “KL2MG-interactions”
(no. 662400) and the DFG Grant under Germany's Excellence Strategy - GZ 2047/1, Projekt-ID 390685813.

Part of this work was also done when both authors were visiting IMPAN in Warsaw during the Simons Semester ``Geometric and Analytic Group Theory" which was partially supported by the grant 346300 for IMPAN from the Simons Foundation and the matching 2015-2019 Polish MNiSW fund. We would like to thank Kai-Uwe Bux for inviting us for a research visit in Bielefeld in May 2019 and many stimulating discussions. Special thanks go to Jonas Flechsig for his comments on preliminary versions of the paper and many useful discussions regarding the finiteness properties of braided Higman-Thompson groups. Furthermore, we want to thank Javier Aramayona and Stefan Witzel for discussions, Andrea Bianchi for comments and Matthew Zaremsky for some helpful communications and comments.

\section{Connectivity tools}\label{sec:connectivitytools}
In this section, we review some of the connectivity tools that we need for calculating the connectivity of our spaces. A good reference is \cite[Section 2]{HV17}, although not all the tools we use can be found there.

\subsection{Discrete Morse theory}
Let $Y$ be a piecewise Euclidean cell complex, and let~$h$ be a map
from the set of vertices of $Y$ to the integers, such that each cell has a unique vertex maximizing~$h$.  Call $h$ a
\emph{height function}, and $h(y)$ the \emph{height} of $y$ for
vertices $y$ in $Y$.  For $t\in\BZ$, define $Y^{\leq t}$ to be
the full subcomplex of~$Y$ spanned by vertices~$y$ satisfying $h(y)\leq
t$.  Similarly, define $Y^{< t}$ and $Y^{=t}$.  The \emph{descending star} $\dst(y)$ of a
vertex $y$ is defined to be the open star of~$y$ in~$Y^{\le h(y)}$.  The
\emph{descending link} $\dlk(y)$ of $y$ is given by the set of ``local
directions'' starting at $y$ and pointing into $\dst(y)$. More details can be found in \cite{BB97}, and the following Morse Lemma is a consequence of \cite[Corollary~2.6]{BB97}.

\begin{lem}[Morse lemma]\label{lemm-Morse}
Let $Y$ be a piecewise Euclidean cell complex and let $h$ be a height function on $Y$.
\begin{enumerate}
    \item Suppose that for any vertex $y$ with $h(y)=t$, $\dlk (y)$ is $(k-1)$-connected. Then the pair $(Y^{\leq t}, Y^{<t})$ is $k$-connected.
    \item Suppose that for  any vertex $y$ with $h(y) \geq t$, $\dlk (y)$ is $(k-1)$-connected. Then $(Y,Y^{<t})$ is $k$-connected.
\end{enumerate}
\end{lem}

Recall that we say a pair of spaces $(X,Y)$ with $Y\subseteq X$ is $k$-connected if the inclusion map $Y \hookrightarrow X$ induces an isomorphism in $\pi_j$ for $j< k$ and an epimorphism in $\pi_k$.

\subsection{Complete join}\label{subsec-com-join} The complete join is another useful tool introduced by Hatcher and Wahl in \cite[Section 3]{HW10} for proving connectivity results. We review the basics here.

\begin{defn}
A surjective simplicial map $\pi: Y\to X$ is called a \emph{complete join} if it satisfies the following properties:
\begin{enumerate} [label=(\arabic*)]
    \item $\pi$ is injective on individual simplices. 
    \item For each $p$-simplex $\sigma = \langle v_0,\cdots,v_p\rangle$ of $X$, $\pi^{-1}(\sigma)$ is the join $\pi^{-1}(v_0)\ast \pi^{-1}(v_1)\ast\cdots\ast \pi^{-1}(v_p)$.
 \end{enumerate}
\end{defn}


\begin{defn}
A simplicial complex  $X$ is called weakly Cohen-Macaulay of dimension $n$ if $X$ is $(n-1)$-connected and the link of each $p$-simplex of $X$ is $(n-p-2)$-connected. We sometimes shorten  weakly Cohen-Macaulay to $wCM$. 
\end{defn}

The main result regarding complete join that we will use in this paper is the following statement.

\begin{prop}\cite[Proposition 3.5]{HW10} \label{prop-join-conn}
If $Y$ is a complete join complex over a  $wCM$ complex $X$ of dimension $n$, then $Y$ is also  $wCM$ of dimension $n$.
\end{prop}

\begin{rem}\label{rem-cjoin}
If $\pi: Y\to X$ is a complete join, then $X$ is a retract of $Y$. In fact, we can define a simplicial map $s:X\to Y$ such that $\pi\circ s = \id_X$ by sending a vertex $v\in X$ to any vertex in $\pi^{-1}(v)$ and then extending it to simplices. The fact that $s$ can be extended to simplices is granted by the condition that $\pi$ is a complete join.  In particular we can also conclude that if $Y$ is $n$-connected, so is $X$.
\end{rem}

\subsection{The mutual link trick}
In the proof of \cite[Theorem 3.10]{BFM+16}, there is a beautiful argument for resolving intersections of arcs inspired by Hatcher's flow argument \cite{Hat91}. They attributed the idea to Andrew Putman. Recall that Hatcher's flow argument allows one to ``flow" a complex to its subcomplex. But in the process, one can only ``flow" a vertex to a new one in its link. The mutual link trick will allow one to ``flow" a vertex to a new one not in its link, provided ``the mutual link" is sufficiently connected.

To apply the  mutual link trick, we first need a lemma that allows us to homotope a simplicial map to a simplexwise injective one \cite[Lemma 3.9]{BFM+16}.  Recall a simplicial map is called \emph{simplexwise injective} if its
restriction to any simplex is injective. See also \cite[Section 2.1]{GRW18} for more information.

\begin{lem}
\label{lem:injectifying}
Let $Y$ be a compact $m$-dimensional combinatorial manifold.  Let $X$ be a
simplicial complex and assume that the link of every $p$-simplex in
$X$ is $(m-p-2)$-connected.  Let $\psi \colon Y \to X$ be a
simplicial map whose restriction to $\partial Y$ is simplexwise
injective.  Then after possibly subdividing the simplicial structure of $Y$, $\psi$ is
homotopic relative $\partial Y$ to a simplexwise injective map.
\end{lem}

Note that, as discussed in \cite[Lemma 5.19]{GLU20}, there is a mistake in the connectivity bound given in \cite{BFM+16} that has been corrected here. 

\begin{lem}[The mutual link trick]\label{lemma-replace-trick}
Let $Y$ be a closed $m$-dimensional combinatorial manifold and $f: Y\to X$ be a simplexwise injective simplicial map. Let $y \in Y$ be a vertex and $f(y) = x$ for some $x\in X$. Suppose $x'$ is another vertex of $X$ satisfying the following condition:
\begin{enumerate}
    \item $f(\Lk_{Y}(y)) \leq \Lk_{X} (x')$,
    \item the mutual link $\Lk_{X} (x) \cap \Lk_{X}(x')$ is $(m-1)$-connected,
\end{enumerate}
Then we can define a new simplexwise injective map $g:Y\to X $ by sending $y$ to $x'$ and all the other vertices $y'$ to $f(y')$ such that $g$ is homotopic to $f$. 
\end{lem}

See \cite[Lemma 1.9]{SW21A} for a proof of the mutual link trick.

\section{Higman--Thompson groups and their braided versions}\label{sec:braidedgroups}
In this section, we first give an introduction to the Higman--Thompson groups and then introduce their braided version. The braided Thompson-like groups in the generality we will consider here were first given by Aroca and Cumplido in \cite{AC20}. Note that Aroca and  Cumplido's exposition closely follows the original introduction of braided Higman--Thompson groups by Brin \cite{Br07} whereas we instead will follow the exposition in \cite{BFM+16}.

\subsection{Higman--Thompson groups} The Higman--Thompson groups were first introduced by Higman as a generalization of the groups \cite{Hi74} given earlier in handwritten, unpublished notes of Thompson.
First let us recall the definition of the Higman--Thompson groups. Although there are a number of equivalent definitions of these groups, we will use the notion of paired forest diagrams. First we define a \emph{finite rooted $d$-ary tree} to be a finite tree such that every vertex has degree $d+1$ except the \emph{leaves} which have degree 1, and the \emph{root}, which has degree $d$ (or degree $1$ if the root is also a leaf). Usually we draw such trees with the root at the top and the nodes descending from it down to the leaves. A vertex $v$ of the tree along with its $d$ adjacent descendants will be called a \emph{caret}. If the leaves of a caret in the tree are leaves of the tree, we will call the caret \emph{elementary}.  A collection of $r$ many $d$-ary trees will be called a $\emph{$(d,r)$-forest}$. When $d$ is clear from the context, we may just call it an $r$-forest.

Define a \emph{paired $(d,r)$-forest diagram} to be a triple $(F_-,\rho,F_+)$
consisting of two $(d,r)$-forests $F_-$ and $F_+$ both with $l$ leaves for some $l$, and a permutation $\rho \in S_l$, the symmetric group on $l$ elements.  We label the leaves
of~$F_-$ with $1,\dots,l$ from left to right, and for each $i$,
the $\rho(i)^{\text{th}}$ leaf of~$F_+$ is labeled $i$.

\begin{figure}[h]\label{fig:ele-V}
\centering
\begin{tikzpicture}[line width=1pt, scale=0.5]
\begin{scope}[xshift= -4cm]
  \draw
   (-12,-2) -- (-9,0) -- (-6,-2)
   (-9,0) -- (-9,-2)
   (-11,-1.33) -- (-11,-2)
   (-11,-1.33) -- (-10,-2)
   (-4,-2) -- (-1,0) -- (2,-2)
   (-1,0)  -- (-1,-2);

\draw[dotted] (3.5,0)  -- (3.5,-3);

  \filldraw
  (-12,-2) circle (1.5pt)
  (-9,0) circle (1.5pt)
   (-6,-2) circle (1.5pt)
   (-9,-2) circle (1.5pt)
  (-11,-1.33)  circle (1.5pt)
   (-11,-2) circle (1.5pt)
   (-10,-2) circle (1.5pt)
    (-4,-2) circle (1.5pt)
   (-1,0) circle (1.5pt)
   (-1,-2) circle (1.5pt)
   (2,-2) circle (1.5pt);
  \node at  (-12,-2.5) {$1$};
  \node at (-11,-2.5) {$2$};
  \node at (-10,-2.5) {$3$};
  \node at (-9,-2.5) {$4$};
  \node at (-6,-2.5) {$5$};
  \node at (-4,-2.5) {$6$};
  \node at (-1,-2.5) {$7$};
  \node at (2,-2.5) {$8$};
\end{scope}

  \begin{scope}[ xshift=13cm]
   \draw
   (-12,-2) -- (-9,0) -- (-6,-2)
   (-9,0) -- (-9,-2)
   
  (-8,-2) -- (-9,-1.33) -- (-10,-2)
   
   (-4,-2) -- (-1,0) -- (2,-2)
   (-1,0)  -- (-1,-2);

  \filldraw
  (-12,-2) circle (1.5pt)
  (-9,0) circle (1.5pt)
   (-6,-2) circle (1.5pt)
   (-9,-2) circle (1.5pt)
  (-9,-1.33)  circle (1.5pt)
   (-10,-2) circle (1.5pt)
   (-8,-2) circle (1.5pt)
    (-4,-2) circle (1.5pt)
   (-1,0) circle (1.5pt)
   (2,-2) circle (1.5pt)
   (-1,-2) circle (1.5pt);
  \node at  (-12,-2.5) {$4$};
  \node at (-6,-2.5) {$6$};
  \node at (-10,-2.5) {$1$};
  \node at (-8,-2.5) {$3$};
  \node at (-9,-2.5) {$2$};
  \node at (-4,-2.5) {$7$};
  \node at (2,-2.5) {$5$};
  \node at (-1,-2.5) {$8$};
 \end{scope}

  \begin{scope}[ yshift=-4.5cm,xshift= -4cm]
    \draw
   (-12,-2) -- (-9,0) -- (-6,-2)
   (-9,0) -- (-9,-2)
   (-4,-2) -- (-1,0) -- (2,-2)
   (-1,0)  -- (-1,-2);

\draw[dotted] (3.5,0)  -- (3.5,-3);

  \filldraw
  (-12,-2) circle (1.5pt)
  (-9,0) circle (1.5pt)
   (-6,-2) circle (1.5pt)
   (-9,-2) circle (1.5pt)
    (-4,-2) circle (1.5pt)
   (-1,0) circle (1.5pt)
   (-1,-2) circle (1.5pt)
   (2,-2) circle (1.5pt);
  \node at  (-12,-2.5) {$1$};
  \node at (-9,-2.5) {$2$};
  \node at (-6,-2.5) {$3$};
  \node at (-4,-2.5) {$4$};
  \node at (-1,-2.5) {$5$};
  \node at (2,-2.5) {$6$};
  \end{scope}

  \begin{scope}[xshift=13cm, yshift=-4.5cm]
   \draw
   (-12,-2) -- (-9,0) -- (-6,-2)
   (-9,0) -- (-9,-2)
   
   
   (-4,-2) -- (-1,0) -- (2,-2)
   (-1,0)  -- (-1,-2);

  \filldraw
  (-12,-2) circle (1.5pt)
  (-9,0) circle (1.5pt)
   (-6,-2) circle (1.5pt)
   (-9,-2) circle (1.5pt)
    (-4,-2) circle (1.5pt)
   (-1,0) circle (1.5pt)
   (2,-2) circle (1.5pt)
   (-1,-2) circle (1.5pt);
  \node at  (-12,-2.5) {$2$};
  \node at (-6,-2.5) {$4$};
  \node at (-9,-2.5) {$1$};
  \node at (-4,-2.5) {$5$};
  \node at (2,-2.5) {$3$};
  \node at (-1,-2.5) {$6$};
  \end{scope}
\end{tikzpicture}

\caption{Reduction, of the top paired $(3,2)$-forest diagram to the bottom one.}
\label{fig:reduction_V}
\end{figure}
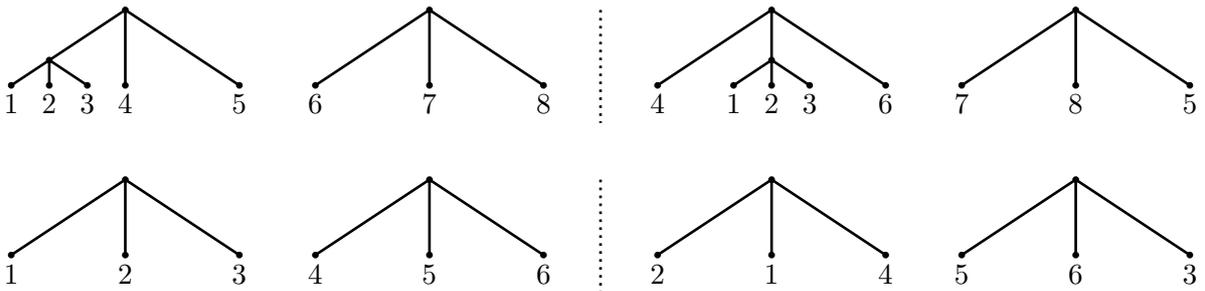

Define a \emph{reduction} of a paired $(d,r)$-forest diagram to be the following: Suppose there is an
elementary caret in~$F_-$ with leaves labeled by $i,\cdots,i+d-1$ from left to right, and an elementary caret in~$F_+$ with  leaves labeled by ~$i, \cdots,i+d-1$ from left to right.  Then we can ``reduce'' the diagram
by removing those carets, renumbering the leaves and replacing
$\rho$ with the permutation~$\rho'\in S_{l-d+1}$ that sends the new leaf
of~$F_-$ to the new leaf of~$F_+$, and otherwise behaves like~$\rho$.
The resulting paired forest diagram~$(F'_-,\rho',F'_+)$ is then said to
be obtained by \emph{reducing}~$(F_-,\rho, F_+)$. See Figure~\ref{fig:reduction_V} below for an idea of reduction
of paired $(3,2)$-forest diagrams. The reverse
operation to reduction is called \emph{expansion}, so $(F_-,\rho,F_+)$
is an expansion of $(F'_-,\rho',F'_+)$.  A paired forest diagram is
called \emph{reduced} if there is no reduction possible.  Define an equivalence relation on the set of paired $(d,r)$-forest diagrams by declaring two paired  forest diagrams to be equivalent if one can be reached by the other through a finite series of reductions and expansions.
Thus an equivalence class of paired forest diagrams consists of all diagrams
having a common reduced representative.  Such reduced representatives
are unique. 

There is a binary operation $\ast$ on the set of equivalence classes
of paired $(d,r)$-forest diagrams.  Let $\alpha =(F_-,\rho,F_+)$ and $\beta=(E_-,\xi,E_+)$
be reduced paired forest diagrams.  By applying repeated expansions to
$\alpha$ and $\beta$, we can find representatives $(F'_-,\rho',F'_+)$ and
$(E'_-,\xi',E'_+)$ of the equivalence classes of $\alpha$ and
$\beta$, respectively, such that~$F'_+ = E'_-$.  Then we
declare $ \alpha \ast  \beta $ to be $(F'_-,\rho'\xi',E'_+)$.  This operation is well defined on the equivalence classes and is a group operation.

\begin{defn}
    \label{def:Higman-V-F-T}
   The \emph{Higman--Thompson group} $V_{d,r}$ is the group of equivalence classes of paired
   $(d,r)$-forest diagrams with the multiplication~$\ast$.  The Higman--Thompson group $F_{d,r}$ is
    the subgroup of~$V_{d,r}$ consisting of elements where the permutation
    is the identity.  The Higman--Thompson group $T_{d,r}$ is
    the subgroup of~$V_{d,r}$ consisting of elements where the permutation
    is cyclic, i.e. there exists some $k$ such that for all $i$, the $i$-th leaf is mapped to the $(i+k)$-th leaf (modulo the number of leaves).
\end{defn}

The usual Thompson's groups $F,T,$ and $V$ are special cases of Higman--Thompson groups. In fact, $F=F_{2,1}$, $T= T_{2,1}$, and $V= V_{2,1}$. Brown and Geoghegan showed in \cite{BG84} that $F$ is of type $F_\infty$ which provided the first example of a torsion-free group of type $F_\infty$ but not of finite cohomological dimension. Later, in \cite[Section 4]{Br87} Brown showed the following.

\begin{thm}
The Higman--Thompson groups  $V_{d,r}$,  $F_{d,r}$ and  $T_{d,r}$  are all of type $F_\infty$.
\end{thm}

\subsection{Braided Higman--Thompson groups with labels}\label{subsect:braidedHTgroups} In this subsection, we introduce braided and labeled braided Higman--Thompson groups. Again, we follow the exposition in \cite[Section 1]{BFM+16} closely to define these groups.

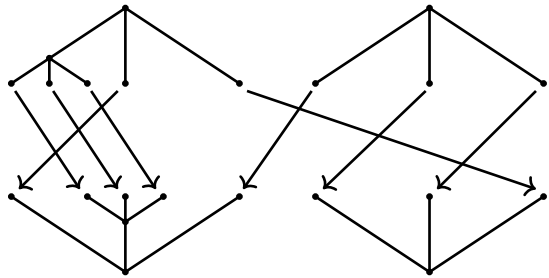
\begin{figure}[h]
\centering
\begin{tikzpicture} [line width=1pt, scale = 0.5]
 \begin{scope}[xshift =10cm]
 
 \draw
   (-12,-2) -- (-9,0) -- (-6,-2)
   (-9,0) -- (-9,-2)
   (-11,-1.33) -- (-11,-2)
   (-11,-1.33) -- (-10,-2)
   (-4,-2) -- (-1,0) -- (2,-2)
   (-1,0)  -- (-1,-2);
   \draw[->](-11.9,-2.2) -> (-10.2,-4.8);
      \draw[->](-10.9,-2.2) -> (-9.2,-4.8);
      \draw[->](-9.9,-2.2) -> (-8.2,-4.8);
            \draw[->](-9.2,-2.2) -> (-11.8,-4.8);
      
        \draw[->](-5.8,-2.2) -> (1.8,-4.8);
         \draw[->](-4.1,-2.2) -> (-5.9,-4.8);
        \draw[->](-1.1,-2.2) -> (-3.8,-4.8);
          \draw[->](1.8,-2.2) -> (-0.8,-4.8);
      
  \filldraw
  (-12,-2) circle (1.5pt)
  (-9,0) circle (1.5pt)
   (-6,-2) circle (1.5pt)
   (-9,-2) circle (1.5pt)
  (-11,-1.33)  circle (1.5pt)
   (-11,-2) circle (1.5pt)
   (-10,-2) circle (1.5pt)
    (-4,-2) circle (1.5pt)
   (-1,0) circle (1.5pt)
   (-1,-2) circle (1.5pt)
   (2,-2) circle (1.5pt);

\end{scope}

 \begin{scope}[ xshift = 10cm, yscale=-1,  yshift=7cm]
  \draw
   (-12,-2) -- (-9,0) -- (-6,-2)
   (-9,0) -- (-9,-2)
   
  (-8,-2) -- (-9,-1.33) -- (-10,-2)
   
   (-4,-2) -- (-1,0) -- (2,-2)
   (-1,0)  -- (-1,-2);

  \filldraw
  (-12,-2) circle (1.5pt)
  (-9,0) circle (1.5pt)
   (-6,-2) circle (1.5pt)
   (-9,-2) circle (1.5pt)
  (-9,-1.33)  circle (1.5pt)
   (-10,-2) circle (1.5pt)
   (-8,-2) circle (1.5pt)
    (-4,-2) circle (1.5pt)
   (-1,0) circle (1.5pt)
   (-1,-2) circle (1.5pt)
   (2,-2) circle (1.5pt);
 \end{scope}   
   
\end{tikzpicture}
\caption{An element of $V_{3,2}$.}
\label{fig:element_of_V}
\end{figure}

For convenience, we will think of 
the forest $F_+$ drawn beneath $F_-$ and upside down, i.e., with the
root at the bottom and the leaves at the top.  The permutation $\rho$
is then indicated by arrows pointing from the leaves of $F_-$ to the
corresponding paired leaves of $F_+$.  See Figure~\ref{fig:element_of_V}
for this visualization of (the unreduced representation of) the
element of $V_{3,2}$ from Figure~\ref{fig:reduction_V}.

Now in the braided version of the Higman--Thompson groups, the permutations of leaves are simply
replaced by braids between the leaves. We will need to go one step farther to define the group $bV_{d,r}(H)$. Here we further replace the permutations by labeled braids as described in the next definition. Recall that an element in the braid group $B_d$ consists of $d$ strings. We enumerate them by their initial points from left to right as $1,2,\cdots,d$.

\begin{defn}\label{defn-label}
    Given any group $H$, an element in the \emph{labeled braid group} $B_l(H)$ is an ordered pair $(b,\lambda)$ where $b$ is in the braid group $B_l$ and $\lambda:\{1,2,\cdots,l\}\to H$ is a map called the labeling map. The group operation is given by stacking the two braids and multiplying the labels using the multiplication in $H$.  In other words, $B_l(H)\cong B_n\ltimes H^l$, where $B_l$ acts on $H^l$ by permuting the coordinates using the canonical map $\rho:B_l\to S_l$ and $S_l$ is the symmetric group on $l$ elements.
\end{defn}

\begin{defn}
    For a group $H$, a \emph{braided paired $(H,d,r)$-forest diagram} is a triple $(F_-,(b,\lambda),F_+)$
    consisting of two $(d,r)$-forests $F_-$ and $F_+$ both with $l$ leaves for some $l$
     and a labeled braid~$(b,\lambda) \in B_l(H)$.
\end{defn}

\begin{figure}[h]
\centering
\begin{tikzpicture}[line width=1pt]
 \begin{scope}[yscale=-1]
\begin{scope}[scale =0.4]
  \draw
   (-12,-2) -- (-9,0) -- (-6,-2)
   (-9,0) -- (-9,-2)
   (-11,-1.33) -- (-11,-2)
   (-11,-1.33) -- (-10,-2)
   (-4,-2) -- (-1,0) -- (2,-2)
   (-1,0)  -- (-1,-2);


 \draw  
(-9,-2) to [out=-80, in=80] (-12,-5);

\draw[white, line width=4pt]  
      (-12,-2) to [out=-90, in=90] (-10,-5) 
   (-11,-2) to [out=-90, in=90] (-9,-5) 
   (-10,-2) to [out=-90, in=90] (-8,-5); 
   
  \draw
   (-12,-2) to [out=-90, in=90] (-10,-5) 
   (-11,-2) to [out=-90, in=90] (-9,-5) 
   (-10,-2) to [out=-90, in=90] (-8,-5);
 
\draw
(-4,-2) to [out=-90, in=90] (-6,-5)
(-1,-2) to [out=-90, in=90] (-4,-5)
(2,-2) to [out=-90, in=90] (-1,-5);

\draw[white, line width=4pt]  
(-6,-2) to [out=-90, in=90] (2,-5);   

\draw  
(-6,-2) to [out=-90, in=90] (2,-5);


  \filldraw
  (-12,-2) circle (1.5pt)
  (-9,0) circle (1.5pt)
   (-6,-2) circle (1.5pt)
   (-9,-2) circle (1.5pt)
  (-11,-1.33)  circle (1.5pt)
   (-11,-2) circle (1.5pt)
   (-10,-2) circle (1.5pt)
    (-4,-2) circle (1.5pt)
   (-1,0) circle (1.5pt)
   (-1,-2) circle (1.5pt)
   (2,-2) circle (1.5pt);
   
   \begin{scope}[yscale=-1,  yshift=7cm]
  \draw
   (-12,-2) -- (-9,0) -- (-6,-2)
   (-9,0) -- (-9,-2)
   
  (-8,-2) -- (-9,-1.33) -- (-10,-2)
   
   (-4,-2) -- (-1,0) -- (2,-2)
   (-1,0)  -- (-1,-2);

\draw[white,fill=white] (-10.2,-2.7) circle (10pt);

\draw (-10.2,-2.7) circle (10pt) node[text=black, scale=.5] {$h_2$};

\draw[white,fill=white] (-9.2,-2.7) circle (10pt);

\draw  (-9.2,-2.7) circle (10pt) node[text=black, scale=.5] {$h_2$};

\draw[white,fill=white] (-8.2,-2.7) circle (10pt);

\draw  (-8.2,-2.7) circle (10pt) node[text=black, scale=.5] {$h_2$};

\draw[white,fill=white] (-11.7,-2.7) circle (10pt);

\draw  (-11.7,-2.7) circle (10pt) node[text=black, scale=.5] {$h_1$};

\draw[white,fill=white] (-5.8,-2.7) circle (10pt);

\draw  (-5.8,-2.7) circle (10pt) node[text=black, scale=.5] {$h_3$};

\draw[white,fill=white] (-3.8,-2.7) circle (10pt);

\draw  (-3.8,-2.7) circle (10pt) node[text=black, scale=.5] {$h_4$};

\draw[white,fill=white] (-0.8,-2.7) circle (10pt);

\draw  (-0.8,-2.7) circle (10pt) node[text=black, scale=.5] {$h_5$};

\draw[white,fill=white] (1.7,-2.7) circle (10pt);

\draw  (1.7,-2.7) circle (10pt) node[text=black, scale=.5] {$h_6$};

 \filldraw[fill=white, draw=black] (-12.2,-4.1) rectangle (-9.5,-4.7);

\draw  (-10.85,-4.4)  node[text=black, scale=.5] {$h_2$};
  \filldraw
  (-12,-2) circle (1.5pt)
  (-9,0) circle (1.5pt)
   (-6,-2) circle (1.5pt)
   (-9,-2) circle (1.5pt)
  (-9,-1.33)  circle (1.5pt)
   (-10,-2) circle (1.5pt)
   (-8,-2) circle (1.5pt)
    (-4,-2) circle (1.5pt)
   (-1,0) circle (1.5pt)
   (-1,-2) circle (1.5pt)
   (2,-2) circle (1.5pt);
 \end{scope} 
   
   \draw[->]
   (3,-3) -> (4,-3);
  
 \end{scope}  
   
\begin{scope}[xshift= 7cm,scale =0.4]

  \draw
   (-12,-2) -- (-9,0) -- (-6,-2)
   (-9,0) -- (-9,-2)
   (-4,-2) -- (-1,0) -- (2,-2)
   (-1,0)  -- (-1,-2);


 \draw  
(-9,-2) to [out=-90, in=90] (-12,-5);

\draw[white, line width=4pt]  

   (-12,-2) to [out=-90, in=90] (-9,-5); 
   
  \draw

   (-12,-2) to [out=-90, in=90] (-9,-5);
 
\draw
(-4,-2) to [out=-90, in=90] (-6,-5)
(-1,-2) to [out=-90, in=90] (-4,-5)
(2,-2) to [out=-90, in=90] (-1,-5); 
\draw[white, line width=4pt]  
(-6,-2) to [out=-90, in=90] (2,-5);   

\draw  
(-6,-2) to [out=-90, in=90] (2,-5);


  \filldraw
  (-12,-2) circle (1.5pt)
  (-9,0) circle (1.5pt)
   (-6,-2) circle (1.5pt)
   (-9,-2) circle (1.5pt)
    (-4,-2) circle (1.5pt)
   (-1,0) circle (1.5pt)
   (-1,-2) circle (1.5pt)
   (2,-2) circle (1.5pt);
   
   \begin{scope}[yscale=-1,  yshift=7cm]
  \draw
   (-12,-2) -- (-9,0) -- (-6,-2)
   (-9,0) -- (-9,-2)
   
   
   (-4,-2) -- (-1,0) -- (2,-2)
   (-1,0)  -- (-1,-2);

\draw[white,fill=white] (-9.2,-2.7) circle (10pt);

\draw  (-9.2,-2.7) circle (10pt) node[text=black, scale=.5] {$h_2$};

\draw[white,fill=white] (-11.7,-2.7) circle (10pt);

\draw  (-11.7,-2.7) circle (10pt) node[text=black, scale=.5] {$h_1$};

\draw[white,fill=white] (-5.8,-2.7) circle (10pt);

\draw  (-5.8,-2.7) circle (10pt) node[text=black, scale=.5] {$h_3$};

\draw[white,fill=white] (-3.8,-2.7) circle (10pt);

\draw  (-3.8,-2.7) circle (10pt) node[text=black, scale=.5] {$h_4$};

\draw[white,fill=white] (-0.8,-2.7) circle (10pt);

\draw  (-0.8,-2.7) circle (10pt) node[text=black, scale=.5] {$h_5$};

\draw[white,fill=white] (1.7,-2.7) circle (10pt);

\draw  (1.7,-2.7) circle (10pt) node[text=black, scale=.5] {$h_6$};

  \filldraw
  (-12,-2) circle (1.5pt)
  (-9,0) circle (1.5pt)
   (-6,-2) circle (1.5pt)
   (-9,-2) circle (1.5pt)
    (-4,-2) circle (1.5pt)
   (-1,0) circle (1.5pt)
   (-1,-2) circle (1.5pt)
   (2,-2) circle (1.5pt);
 \end{scope} 
\end{scope}

\end{scope} 
\end{tikzpicture}
\caption{Reduction of braided paired forest diagrams.}
\label{fig:reduction_Vbr}
\end{figure}
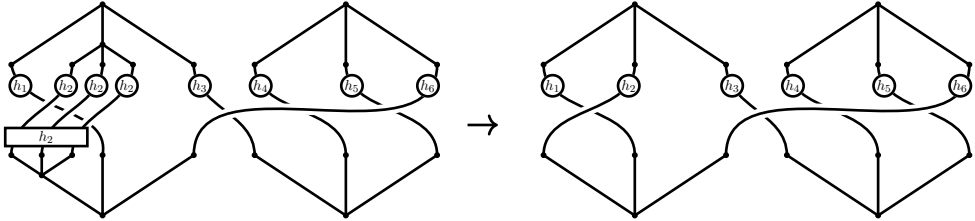

We draw braided paired forest diagrams with $F_+$ upside down and below
$F_-$ with the strands of the braid connecting the leaves and with each stand labeled by an element in $H$.  This is
analogous to the visualization of paired forest diagrams in
Figure~\ref{fig:element_of_V} and examples of braided paired forest
diagrams can be seen in Figure~\ref{fig:reduction_Vbr}.

Now to define the group $bV_{d,r}(H)$, we will restrict ourselves to the case $H\leq B_d$, although the definition works as long as we have a homomorphism $s:H\to B_d$.

As in the Higman--Thompson group case, we can define an equivalence relation on the set of
braided paired forest diagrams using the notions of reduction and
expansion.  This time, it is easier to first define expansion and then take
reduction as the reverse of expansion.  Let $\rho_b\in S_l$ denote the
permutation corresponding to the braid $b\in B_l$.  Let $(F_-,(b,\lambda),F_+)$
be a braided paired forest diagram.  Label the leaves of $F_-$ from $1$
to $l$, left to right, and for each $i$ label the
$\rho_b(i)^{\text{th}}$ leaf of~$F_+$ by~$i$.  By the~$i^{\text{th}}$
strand of the braid we will always mean the strand that begins at
the~$i^{\text{th}}$ leaf of $F_-$, i.e., we count the strands from the
top. The label for the $i^{\text{th}}$
strand of the braid is given by $\lambda(i)$. An \emph{expansion} of $(F_-,(b, \lambda),F_+)$ is the following:
For some $1\le i\le l$, replace~$F_\pm$ with forests $F_\pm'$ obtained from
$F_\pm$ by adding a caret to the leaf labeled~$i$.  Then replace $b$
with a braid $b' \in B_{l+d-1}$, obtained from replacing the
$i^{\text{th}}$ strand of $b$ with the braid $\lambda(i)$. Finally, we label the $d$ new strands all by $\lambda(i)$. We denote the new labeling system by $\lambda'$ so that the triple $(F_-',(b',\lambda'),F_+')$ is an
\emph{expansion} of $(F_-,(b, \lambda),F_+)$. As with paired forest diagrams,
\emph{reduction} is the reverse of expansion, so $(F_-,(b, \lambda) ,F_+)$ is a reduction of
$(F_-',(b', \lambda'),F_+')$.  See
Figure~\ref{fig:reduction_Vbr} for an idea of reduction of braided
paired forest diagrams. Note that in the picture, we draw a small circle on each string inside which we write the corresponding label and we use a box with a label to indicate that the corresponding strings are braided according to that label.

Two braided paired forest diagrams with labels in $H$ are equivalent if one is obtained
from the other by a sequence of reductions or expansions.  The
multiplication operation~$\ast$ on the equivalence classes is defined
the same way as for~$V_{d,r}$. In more detail, let $\alpha =(F_-,(b_1,\lambda_1),F_+)$ and $\beta=(E_-,(b_2,\lambda_2),E_+)$
be reduced braided paired $(H,d,r)$-forest diagrams. By applying repeated expansions to
$\alpha$ and $\beta$ we can find representatives $(F'_-,(b_1',\lambda_1'),F'_+)$ and
$(E'_-,(b_2',\lambda_2'),E'_+)$ of the equivalence classes of $\alpha$ and
$\beta$, respectively, such that~$F'_+ = E'_-$.  Then we
declare $  \alpha \ast \beta $ to be $(F'_-, (b_1',\lambda_1') (b_2',\lambda_2'),E'_+)$.  This operation is well defined on the equivalence classes and is a group operation as proved in \cite[Section 3]{AC20}.

\begin{defn}\label{defn-bht-lb}
    Given any subgroup $H\leq B_d$, the \emph{braided Higman--Thompson group} $bV_{d,r}(H)$ is the group of equivalence classes of braided paired $(H,d,r)$-forests diagrams with the
    multiplication~$\ast$. For any $H\leq PB_d$, the \emph{braided Higman--Thompson group} $bF_{d,r}(H)$ is the group of equivalence classes of braided paired $(H,d,r)$-forest diagrams where the braids are all pure. Finally, for any $H\leq PB_d$, the \emph{braided Higman--Thompson group} $bT_{d,r}(H)$ is the group of equivalence classes of braided paired $(H,d,r)$-forest diagrams where the permutations corresponding to the braids are all cyclic.
\end{defn}

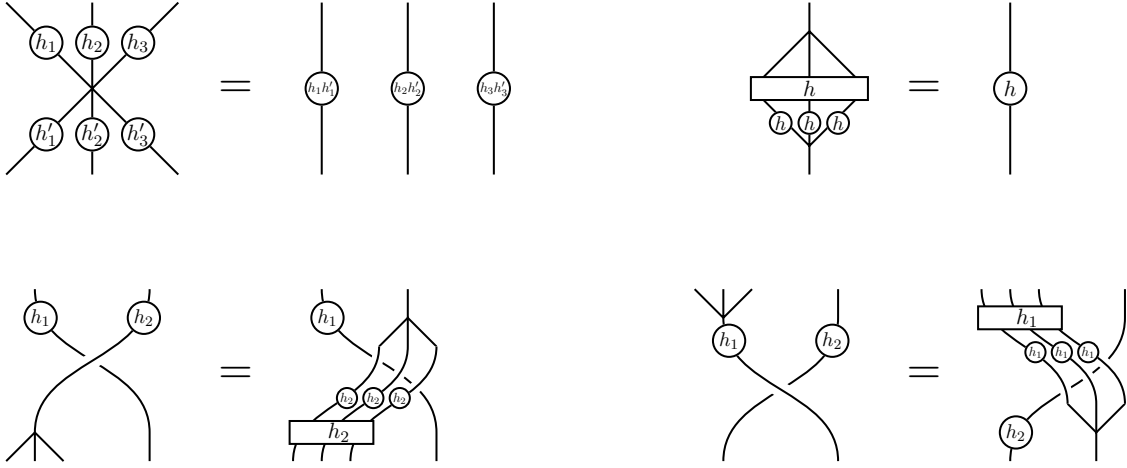
\begin{figure}[h]
\centering

\resizebox{15cm}{6.2cm}{%
\begin{tikzpicture}[line width=.5pt] 

  \draw
   (0,0.75) -- (1.5,-0.75)
   (0.75,0.75) -- (0.75,-0.75)
   (1.5,0.75) -- (0,-0.75);
  \node at (2,0) {$=$};
  
      \draw[white,fill=white] (0.35,.4) circle (4pt);

\draw  (0.35,.4)  circle (4pt) node[text=black, scale=.5] {$h_1$};

    \draw[white,fill=white] (0.35,-.4) circle (4pt);

\draw  (0.35,-.4)  circle (4pt) node[text=black, scale=.5] {$h_1'$};

      \draw[white,fill=white] (0.75, .4) circle (4pt);

\draw (0.75, .4)  circle (4pt) node[text=black, scale=.5] {$h_2$};

    \draw[white,fill=white] (0.75, -.4) circle (4pt);

\draw  (0.75, -.4)  circle (4pt) node[text=black, scale=.5] {$h_2'$};

      \draw[white,fill=white] (1.15, .4) circle (4pt);

\draw (1.15, .4)  circle (4pt) node[text=black, scale=.5] {$h_3$};

    \draw[white,fill=white] (1.15, -.4) circle (4pt);

\draw  (1.15, -.4)  circle (4pt) node[text=black, scale=.5] {$h_3'$};

 \begin{scope}[xshift=2.75cm]
  \draw
   (0,0.75) -- (0,-0.75)
    (0.75,0.75) -- (0.75,-0.75)
   (1.5,0.75) -- (1.5,-0.75);
   
    \draw[white,fill=white] (0,0) circle (4pt);

\draw   (0,0) circle (4pt) node[text=black, scale=.31] {$h_1h_1'$};

    \draw[white,fill=white] (0.75,0) circle (4pt);

\draw   (0.75,0) circle (4pt) node[text=black, scale=.31] {$h_2h_2'$};

    \draw[white,fill=white] (1.5,0) circle (4pt);

\draw   (1.5,0) circle (4pt) node[text=black, scale=.31] {$h_3h_3'$};

 \end{scope}
 \begin{scope}[xshift=6.5cm]
  \draw
   (0.5,0.75) -- (0.5,0.5) -- (0,0) -- (0.5,-0.5) -- (1,0) -- (0.5,0.5)
   (0.5,-0.5) -- (0.5,0.5)
(0.5,-0.5) -- (0.5,-0.75);
  \node at (1.5,0) {$=$};
  \draw
   (2.25,0.75) -- (2.25,-0.75);
   
      \draw[white,fill=white] (0.25,-.3) circle (2.7pt);

\draw    (0.25,-.3)  circle (2.7pt) node[text=black, scale=.45] {$h$};

    \draw[white,fill=white] (0.5,-.3) circle (2.7pt);

\draw   (0.5,-.3) circle (2.7pt) node[text=black, scale=.45] {$h$};

    \draw[white,fill=white] (0.75,-.3) circle (2.7pt);

\draw   (0.75,-.3) circle (2.7pt) node[text=black, scale=.45] {$h$}; 
   
       \draw[white,fill=white] (2.25,0) circle (4pt);

 \filldraw[fill=white, draw=black] (-0.01,-.1) rectangle (1.01,0.1);
\draw  (0.5,0)  node[text=black, scale=.5] {$h$};

\draw   (2.25,0) circle (4pt) node[text=black, scale=.5] {$h$};

 \end{scope}

 \begin{scope}[yshift=-2.5cm]
  \draw
   (0.25,0.75) to [out=-90, in=90, looseness=1] (1.25,-0.5) -- (1.25,-0.75);
  \draw[white, line width=4pt]
   (1.25,0.75) to [out=-90, in=90, looseness=1] (0.25,-0.5);
  \draw
   (1.25,0.75) to [out=-90, in=90, looseness=1] (0.25,-0.5);
  \draw
   (0,-0.75) -- (0.25,-0.5) -- (0.5,-0.75)
   (0.25,-0.5) --(0.25,-0.75);
   
\draw[white,fill=white]  (0.3,0.5) circle (4pt);
\draw   (0.3,0.5) circle (4pt) node[text=black, scale=.45] {$h_1$};
   
   \draw[white,fill=white] (1.2,0.5) circle (4pt);
\draw  (1.2,0.5) circle (4pt) node[text=black, scale=.45] {$h_2$};
   
  \node at (2,0) {$=$};
 \end{scope}

 \begin{scope}[yshift=-2.5cm, xshift=2.5cm]
  \draw
   (0.25,0.75) to [out=-90, in=90, looseness=1] (1.25,-0.5) -- (1.25,-0.75);
  \draw
   (0.75,0.25) -- (1,0.5) -- (1.25,0.25)  
   (1,0.5) -- (1,0.75);
  \draw[white, line width=3pt]
   (0.75,0.25) to [out=-90, in=90, looseness=1] (0,-0.75)
   (1.25,0.25) to [out=-90, in=90, looseness=1] (0.5,-0.75)
   (1,0.25) to [out=-90, in=90, looseness=1] (0.25,-0.75);
  \draw
   (0.75,0.25) to [out=-90, in=90, looseness=1] (0,-0.75)
   (1.25,0.25) to [out=-90, in=90, looseness=1] (0.5,-0.75)
    (1,0.5) -- (1,0.25) to [out=-90, in=90, looseness=1] (0.25,-0.75);
    
    \draw[white,fill=white]  (0.3,0.5) circle (4pt);
\draw   (0.3,0.5) circle (4pt) node[text=black, scale=.45] {$h_1$};
    
 \filldraw[fill=white, draw=black] (-0.03,-.6) rectangle (.7,-0.4);
\draw  (0.4,-0.5)  node[text=black, scale=.5] {$h_2$};

   \draw[white,fill=white] (.93,-.2) circle (2.5pt);
\draw  (.93,-.2) circle (2.5pt) node[text=black, scale=.3] {$h_2$};

   \draw[white,fill=white] (.47,-.2) circle (2.5pt);
\draw  (.47,-.2) circle (2.5pt) node[text=black, scale=.3] {$h_2$};

   \draw[white,fill=white] (.7,-.2) circle (2.5pt);
\draw  (.7,-.2) circle (2.5pt) node[text=black, scale=.3] {$h_2$};

 \end{scope}

 \begin{scope}[yshift=-2.5cm, xshift=6cm, yscale=-1]
  \draw
   (0.25,0.75) to [out=-90, in=90, looseness=1] (1.25,-0.5) -- (1.25,-0.75);
  \draw[white, line width=4pt]
   (1.25,0.75) to [out=-90, in=90, looseness=1] (0.25,-0.5);
  \draw
   (1.25,0.75) to [out=-90, in=90, looseness=1] (0.25,-0.5);
  \draw
   (0,-0.75) -- (0.25,-0.5) -- (0.5,-0.75)
    (0.25,-0.75) -- (0.25,-0.5);
    
    \draw[white,fill=white]  (0.3,-0.3) circle (4pt);
\draw   (0.3,-0.3) circle (4pt) node[text=black, scale=.45] {$h_1$};
   
   \draw[white,fill=white] (1.2,-.3) circle (4pt);
\draw  (1.2,-0.3) circle (4pt) node[text=black, scale=.45] {$h_2$}; 
  \node at (2,0) {$=$};
 \end{scope}

 \begin{scope}[yshift=-2.5cm, xshift=8.5cm, yscale=-1]
  \draw
   (0.25,0.75) to [out=-90, in=90, looseness=1] (1.25,-0.5) -- (1.25,-0.75);
  \draw
   (0.75,0.25) -- (1,0.5) -- (1.25,0.25)   (1,0.5) -- (1,0.75);
  \draw[white, line width=3pt]
   (0.75,0.25) to [out=-90, in=90, looseness=1] (0,-0.75)
   
   (1,0.25) to [out=-90, in=90, looseness=1] (0.25,-0.75)
   
   (1.25,0.25) to [out=-90, in=90, looseness=1] (0.5,-0.75);
  \draw
   (0.75,0.25) to [out=-90, in=90, looseness=1] (0,-0.75)
  (1,0.5) -- (1,0.25)   to [out=-90, in=90, looseness=1] (0.25,-0.75) 
   (1.25,0.25) to [out=-90, in=90, looseness=1] (0.5,-0.75);
   
    \filldraw[fill=white, draw=black] (-0.03,-.6) rectangle (.7,-0.4);
\draw  (0.4,-0.5)  node[text=black, scale=.5] {$h_1$};      

    \draw[white,fill=white]  (0.3,0.5) circle (4pt);
\draw   (0.3,0.5) circle (4pt) node[text=black, scale=.45] {$h_2$};

   \draw[white,fill=white] (.93,-.2) circle (2.5pt);
\draw  (.93,-.2) circle (2.5pt) node[text=black, scale=.3] {$h_1$};

   \draw[white,fill=white] (.47,-.2) circle (2.5pt);
\draw  (.47,-.2) circle (2.5pt) node[text=black, scale=.3] {$h_1$};

   \draw[white,fill=white] (.7,-.2) circle (2.5pt);
\draw  (.7,-.2) circle (2.5pt) node[text=black, scale=.3] {$h_1$};

 \end{scope}
\end{tikzpicture}
}
\caption{Moves to reduce braided paired forest diagrams after stacking.}
\label{fig:reduction_moves}

\end{figure}

A convenient way to visualize multiplication in~$bV_{d,r}(H), bF_{d,r}(H)$, and $bT_{d,r}(H)$ is via ``stacking'' braided paired forest diagrams. For $g,h$ in $bV_{d,r}(H), bF_{d,r}(H)$, or $bT_{d,r}(H)$, each pictured as a forest-braid-forest as before,~$g\ast h$ is obtained by attaching the top of $h$ to the bottom of~$g$ and then reducing the picture via certain moves. We indicate four of these moves in Figure~\ref{fig:reduction_moves} for $d=3$. A merge followed immediately by a split, or a split followed immediately by a merge, is equivalent to doing nothing except multiplying the labels, as seen in the top two pictures. Also, splits and merges interact with braids in the ways indicated by the bottom two pictures. We leave it to the reader to further inspect the details of this visualization of multiplication in these groups. This is closely related to the strand diagram model for Thompson's groups in \cite{BM14}. See also \cite[Section 1.2]{Br07}.

From now on, we will just refer to the braided $(H,d,r)$-forest diagrams as being the elements of~$bV_{d,r}(H), bF_{d,r}(H)$, or $bT_{d,r}(H)$, though one should keep in mind that the
elements are actually equivalence classes under the reduction and expansion operations. When $H$ is the trivial group, we denote the groups simply by $bV_{d,r}, bF_{d,r}$, or $bT_{d,r}$.

\begin{thm}\cite{BFM+16, Wi19}\label{thm-fin-tms}
The groups $bV_{2,1}, bF_{2,1}$, and $bT_{2,1}$  are of type $F_{\infty}$.
\end{thm}

Another interesting class of relatives of the braided Thompson groups is the ribbon Higman--Thompson groups. Let us explain this in more detail.

\begin{defn}
Let  $\CI = \amalg_{i=1}^l I_i: [0,1] \times \{1,\cdots,l\} \to \BR^2 $ be an embedding which we refer to as the \emph{marked bands}. A \emph{ribbon braid} is a map  $R: ([0,1] \times \{0,1,\cdots,l\})  \times [0,1] \to  \BR^2$  such that for any $0\leq t\leq 1$, $R_t: [0,1] \times \{1,\cdots,l\} \to  \BR^2$ is an embedding,  $R_0 = \CI$, and there exists $\sigma\in S_l$ such that $R_1(t) \mid_{I_i} = I_{\sigma(i)}(t)$ or $R_1(t)|_{I_i} = I_{\sigma(i)}(1-t)$. The usual product of paths defines a group structure on the set of ribbon braids up to homotopy among ribbon braids. This group, denoted by $RB_l$, does not depend on the
choice of the marked bands and it is called the ribbon braid group  with $l$ bands. A ribbon braid is \emph{pure} if $\sigma$ is trivial and we define $PRB_l$ to be the \emph{pure ribbon braid group} with $l$ bands. If we further assume $R_1(t) \mid_{I_i} = I_{\sigma(i)}(t)$, this subgroup is called the \emph{oriented ribbon braid group} $RB_l^+$. Similarly, we have the \emph{oriented pure ribbon braid group} $PRB_l^+$. 
\end{defn}

\begin{rem}
Note that $RB_l \cong \BZ^l\rtimes B_l$, where the action of $B_l$ is induced by the symmetric group action on the coordinates of $\BZ^l$. In particular, for the pure ribbon braid group $PRB_l$, we have $PRB_l  \cong \BZ^l\times PB_l$. Under this isomorphism, $RB_l^+ \cong (2\BZ)^l\rtimes B_l$ and $PRB^+_l\cong (2\BZ)^l\times PB_l$.
\end{rem}

\begin{defn}
    A \emph{ribbon braided paired $(d,r)$-forest diagram} is a triple $(F_-,\mathfrak{r},F_+)$
    consisting of two $(d,r)$-forests $F_-$ and $F_+$ both with $l$ leaves for some $l$
     and a ribbon braid~$\mathfrak{r} \in RB_l$ connecting the leaves of $F_-$ to the leaves of $F_+$.
\end{defn}

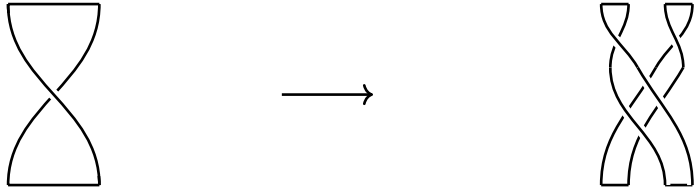
\begin{figure}[h]
\centering
\begin{tikzpicture}[line width=1pt,scale=1.2]
 \draw
   (0,0) -- (1,0)   (0,-2) -- (1,-2) 
  (0,0) to [out=-90, in=90] (1,-2)   (1,0) to [out=-90, in=90] (0,-2);
  
    \draw[white, line width=4pt]
    (0,-0.1) to [out=-90, in=90] (1,-1.9);
    
   \draw   (0,0) to [out=-90, in=90] (1,-2);  
    
  \draw[->]
   (3,-1) -> (4,-1);

 \begin{scope}[xshift=6.5cm]

   \draw (0,0) -- (0.3,0) (0.7,0) -- (1,0)   (0,-2) -- (0.3,-2) (0.7,-2) -- (1,-2);

      
      \draw (0.3,0) to [out=-90, in=90] (0.1,-0.7)
    (1,0) to [out=-90, in=60] (0.6,-.7);
  
    \draw[white, line width=4pt]
 (0,-.1) to [out=-90, in=110] (0.4,-0.7)
 (.7,-.1) to [out=-90, in=80] (0.9,-.7) ;
 
\draw (0,0) to [out=-90, in=120] (0.4,-0.7)  (0.7,0) to [out=-90, in=90]  (0.9,-.7)  ;

 \draw   
 (0.9,-0.7) to [out=-120, in=90] (0.3,-2)
 (0.6,-0.7) to [out=-120, in=90] (0,-2);
    
   \draw[white, line width=4pt]
(0.1,-0.7) to [out=-90, in=90] (0.7,-2)
 (0.4,-0.7) to [out=-60, in=90] (1,-2) ;

     \draw   
 (0.1,-0.7) to [out=-90, in=90] (0.7,-2)
 (0.4,-0.7) to [out=-60, in=90] (1,-2);

 \end{scope}
\end{tikzpicture}
\caption{Splitting a ribbon into $2$ ribbons.}
\label{fig:Splitting_ribb}
\end{figure}

The expansion and reduction rules for the ribbon braids just come from the natural way of splitting a ribbon band into $d$ components and the inverse operation to this. See Figure \ref{fig:Splitting_ribb} for how to split a half twisted band when $d=2$. Note that not only are the two bands themselves twisted but the bands are also braided.  Everything else will be the same as in the braided case, so we omit the details here.  As usual, we define two ribbon braided paired forest diagrams to be equivalent if one is obtained
from the other by a sequence of reductions or expansions.  The
multiplication operation~$\ast$ on the equivalence classes is defined
the same way as for~$bV_{d,r}$.

\begin{defn}
    The \emph{ribbon Higman--Thompson group} $RV_{d,r}$ (resp. the \emph{oriented ribbon Higman--Thompson group} $RV^+_{d,r}$) is the group of equivalence
    classes of (resp. oriented) ribbon braided paired $(d,r)$-forests diagrams with the
    multiplication~$\ast$. The \emph{oriented ribbon Higman--Thompson group} $RF^+_{d,r}$ is the group of equivalence classes of oriented ribbon braided paired $(d,r)$-forest diagrams where the ribbon braids are all pure. Finally, the \emph{oriented ribbon Higman--Thompson group} $RT^+_{d,r}$ is the group of equivalence classes of oriented ribbon braided paired $(d,r)$-forest diagrams where the permutations corresponding to the braids are all cyclic.
\end{defn}

\begin{rem}
As in Definition \ref{defn-bht-lb}, in order for the definition of the ribbon Higman--Thompson groups to work for $F$ and $T$, we need the ribbon braids to stay pure under the expansion and hence the ribbon braids must be oriented. 
\end{rem}

View the braid group $B_d$ as the mapping class group of the disk with $d$ marked points. Let $C=\langle \sigma \rangle $ be the subgroup of $B_d$ generated by the  (counterclockwise) half Dehn twist around the boundary, then the corresponding group $bV_{d,r}(C)$ can be identified with the ribbon Higman--Thompson group $RV_{d,r}$ as follows. The group $B_d(C)$ can be naturally identified with the ribbon braid group $RB_d$ by mapping the label $\sigma^k$ in each string to a band twisted counterclockwise with angle $k\pi$. Moreover, the expansion and reduction and multiplication rule for braided paired $(H,d,r)$-forests diagrams and the ribbon braided paired $(d,r)$-forests diagrams are exactly the same. Hence, we have identified  $RV_{d,r}$ with $ bV_{d,r}(C)$. The argument in fact shows the following.

\begin{prop}\label{lem-idf-rb-br}
$RV_{d,r}\cong bV_{d,r}(C)$, $RV^+_{d,r}\cong bV_{d,r}(2C)$;  $RF^+_{d,r}\cong bF_{d,r}(2C)$;  $RT^+_{d,r}\cong bT_{d,r}(2C)$.
\end{prop}

Thumann showed the following in \cite[Section 4.6.2]{Th17}.
\begin{thm}
The  ribbon Higman--Thompson group $RV_{2,1}$ is of type $F_\infty$.
\end{thm}

\section{Finiteness properties of braided Higman--Thompson groups}\label{sec:finitenessprop}
In this section, we will determine the finiteness properties of the braided Higman--Thompson groups $bV_{d,r}(H)$, $bF_{d,r}(H)$ and $bT_{d,r}(H)$. First, we will generalize the braided paired forest diagrams to allow for the forests to each have an arbitrary number of trees. This will be used to build a complex which the groups act on that will then allow us to induce the finiteness properties of the corresponding braided Higman--Thompson groups. Recall that in Definition \ref{defn-bht-lb}, the label group $H$ for $bV_{d,r}(H)$ is a subgroup of $B_d$, while for $bF_{d,r}(H)$ and $bT_{d,r}(H)$, it lies in $PB_d$, although $H$ will not play a big role in our proof.

By the terminology in \cite{BFM+16}, given a braided paired
forest diagram $(F_-,(b, \lambda),F_+)$ where $F_-,F_+$ are forests with $l$ leaves and  $(b, \lambda) \in B_l(H)$, we call a $d$-caret in~$F_-$ a
\emph{split}.  Similarly a \emph{merge} is a $d$-caret in $F_+$.
With this terminology, the picture representing the
braided paired forest diagram is called a \emph{split-braid-merge diagram}, abbreviated
\emph{spraige}.  We first draw one strand splitting up into $l$
strands in a certain way, representing $F_-$.  Then the $l$ strands braid and are labeled with labels defined by $\lambda$, representing $(b, \lambda)$, and finally according to $F_+$ we merge the
strands back together.

\begin{defn}\label{def:spraiges}
An \emph{$(n,m)$-spraige} is a spraige that begins
on $n$ strands, the \emph{heads}, and ends on $m$ strands, the
\emph{feet}.  As indicated above, we can equivalently think of an
$(n,m)$-spraige as a \emph{braided paired forest diagram} $(F_-, (b, \lambda),
F_+)$, where~$F_-$ has $n$ roots,~$F_+$ has $m$ roots and both have
the same number of leaves.  By an~$n$-\emph{spraige} we mean an
$(n,m)$-spraige for some $m$, and by a \emph{spraige} we mean an~$(n,m)$-spraige
for some $n$ and $m$.  Let $\spraige$ denote the set of all spraiges,~$\spraige_{n,m}$ the set of all $(n,m)$-spraiges, and $\spraige_n$ the set of
all $n$-spraiges.
\end{defn}

Note that an $n$-spraige has $n$ heads, but can have any number of
feet.  This gives a natural function, namely the ``number of feet'' function
$f\colon\spraige\to\BN$ given by $f(\sigma)=m$ if
$\sigma\in\spraige_{n,m}$ for some $n$.

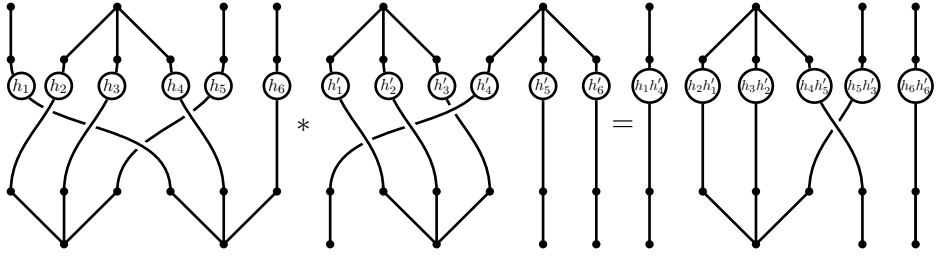
\begin{figure}[h]
\centering
\begin{tikzpicture}[line width=1pt,scale=0.7]
  
  \draw
   (-1,0) -- (-1,-1)   (3,0) -- (3,-1)   (4,0) -- (4,-3.5) (1,0) -- (1,-1)
   (0,-1) -- (1,0) -- (2,-1);
  \draw
   (3,-1) to [out=270, in=90] (1,-3.5);
  \draw[line width=4pt, white]
   (-1,-1) to [out=270, in=90] (2,-3.5);
  \draw
   (-1,-1) to [out=270, in=90] (2,-3.5);
  \draw[line width=4pt, white]
   (1,-1) to [out=270, in=90] (0,-3.5);
  \draw
   (1,-1) to [out=270, in=90] (0,-3.5);
  \draw[line width=4pt, white]
   (0,-1) to [out=270, in=90] (-1,-3.5);
  \draw
   (0,-1) to [out=270, in=90] (-1,-3.5);

  \draw[line width=4pt, white]
   (2,-1) to [out=270, in=90] (3,-3.5);
  \draw
   (2,-1) to [out=270, in=90] (3,-3.5) -- (3,-4.5)
   (2,-3.5) -- (3,-4.5) -- (4,-3.5) (-1,-3.5)--(0,-4.5)
   (0,-3.5) -- (0,-4.5) -- (1,-3.5);
  \node at (4.5,-2.25) {$\ast$};
  \filldraw
   (-1,0) circle (1.5pt)   (-1,-1) circle (1.5pt) (1,0) circle (1.5pt)  (3,0) circle (1.5pt)   (4,0)
circle (1.5pt)   (0,-1) circle (1.5pt)   (1,-1) circle (1.5pt)   (2,-1) circle
(1.5pt)   (3,-1) circle (1.5pt)   (4,-1) circle (1.5pt)   (0,-3.5) circle
(1.5pt)   (1,-3.5) circle (1.5pt)   (2,-3.5) circle (1.5pt)   (3,-3.5) circle
(1.5pt)   (4,-3.5) circle (1.5pt)   (-1,-3.5) circle (1.5pt)  (0,-4.5) circle 
(1.5pt)   (3,-4.5) circle (1.5pt);

   \draw[white,fill=white] (-.8,-1.5) circle (7 pt);
\draw   (-.8,-1.5) circle (7pt) node[text=black, scale=.6] {$h_1$};
  
     \draw[white,fill=white] (-.1,-1.5) circle (7 pt);
\draw   (-.1,-1.5) circle (7pt) node[text=black, scale=.6] {$h_2$};

       \draw[white,fill=white] (.9,-1.5) circle (7 pt);
\draw   (.9,-1.5) circle (7pt) node[text=black, scale=.6] {$h_3$};

     \draw[white,fill=white] (2.1,-1.5) circle (7 pt);
\draw   (2.1,-1.5) circle (7pt) node[text=black, scale=.6] {$h_4$};

     \draw[white,fill=white] (2.9,-1.5) circle (7 pt);
\draw   (2.9,-1.5) circle (7pt) node[text=black, scale=.6] {$h_5$};

     \draw[white,fill=white] (4,-1.5) circle (7 pt);
\draw   (4,-1.5) circle (7pt) node[text=black, scale=.6] {$h_6$};

 \begin{scope}[xshift=5cm]
  \draw
   (0,-1) -- (1,0) -- (2,-1) (1,0)--(1,-1)
   (3,-1) -- (4,0) -- (5,-1) -- (5,-4.5) (4,0)--(4,-1)
   (2,-1) to [out=270, in=90] (3,-3.5);
  \draw[line width=4pt, white]
   (3,-1) to [out=270, in=90] (0,-3.5);
  \draw
   (3,-1) to [out=270, in=90] (0,-3.5);
  \draw[line width=4pt, white]
  (1,-1) to [out=270, in=90] (2,-3.5)
   (0,-1) to [out=270, in=90] (1,-3.5);
  \draw
   (0,-1) to [out=270, in=90] (1,-3.5)
   (1,-1) to [out=270, in=90] (2,-3.5)
   (4,-1)--(4,-4.5)
   (1,-3.5) -- (2,-4.5) -- (3,-3.5) (2,-4.5)--(2,-3.5)   (0,-3.5) -- (0,-4.5);
  \node at (5.5,-2.25) {$=$};
  \filldraw
   (0,-1) circle (1.5pt)   (1,0) circle (1.5pt)    (4,0) circle (1.5pt)  (5,-4.5) circle (1.5pt) (4,-3.5) circle (1.5pt)
   (4,-1) circle (1.5pt) (5,-1) circle (1.5pt) (5,-3.5) circle (1.5pt)
(1,-1) circle (1.5pt)   (2,-1) circle (1.5pt)   (3,-1) circle (1.5pt)   (0,-3.5)
circle (1.5pt)   (1,-3.5) circle (1.5pt)   (2,-3.5) circle (1.5pt)   (3,-3.5)
circle (1.5pt)   (0,-4.5) circle (1.5pt)   (2,-4.5) circle (1.5pt)   (4,-4.5)
circle (1.5pt);

     \draw[white,fill=white] (.1,-1.5) circle (7 pt);
\draw   (.1,-1.5) circle (7pt) node[text=black, scale=.55] {$h_1'$};

       \draw[white,fill=white] (1.1,-1.5) circle (7 pt);
\draw   (1.1,-1.5) circle (7pt) node[text=black, scale=.55] {$h_2'$};

     \draw[white,fill=white] (2.1,-1.5) circle (7 pt);
\draw   (2.1,-1.5) circle (7pt) node[text=black, scale=.55] {$h_3'$};

     \draw[white,fill=white] (2.9,-1.5) circle (7 pt);
\draw   (2.9,-1.5) circle (7pt) node[text=black, scale=.55] {$h_4'$};

     \draw[white,fill=white] (4,-1.5) circle (7 pt);
\draw   (4,-1.5) circle (7pt) node[text=black, scale=.55] {$h_5'$};
  
       \draw[white,fill=white] (5,-1.5) circle (7 pt);
\draw   (5,-1.5) circle (7pt) node[text=black, scale=.55] {$h_6'$};

 \end{scope}

 \begin{scope}[xshift=11cm]
  \draw
   (0,0) -- (0,-4.5)   (2,0) -- (2,-3.5)   (4,0) -- (4,-1) (5,0) -- (5,-4.5)
   (1,-3.5) -- (1,-1) -- (2,0) -- (3,-1) 
   (4,-1) to [out=270, in=90] (3,-3.5);
  \draw[line width=4pt, white]
   (3,-1) to [out=270, in=90] (4,-3.5);
  \draw
   (3,-1) to [out=270, in=90] (4,-3.5)
   (1,-3.5) -- (2,-4.5) -- (3,-3.5) (2,-3.5)--(2,-4.5)
   (4,-3.5) -- (4,-4.5)  (5,-4.5) -- (5,-3.5);
  \filldraw
   (0,0) circle (1.5pt)   (5,0) circle (1.5pt)   (5,-1) circle (1.5pt)  (2,0) circle (1.5pt)  
(4,0) circle (1.5pt)   (0,-1) circle (1.5pt)   (1,-1) circle (1.5pt)   (2,-1)
circle (1.5pt)   (3,-1) circle (1.5pt)   (4,-1) circle (1.5pt)   (0,-3.5)
circle (1.5pt)   (1,-3.5) circle (1.5pt)   (2,-3.5) circle (1.5pt)   (3,-3.5)
circle (1.5pt)   (4,-3.5) circle (1.5pt)   (0,-4.5) circle (1.5pt)  (5,-4.5) circle (1.5pt)   (5,-3.5) circle (1.5pt)
(2,-4.5) circle (1.5pt)   (4,-4.5) circle (1.5pt);

     \draw[white,fill=white] (0,-1.5) circle (9 pt);
\draw   (0,-1.5) circle (9pt) node[text=black, scale=.5] {$h_1h_4'$};

       \draw[white,fill=white] (1,-1.5) circle (9 pt);
\draw   (1,-1.5) circle (9pt) node[text=black, scale=.5] {$h_2h_1'$};

     \draw[white,fill=white] (2,-1.5) circle (9 pt);
\draw   (2,-1.5) circle (9pt) node[text=black, scale=.5] {$h_3h_2'$};

     \draw[white,fill=white] (3.1,-1.5) circle (9 pt);
\draw   (3.1,-1.5) circle (9pt) node[text=black, scale=.5] {$h_4h_5'$};

     \draw[white,fill=white] (4,-1.5) circle (9 pt);
\draw   (4,-1.5) circle (9pt) node[text=black, scale=.5] {$h_5h_3'$};
  
       \draw[white,fill=white] (5,-1.5) circle (9 pt);
\draw   (5,-1.5) circle (9pt) node[text=black, scale=.5] {$h_6h_6'$};

 \end{scope}
\end{tikzpicture}
\caption{Multiplication of spraiges.}
\label{fig:spraiges-multiplication}
\end{figure}

The pictures in Figure~\ref{fig:spraiges-multiplication} are examples
of spraiges. It is clear that the notion of reduction and expansion generalizes to 
diagrams of arbitrary spraiges, and one can consider equivalence classes under reduction
and expansion.  As is the case with paired forest diagrams and braided paired forest diagrams, each such class has a unique reduced representative.  We will just
call an equivalence class of spraiges a spraige, so in particular the elements
of $bV_{d,r}(H), bF_{d,r}(H)$, and $bT_{d,r}(H)$ are all sets of $(r,r)$-spraiges.

The operation $\ast$ defined for braided Higman--Thompson groups can be defined in general
for spraiges via concatenation of diagrams.  It is only defined for
certain pairs of spraiges, namely we can multiply
$\sigma_1\ast\sigma_2$ for $\sigma_1\in\spraige_{n_1,m_1}$ and
$\sigma_2 \in \spraige_{n_2,m_2}$ if and only if $m_1=n_2$.  In this
case, we obtain $\sigma_1 \ast \sigma_2 \in \spraige_{n_1,m_2}$. In the figures, we will sometimes lengthen a single-node tree to an edge for aesthetic reasons.

Note that for every $n\in \BN$ there is an identity
$(n,n)$-spraige $1_n$ with respect to $\ast$, namely the spraige
represented by $(1_n,(\id, \iota),1_n)$ where $\iota$ is the trivial function which chooses the identity in $H$ as the label for each strand.  By abuse of notation, we are using $1_n$ to also denote the trivial forest with $n$ roots. Note also, given any $(n,m)$-spraige $(F_-,(b, \lambda),F_+)$, there exists an
inverse $(m,n)$-spraige $(F_+,(b, \lambda)^{-1},F_-)$ with
\begin{align*}
(F_-,(b, \lambda),F_+) \ast (F_+,(b, \lambda)^{-1},F_-) & = 1_n \\
\intertext{and}
(F_+,(b, \lambda)^{-1},F_-)\ast(F_-,(b, \lambda),F_+) & = 1_m 
\end{align*}
These two together give that $\spraige$ is a groupoid under the operation $\ast$.

\begin{figure}[h]
\centering
\begin{tikzpicture}
  \draw[line width=1pt]
   (0.5,-0.5) -- (1,0) -- (1.5,-0.5) (1,0)--(1,-.5) (4,0)--(4,-.5)
   (3.5,-0.5) -- (4,0) -- (4.5,-0.5);
  \filldraw 
   (0,0) circle (1.5pt)
   (1,0) circle (1.5pt)
   (2,0) circle (1.5pt)
   (3,0) circle (1.5pt)
   (4,0) circle (1.5pt)

   (0.5  ,-0.5) circle (1.5pt)
      (1  ,-0.5) circle (1.5pt)
      (4  ,-0.5) circle (1.5pt)
   (1.5  ,-0.5) circle (1.5pt)
   (3.5  ,-0.5) circle (1.5pt)
   (4.5  ,-0.5) circle (1.5pt);
	
  \node at (-0.7,0) {$F_{\{2,5\}}^{(5)}$};

 \begin{scope}[yshift=-1.3cm]
  \draw[line width=1pt]
   (0.5,-1) -- (0.5,-0.5) -- (1,0) -- (1.5,-0.5) -- (1.5,-1) (1,0)--(1,-1) (4,0)--(4,-1)
   (3.5,-1) -- (3.5,-0.5) -- (4,0) -- (4.5,-0.5) -- (4.5,-1)
   (0,0) -- (0,-1)   (2,0) -- (2,-1)   (3,0) -- (3,-1);
  \filldraw 
   (0,0) circle (1.5pt)   (0,-0.5) circle (1.5pt)   (0,-1) circle (1.5pt)
   (1,0) circle (1.5pt)
   (2,0) circle (1.5pt)   (2,-0.5) circle (1.5pt)   (2,-1) circle (1.5pt)
   (3,0) circle (1.5pt)   (3,-0.5) circle (1.5pt)   (3,-1) circle (1.5pt)
   (4,0) circle (1.5pt)

     (1  ,-0.5) circle (1.5pt)
      (4  ,-0.5) circle (1.5pt)
           (1  ,-1) circle (1.5pt)
      (4  ,-1) circle (1.5pt)

   (0.5,-0.5) circle (1.5pt)   (0.5,-1) circle (1.5pt)
   (1.5,-0.5) circle (1.5pt)   (1.5,-1) circle (1.5pt)
   (3.5,-0.5) circle (1.5pt)   (3.5,-1) circle (1.5pt)
   (4.5,-0.5) circle (1.5pt)   (4.5,-1) circle (1.5pt);
	
  \node at (-0.7,-0.5) {$\lambda_{\{2,5\}}^{(5)}$};
 \end{scope}

 \begin{scope}[yshift=-4cm, yscale=-1]
  \draw[line width=1pt]
   (0.5,-1) -- (0.5,-0.5) -- (1,0) -- (1.5,-0.5) -- (1.5,-1) (1,0)--(1,-1) (4,0)--(4,-1)
   (3.5,-1) -- (3.5,-0.5) -- (4,0) -- (4.5,-0.5) -- (4.5,-1)
   (0,0) -- (0,-1)   (2,0) -- (2,-1)   (3,0) -- (3,-1);
  \filldraw 
   (0,0) circle (1.5pt)   (0,-0.5) circle (1.5pt)   (0,-1) circle (1.5pt)
   (1,0) circle (1.5pt)
   (2,0) circle (1.5pt)   (2,-0.5) circle (1.5pt)   (2,-1) circle (1.5pt)
   (3,0) circle (1.5pt)   (3,-0.5) circle (1.5pt)   (3,-1) circle (1.5pt)
   (4,0) circle (1.5pt)

     (1  ,-0.5) circle (1.5pt)
      (4  ,-0.5) circle (1.5pt)
           (1  ,-1) circle (1.5pt)
      (4  ,-1) circle (1.5pt)

   (0.5,-0.5) circle (1.5pt)   (0.5,-1) circle (1.5pt)
   (1.5,-0.5) circle (1.5pt)   (1.5,-1) circle (1.5pt)
   (3.5,-0.5) circle (1.5pt)   (3.5,-1) circle (1.5pt)
   (4.5,-0.5) circle (1.5pt)   (4.5,-1) circle (1.5pt);
	
  \node at (-0.7,-0.5) {$\mu_{\{2,5\}}^{(5)}$};
 \end{scope}
\end{tikzpicture}
\caption{The elementary forest $F^{(5)}_{\{2,5\}}$, and the spraiges
$\lambda^{(5)}_{\{2,5\}}$ and $\mu^{(5)}_{\{2,5\}}$.}
\label{fig:elem_forests}
\end{figure}
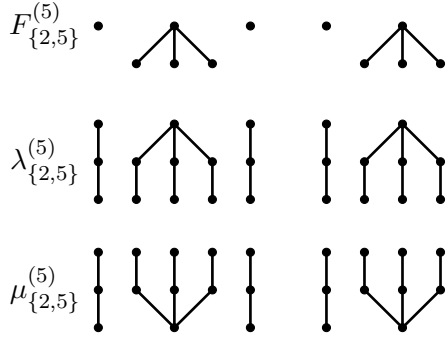

Some forests will be important enough to the construction of the Stein space in Section \ref{sec:def_stein_space} that we name them now.
 For $n\in \BN$ and $J \subseteq \{1,\dots, n\}$, define $F^{(n)}_J$
to be the forest with $n$ roots and $|J|$ carets, with a caret
attached to the $i^{\text{th}}$ root for each $i\in J$.  A characterizing property of these forests is that every caret is elementary and so we will call such forests \emph{elementary}. Define the spraige
$\lambda^{(n)}_J$ to be the $(n,n+(d-1)|J|)$-spraige $(F^{(n)}_J, (\id, \iota),
1_{n+(d-1)|J|})$, and the spraige $\mu^{(n)}_J$ to be its inverse.  If $J =
\{i\}$, we will write $F^{(n)}_i$, $\lambda^{(n)}_i$ and $\mu^{(n)}_i$ instead.
See Figure~\ref{fig:elem_forests} for an example of an elementary
forest and the corresponding spraiges. Note that we did not draw the labels as they are all $1\in H$.

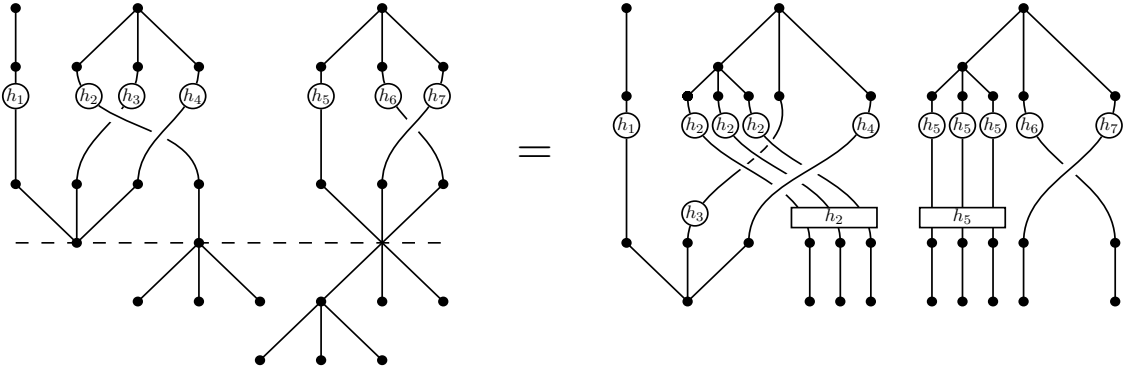
\begin{figure}[h]
\centering

\resizebox{15cm}{4.8cm}{%
\begin{tikzpicture}[line width=.4pt]
  
  \draw
(-0.5,0)--(-0.5,-.5)   (0,-0.5) -- (0.5,0) -- (1,-0.5)   (0.5,0) -- (0.5,-0.5)    (2.5,0) -- (3,-.5)  (2.5,0) -- (2,-.5)--
(2,-1.5)   (2.5,0) -- (2.5,-0.5)    
   (0.5,-0.5) to [out=-90, in=90] (0,-1.5)   (2.5,-0.5) to [out=-90, in=90]
(3,-1.5);

  \draw[white, line width=4pt]
   (0,-0.5) to [out=-90, in=90] (1,-1.5)   (3,-0.5) to [out=-90, in=90]
(2.5,-1.5);
  \draw
   (0,-0.5) to [out=-90, in=90] (1,-1.5)   (3,-0.5) to [out=-90, in=90]
(2.5,-1.5);
  \draw[white, line width=4pt]
   (1,-0.5) to [out=-90, in=90] (0.5,-1.5);
  \draw
   (1,-0.5) to [out=-90, in=90] (0.5,-1.5);

  \draw
   (0,-1.5) -- (0,-2) -- (0.5,-1.5)   (1,-1.5) -- (1,-2)   (2,-1.5) --
(2.5,-2) -- (3,-1.5)   (2.5,-2) -- (2.5,-2.5)    (2.5,-1.5) -- (2.5,-2);

  \draw[dashed]
   (-0.5,-2) -- (3,-2);

  \draw
   (0.5,-2.5) -- (1,-2) -- (1.5,-2.5)  (1,-2)--(1,-2.5)  (2,-2.5) -- (2.5,-2) -- (3,-2.5)  (-.5,-0.5)-- (-.5,-1.5)--(0,-2)
   (2,-2.5) -- (2.5,-3)
   (1.5,-3) -- (2,-2.5) -- (2,-3);

  \node at (3.75,-1.25) {$=$};

  \filldraw
  
    (-.5,0) circle (1pt)   (-.5,-0.5) circle (1pt)  (-.5,-1.5) circle (1pt)
  (3,-2.5) circle (1pt) (2.5,-2.5) circle (1pt)  (2,-2.5) circle (1pt) 
  (2.5,-3) circle (1pt) 
   (0,-0.5) circle (1pt)   (0.5,-0.5) circle (1pt)   (1,-0.5) circle (1pt)  
   (2.5,-0.5) circle (1pt)  (2.5,0)circle (1pt)   (2,-0.5) circle (1pt)

     (0.5,0) circle (1pt)  
(3,-0.5) circle (1pt)   
   (0,-1.5) circle (1pt)   (0.5,-1.5) circle (1pt)   (1,-1.5) circle (1pt)  
 (2,-1.5) circle (1pt)   (2.5,-1.5) circle (1pt)  
 (3,-1.5) circle (1pt) 
   (0,-2) circle (1pt)   (1,-2) circle (1pt)    

   (0.5,-2.5) circle (1pt)   (1,-2.5) circle (1pt)   (1.5,-3) circle (1pt)  
(1.5,-2.5) circle (1pt)   (2,-3) circle (1pt)     ;

     \draw[white,fill=white] (-.5,-.75) circle (3 pt);
\draw  (-.5,-.75) circle (3 pt) node[text=black, scale=.4] {$h_1$};

     \draw[white,fill=white] (.1,-.75) circle (3 pt);
\draw  (.1,-.75) circle (3 pt) node[text=black, scale=.4] {$h_2$};

     \draw[white,fill=white] (.45,-.75) circle (3 pt);
\draw  (.45,-.75) circle (3 pt) node[text=black, scale=.4] {$h_3$};

     \draw[white,fill=white] (.95,-.75) circle (3 pt);
\draw  (.95,-.75) circle (3 pt) node[text=black, scale=.4] {$h_4$};

     \draw[white,fill=white] (2,-.75) circle (3 pt);
\draw  (2,-.75) circle (3 pt) node[text=black, scale=.4] {$h_5$};
  
     \draw[white,fill=white] (2.55,-.75) circle (3 pt);
\draw  (2.55,-.75) circle (3 pt) node[text=black, scale=.4] {$h_6$};

     \draw[white,fill=white] (2.95,-.75) circle (3 pt);
\draw  (2.95,-.75) circle (3 pt) node[text=black, scale=.4] {$h_7$};

 \begin{scope}[xshift=5cm]
  \draw
  (-0.5,0)--(-0.5,-2)

   (0,-0.75) -- (0.75,0) -- (1.5,-0.75)   (0.5,-0.75) -- (0.25,-0.5)--(0.25,-0.75) (0.75,0)--(.75,-.75)   
   (2,-.75) -- (2.75,0)   (2,-.75)--(2,-2)  (2.25,-.75)--(2.25,-2)  (2.25,-.75) -- (2.25,-.5)  -- (2.5,-0.75)  (2.5,-.75)--(2.5,-2) 
 (3.5,-0.75) -- (2.75,0) (2.75,0)--(2.75,-.75)
   (.75,-0.75) to [out=-70, in=90] (0,-2)   (2.75,-0.75) to [out=-90, in=90]
(3.5,-2);

  \draw[white, line width=2pt]
   (0,-0.75) to [out=-90, in=90] (1,-2)   (0.25,-0.75) to [out=-90, in=90] (1.25,-2)  (0.5,-0.75) to [out=-90, in=90]
(1.5,-2);

  \draw[white, line width=4pt]
  (3.5,-0.75) to [out=-90, in=90] (2.75,-2);
  \draw
   (0,-0.75) to [out=-90, in=90] (1,-2)  (0.25,-0.75) to [out=-90, in=90] (1.25,-2)   (0.5,-0.75) to [out=-90, in=90]
(1.5,-2)   (3.5,-0.75) to [out=-90, in=90] (2.75,-2);
  \draw[white, line width=4pt]
   (1.5,-0.75) to [out=-90, in=90] (0.5,-2);
  \draw
   (1.5,-0.75) to [out=-90, in=90] (0.5,-2);

  \draw
   (0,-2) -- (0,-2.5) -- (0.5,-2) (-.5,-2)--(0,-2.5)
   (1.25,-2)--(1.25,-2.5)   (2.25,-2)--(2.25,-2.5)
   (1,-2) -- (1,-2.5)   (1.5,-2) --
(1.5,-2.5)   (2,-2) -- (2,-2.5)   (2.5,-2) -- (2.5,-2.5)   (2.75,-2) -- (2.75,-2.5)  
(3.5,-2) -- (3.5,-2.5)   ;

  \filldraw
  
    (-0.5,0)circle (1pt)  (-0.5,-0.75)circle (1pt) (-0.5,-2)circle (1pt)
  
   (0,-0.75) circle (1pt)  (0,-0.75) circle (1pt)  (0,-0.75) circle (1pt)
  
   (1.25,-2) circle (1pt)   (1.25,-2.5) circle (1pt)
    (2.25,-.75)circle (1pt)  (2.75,-.75) circle (1pt)
   (2.25,-2)circle (1pt)  (2.25,-2.5)circle (1pt) 
   (0,-0.75) circle (1pt)   (0.5,-0.75) circle (1pt)   (.75,-0.75) circle (1pt)  
(1.5,-0.75) circle (1pt)   (2,-0.75) circle (1pt)   (2.5,-0.75) circle (1pt)  
  (3.5,-0.75) circle (1pt)   
   (0.25,-0.5) circle (1pt)   (0.25,-0.75) circle (1pt)   (0.75,0) circle (1pt)  
(2.75,0)  circle (1pt)   (2.25,-0.5) circle (1pt)   

   (0,-2) circle (1pt)   (0.5,-2) circle (1pt)   (1,-2) circle (1pt)   (1.5,-2)
circle (1pt)   (2,-2) circle (1pt)   (2.5,-2) circle (1pt)   (2.75,-2) circle (1pt)
  (3.5,-2) circle (1pt)   
   (0,-2.5) circle (1pt)   (1,-2.5) circle (1pt)   (1.5,-2.5) circle (1pt)  
(2,-2.5) circle (1pt)   (2.5,-2.5) circle (1pt)   (2.75,-2.5) circle (1pt)  
(3.5,-2.5) circle (1pt)   ;

     \draw[white,fill=white] (-.5,-1) circle (3 pt);
\draw  (-.5,-1) circle (3 pt) node[text=black, scale=.4] {$h_1$};

     \draw[white,fill=white] (0.06,-1) circle (3 pt);
\draw  (0.06,-1) circle (3 pt) node[text=black, scale=.4] {$h_2$};

     \draw[white,fill=white] (.31,-1) circle (3 pt);
\draw  (.31,-1) circle (3 pt) node[text=black, scale=.4] {$h_2$};

     \draw[white,fill=white] (.56,-1) circle (3 pt);
\draw  (.56,-1) circle (3 pt) node[text=black, scale=.4] {$h_2$};

     \draw[white,fill=white] (.06,-1.75) circle (3 pt);
\draw  (.06,-1.75) circle (3 pt) node[text=black, scale=.4] {$h_3$};

     \draw[white,fill=white] (1.46,-1) circle (3 pt);
\draw  (1.46,-1) circle (3 pt) node[text=black, scale=.4] {$h_4$};

     \draw[white,fill=white] (2,-1) circle (3 pt);
\draw  (2,-1) circle (3 pt) node[text=black, scale=.4] {$h_5$};

     \draw[white,fill=white] (2.25,-1) circle (3 pt);
\draw  (2.25,-1) circle (3 pt) node[text=black, scale=.4] {$h_5$};

     \draw[white,fill=white] (2.5,-1) circle (3 pt);
\draw  (2.5,-1) circle (3 pt) node[text=black, scale=.4] {$h_5$};
  
     \draw[white,fill=white] (2.8,-1) circle (3 pt);
\draw  (2.8,-1) circle (3 pt) node[text=black, scale=.4] {$h_6$};

     \draw[white,fill=white] (3.45,-1) circle (3 pt);
\draw  (3.45,-1) circle (3 pt) node[text=black, scale=.4] {$h_7$};

    \filldraw[fill=white, draw=black] (0.85,-1.7) rectangle (1.55,-1.87);
\draw  (1.2,-1.77)  node[text=black, scale=.4] {$h_2$};   

    \filldraw[fill=white, draw=black] (1.9,-1.7) rectangle (2.6,-1.87);
\draw  (2.25,-1.77)  node[text=black, scale=.4] {$h_5$};  

 \end{scope}
\end{tikzpicture}

}
\caption{A splitting of a spraige.}
\label{fig:spraige_splitting}
\end{figure}

\begin{figure}[h]
\centering

\resizebox{12cm}{4.1cm}{%

\begin{tikzpicture}[line width=.4pt]
  
  \draw
 (-0.5,-2)--(-.5,0)  (0,-0.5) -- (0.5,0) -- (1,-0.5)   (0.5,0) -- (0.5,-0.5)   (1.5,-0.5) --
(2,0) -- (2.5,-0.5)   (2,0) -- (2,-0.5)
   (0.5,-0.5) to [out=-90, in=90] (0,-1.5)   (1.5,-0.5) to [out=-90, in=90]
(2,-1.5)   (2,-0.5) to [out=-90, in=90] (2.5,-1.5) 
 (2.5,-0.5) to [out=-90, in=90] (3,-1.5);

  \draw[white, line width=4pt]
   (0,-0.5) to [out=-90, in=90] (1,-1.5)   (3,-0.5) to [out=-90, in=90]
(1.5,-1.5);
  \draw
   (0,-0.5) to [out=-90, in=90] (1,-1.5)   (3,-0.5) to [out=-90, in=90]
(1.5,-1.5);
  \draw[white, line width=4pt]
   (1,-0.5) to [out=-90, in=90] (0.5,-1.5);
  \draw
   (1,-0.5) to [out=-90, in=90] (0.5,-1.5);

  \draw
  (3,-1.5)--(3,-2)
   (0,-1.5) -- (0,-2)   (0.5,-1.5) -- (1,-2) -- (1.5,-1.5)   (1,-1.5) --
(1,-2)   (2,-1.5) -- (2,-2)   (2.5,-1.5) -- (2.5,-2);

  \draw[dashed]
   (-0.5,-2) -- (3.2,-2);

  \draw
    (-.5,-2)--(0,-2.5)  
   (0,-2) -- (0,-2.5) -- (1,-2) (3,0) --  (3,-0.5) 
   (2.5,-2.5)-- (2.5,-2)
   (2,-2) -- (2.5,-2.5) -- (3,-2);

  \node at (3.5,-1.25) {$=$};

  \filldraw
  (-0.5,-1.5) circle (1pt) (-0.5,-2) circle (1pt)
  (-0.5,-.5)circle (1pt) (-.5,0)circle (1pt)
   (0,-0.5) circle (1pt)   (0.5,0) circle (1pt)   (0.5,-0.5) circle (1pt)  
  (1,-0.5) circle (1pt)   (1.5,-0.5) circle (1pt)   (2,0)
circle (1pt)   (2,-0.5) circle (1pt)     (2.5,-0.5) circle (1pt)
(3,0) circle (1pt)  (3,-0.5) circle (1pt) 
   (0,-1.5) circle (1pt)   (0,-2) circle (1pt)   (0.5,-1.5) circle (1pt)  
(1,-2) circle (1pt)   (1.5,-1.5) circle (1pt)   (1,-1.5) circle (1pt)  
   (2,-1.5) circle (1pt)   (2,-2) circle (1pt) (3,-2)circle (1pt) 
(2.5,-1.5) circle (1pt)   (2.5,-2) circle (1pt) (2.5,-2.5) circle (1pt)
   (0,-2.5) circle (1pt) (3,-1.5)circle (1pt);

  \draw[white,fill=white] (-.5,-.7) circle (3 pt);
\draw  (-.5,-.7) circle (3 pt) node[text=black, scale=.4] {$h_1$};

     \draw[white,fill=white] (0.03,-.7) circle (3 pt);
\draw  (0.03,-.7) circle (3 pt) node[text=black, scale=.4] {$h_2$};

     \draw[white,fill=white] (.46,-.7) circle (3 pt);
\draw  (.46,-.7) circle (3 pt) node[text=black, scale=.4] {$h_3$};

     \draw[white,fill=white] (.96,-.7) circle (3 pt);
\draw  (.95,-.7) circle (3 pt) node[text=black, scale=.4] {$h_4$};

     \draw[white,fill=white] (1.55,-.7) circle (3 pt);
\draw  (1.55,-.7) circle (3 pt) node[text=black, scale=.4] {$h_5$};

     \draw[white,fill=white] (2.05,-.7) circle (3 pt);
\draw  (2.05,-.7) circle (3 pt) node[text=black, scale=.4] {$h_5$};

     \draw[white,fill=white] (2.55,-.7) circle (3 pt);
\draw  (2.56,-.7) circle (3 pt) node[text=black, scale=.4] {$h_5$};   
 
      \draw[white,fill=white] (2.93,-.7) circle (3 pt);
\draw  (2.93,-.7) circle (3 pt) node[text=black, scale=.4] {$h_6$};

    \filldraw[fill=white, draw=black] (1.84,-1.2) rectangle (3.1,-1.4);
\draw  (2.5,-1.3)  node[text=black, scale=.4] {$h_5$};

 \begin{scope}[xshift=4.5cm]
  \draw
 (-0.5,-1.5)--(-.5,0)  (0,-0.5) -- (0.5,0) -- (1,-0.5)   (0.5,0) -- (0.5,-0.5)   (1.5,-0.5) --
(1.5,0)  (2,0) -- (2,-0.5)
   (0.5,-0.5) to [out=-90, in=90] (0,-1.5)   (1.5,-0.5) to [out=-90, in=90]
(2,-1.5);

  \draw[white, line width=4pt]
   (0,-0.5) to [out=-90, in=90] (1,-1.5)   (2,-0.5) to [out=-90, in=90]
(1.5,-1.5);
  \draw
   (0,-0.5) to [out=-90, in=90] (1,-1.5)   (2,-0.5) to [out=-90, in=90]
(1.5,-1.5);
  \draw[white, line width=4pt]
   (1,-0.5) to [out=-90, in=90] (0.5,-1.5);
  \draw
   (1,-0.5) to [out=-90, in=90] (0.5,-1.5);

  \draw
  
   (0,-1.5) -- (0,-2)   (0.5,-1.5) -- (1,-2) -- (1.5,-1.5)   (1,-1.5) --
(1,-2)   (2,-1.5) -- (2,-2)   ;


  \draw
    (-.5,-1.5)--(0,-2.5)  
   (0,-2) -- (0,-2.5) -- (1,-2)  
   (2,-2) -- (2,-2.5) ;


  \filldraw
  (-0.5,-1.5) circle (1pt)

  (-0.5,-.5)circle (1pt) (-.5,0)circle (1pt)
   (0,-0.5) circle (1pt)   (0.5,0) circle (1pt)   (0.5,-0.5) circle (1pt)  
  (1,-0.5) circle (1pt)   (1.5,-0.5) circle (1pt)   (2,0)
circle (1pt)   (2,-0.5) circle (1pt)     

   (0,-1.5) circle (1pt)     (0.5,-1.5) circle (1pt)  
(1,-2) circle (1pt)   (1.5,-1.5) circle (1pt)   (1,-1.5) circle (1pt)  
   (2,-1.5) circle (1pt)   
   (0,-2.5) circle (1pt)
    (1.5,0) circle (1pt)
   (2,-2.5) circle (1pt);

  \draw[white,fill=white] (-.5,-.7) circle (3 pt);
\draw  (-.5,-.7) circle (3 pt) node[text=black, scale=.4] {$h_1$};

     \draw[white,fill=white] (0.03,-.7) circle (3 pt);
\draw  (0.03,-.7) circle (3 pt) node[text=black, scale=.4] {$h_2$};

     \draw[white,fill=white] (.46,-.7) circle (3 pt);
\draw  (.46,-.7) circle (3 pt) node[text=black, scale=.4] {$h_3$};

     \draw[white,fill=white] (.96,-.7) circle (3 pt);
\draw  (.95,-.7) circle (3 pt) node[text=black, scale=.4] {$h_4$};

     \draw[white,fill=white] (1.55,-.7) circle (3 pt);
\draw  (1.55,-.7) circle (3 pt) node[text=black, scale=.4] {$h_5$};

     \draw[white,fill=white] (1.95,-.7) circle (3 pt);
\draw  (1.95,-.7) circle (3 pt) node[text=black, scale=.4] {$h_6$};

 \end{scope}
\end{tikzpicture}

}
\caption{An elementary merging of a spraige.}
\label{fig:spraige_merging}
\end{figure}
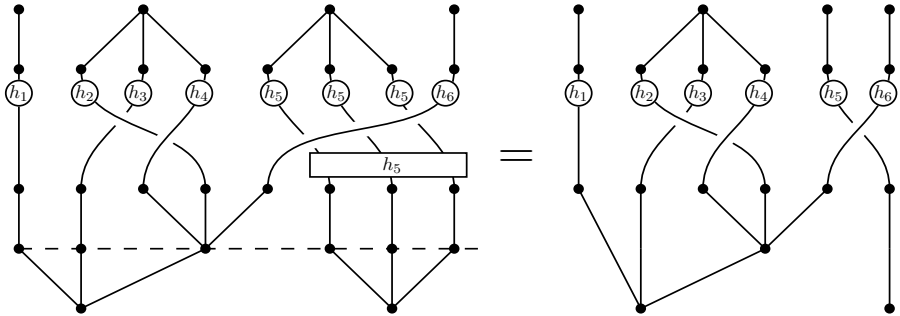

Fix an $(n,m)$-spraige $\sigma$.  For any forest~$F$ with~$m$ roots
and~$l$ leaves, define the \emph{splitting of $\sigma$ by $F$} as
multiplying $\sigma$ by the spraige $(F,(\id, \iota),1_{l})$ from the right.
Similarly, a \emph{merging of $\sigma$ by $F'$} is right multiplication
by the spraige $(1_m,(\id, \iota),F')$, where $F'$ now has $l$ roots and $m$
leaves.  In the case where $F$ (resp.~$F'$) is an elementary
forest, we call this operation \emph{elementary splitting}
(resp.~\emph{elementary merging)}.  See
Figure~\ref{fig:spraige_splitting} for an idea of splitting and
Figure~\ref{fig:spraige_merging} for an idea of elementary merging.

In the special case that $F=F^{(n)}_i$ for $i\in \{1,\dots,n\}$, we can
think of a splitting by~$F$ as simply attaching a single caret to the
$i^{\text{th}}$ foot of a spraige, possibly followed by reductions such as a splitting move as seen in Figure~\ref{fig:reduction_moves} .
Similarly, a merging by $F$ in this case can be thought of as merging
the $i$-th through $(i+d-1)$-th feet together.  In
these cases we will also speak of \emph{adding a split}
(respectively~\emph{merge}) to the spraige.

The following types of spraiges will prove to be particularly
important.  First, a \emph{braige} is defined to be a spraige where
there are no splits, i.e., a spraige of the form~$(1_n,(b, \lambda),F)$ for $b\in
B_n$ and $F$ having $n$ leaves.  Also, when $F$ is elementary, we will
call $(1_n,(b, \lambda),F)$ an \emph{elementary braige}.  Analogously to
spraiges, we define \emph{$n$-braiges} and \emph{elementary
$n$-braiges}.

To deal with $bF_{d,r}(H)$ and $bT_{d.r}(H)$, we make the following convention:  Whenever we
want to only consider pure or cyclic labeled braids, we will attach the modifier
``pure'' or ``cyclic'', e.g., we can talk about pure $n$-spraiges or elementary
cyclic $n$-braiges.

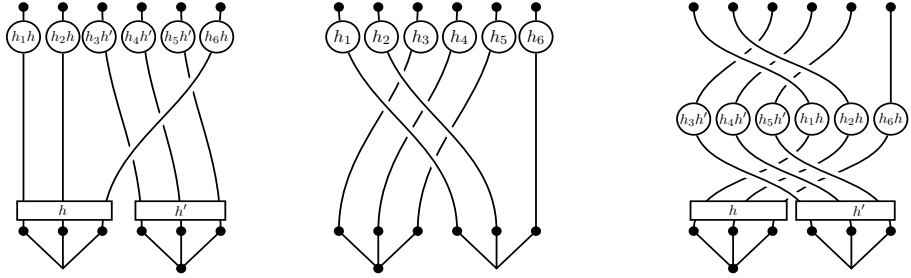
\begin{figure}[t]
\centering

\resizebox{12cm}{3.6cm}{%
\begin{tikzpicture}[line width=.65pt]
  
  \draw
  (0,-1.5) to [out=270, in=90] (-1,-4.5)
   (.5,-1.5) to [out=270, in=90] (-.5,-4.5)
  (1,-1.5) to [out=270, in=90] (0,-4.5)
  (1.5,-1.5) -- (1.5,-4.5);
  \draw [line width=4pt, white]
  (-.5,-1.5) to [out=270, in=90] (1,-4.5)
  
  (-1,-1.5) to [out=270, in=90] (.5,-4.5);
  \draw
  (-1,-1.5) to [out=270, in=90] (.5,-4.5)
  (-.5,-1.5) to [out=270, in=90] (1,-4.5)
  (-1,-4.5) -- (-0.5,-5) -- (0,-4.5)  
  (-.5,-4.5) -- (-0.5,-5)
  (.5,-4.5) -- (1,-5) -- (1.5,-4.5) (1,-5)-- (1,-4.5);
  \filldraw
  (.5,-1.5)circle (1.5pt)   (.5,-4.5)circle (1.5pt) 
  (-1,-1.5) circle (1.5pt)   (0,-1.5) circle (1.5pt)   (1,-1.5) circle (1.5pt)  
(1.5,-1.5) circle (1.5pt) (-.5,-1.5)circle (1.5pt) 
  (-1,-4.5) circle (1.5pt)   (0,-4.5) circle (1.5pt)   (1,-4.5) circle (1.5pt)  
 (1.5,-4.5)circle (1.5pt)

(-.5,-4.5) circle (1.5pt)
  (-0.5,-5) circle (1.5pt);

    \draw[white,fill=white] (-.95,-1.9) circle (6 pt);
\draw  (-.95,-1.9) circle (6 pt) node[text=black, scale=.6] {$h_1$};

    \draw[white,fill=white] (-.46,-1.9) circle (6 pt);
\draw  (-.46,-1.9) circle (6 pt) node[text=black, scale=.6] {$h_2$};

    \draw[white,fill=white] (.04,-1.9) circle (6 pt);
\draw  (.04,-1.9) circle (6 pt) node[text=black, scale=.6] {$h_3$};

    \draw[white,fill=white] (.53,-1.9) circle (6 pt);
\draw  (.53,-1.9) circle (6 pt) node[text=black, scale=.6] {$h_4$};

    \draw[white,fill=white] (1.03,-1.9) circle (6 pt);
\draw  (1.03,-1.9) circle (6 pt) node[text=black, scale=.6] {$h_5$};

    \draw[white,fill=white] (1.5,-1.9) circle (6 pt);
\draw  (1.5,-1.9) circle (6 pt) node[text=black, scale=.6] {$h_6$};

 \begin{scope}[xshift=-3.5cm,xscale=-1]
   \draw

  (1,-1.5)-- (1,-4.5)
  (1.5,-1.5) -- (1.5,-4.5);

  \draw
  (0,-1.5) to [out=270, in=90] (-.5,-4.5)
  (.5,-1.5) to [out=270, in=90] (0,-4.5)
  (-.5,-1.5)to [out=270, in=90](-1,-4.5)
  (-1,-4.5) -- (-0.5,-5) -- (0,-4.5)  
  (-.5,-4.5) -- (-0.5,-5)
  (.5,-4.5) -- (1,-5) -- (1.5,-4.5) (1,-5)-- (1,-4.5);
  
    \draw [line width=4pt, white]
    (-1,-1.5) to [out=270, in=90] (.5,-4.5);
  
  \draw (-1,-1.5) to [out=270, in=90] (.5,-4.5);
  
  \filldraw
  (.5,-1.5)circle (1.5pt)   (.5,-4.5)circle (1.5pt) 
  (-1,-1.5) circle (1.5pt)   (0,-1.5) circle (1.5pt)   (1,-1.5) circle (1.5pt)  
(1.5,-1.5) circle (1.5pt) (-.5,-1.5)circle (1.5pt) 
  (-1,-4.5) circle (1.5pt)   (0,-4.5) circle (1.5pt)   (1,-4.5) circle (1.5pt)  
 (1.5,-4.5)circle (1.5pt)

(-.5,-4.5) circle (1.5pt)
  (-0.5,-5) circle (1.5pt);

    \draw[white,fill=white] (-.95,-1.9) circle (6 pt);
\draw  (-.95,-1.9) circle (6 pt) node[text=black, scale=.5] {$h_6h$};

    \draw[white,fill=white] (-.46,-1.9) circle (6 pt);
\draw  (-.46,-1.9) circle (6 pt) node[text=black, scale=.5] {$h_5h'$};

    \draw[white,fill=white] (.04,-1.9) circle (6 pt);
\draw  (.04,-1.9) circle (6 pt) node[text=black, scale=.5] {$h_4h'$};

    \draw[white,fill=white] (.54,-1.9) circle (6 pt);
\draw  (.54,-1.9) circle (6 pt) node[text=black, scale=.5] {$h_3h'$};

    \draw[white,fill=white] (1,-1.9) circle (6 pt);
\draw  (1,-1.9) circle (6 pt) node[text=black, scale=.5] {$h_2h$};

    \draw[white,fill=white] (1.5,-1.9) circle (6 pt);
\draw  (1.5,-1.9) circle (6 pt) node[text=black, scale=.5] {$h_1h$};

   \filldraw[fill=white, draw=black] (1.58,-4.1) rectangle (.38,-4.35);
\draw  (1,-4.22)  node[text=black, scale=.5] {$h$};  
   
   \filldraw[fill=white, draw=black] (-1.06,-4.1) rectangle (0.08,-4.35);
\draw  (-.5,-4.22)  node[text=black, scale=.5] {$h'$};  

 \end{scope}

 \begin{scope}[xshift=4.5cm]
  \draw
  (0,-1.5) to [out=270, in=90] (-1,-3)
 
 (0.5,-3) to [out=270, in=90] (-1,-4.5)
   (1,-3) to [out=270, in=90] (-.5,-4.5)
    (1.5,-3) to [out=270, in=90] (0,-4.5)
  
   (.5,-1.5) to [out=270, in=90] (-.5,-3)
  (1,-1.5) to [out=270, in=90] (0,-3)
  (1.5,-1.5) -- (1.5,-3);
  \draw [line width=4pt, white]
  (-.5,-1.5) to [out=270, in=90] (1,-3)
  
  (-1,-1.5) to [out=270, in=90] (.5,-3)
  
      (-1,-3) to [out=270, in=90] (.5,-4.5)
       (-.5,-3) to [out=270, in=90] (1,-4.5)  
   (0,-3) to [out=270, in=90] (1.5,-4.5);
  \draw
  (-1,-1.5) to [out=270, in=90] (.5,-3)
  (-.5,-1.5) to [out=270, in=90] (1,-3)

     (-1,-3) to [out=270, in=90] (.5,-4.5)
       (-.5,-3) to [out=270, in=90] (1,-4.5)  
   (0,-3) to [out=270, in=90] (1.5,-4.5)

  (-1,-4.5) -- (-0.5,-5) -- (0,-4.5)  
  (-.5,-4.5) -- (-0.5,-5)
  (.5,-4.5) -- (1,-5) -- (1.5,-4.5) (1,-5)-- (1,-4.5);
  \filldraw
  (.5,-1.5)circle (1.5pt)   (.5,-4.5)circle (1.5pt) 
  (-1,-1.5) circle (1.5pt)   (0,-1.5) circle (1.5pt)   (1,-1.5) circle (1.5pt)  
(1.5,-1.5) circle (1.5pt) (-.5,-1.5)circle (1.5pt) 
  (-1,-4.5) circle (1.5pt)   (0,-4.5) circle (1.5pt)   (1,-4.5) circle (1.5pt)  
 (1.5,-4.5)circle (1.5pt)

(-.5,-4.5) circle (1.5pt)
  (-0.5,-5) circle (1.5pt);
  
     \draw[white,fill=white] (-1,-3) circle (6 pt);
\draw  (-1,-3) circle (6 pt) node[text=black, scale=.5] {$h_3h'$}; 
     \draw[white,fill=white] (-.5,-3) circle (6 pt);
\draw  (-.5,-3) circle (6 pt) node[text=black, scale=.5] {$h_4h'$}; 
     \draw[white,fill=white] (0,-3) circle (6 pt);
\draw  (0,-3) circle (6 pt) node[text=black, scale=.5] {$h_5h'$}; 
     \draw[white,fill=white] (0.5,-3) circle (6 pt);
\draw  (0.5,-3) circle (6 pt) node[text=black, scale=.5] {$h_1h$}; 
     \draw[white,fill=white] (1,-3) circle (6 pt);
\draw  (1,-3) circle (6 pt) node[text=black, scale=.5] {$h_2h$}; 
     \draw[white,fill=white] (1.5,-3) circle (6 pt);
\draw  (1.5,-3) circle (6 pt) node[text=black, scale=.5] {$h_6h$};

   \filldraw[fill=white, draw=black] (-1.04,-4.1) rectangle (.18,-4.35);
\draw  (-.5,-4.22)  node[text=black, scale=.5] {$h$};  
   
   \filldraw[fill=white, draw=black] (0.3,-4.1) rectangle (1.55,-4.35);
\draw  (1.1,-4.22)  node[text=black, scale=.5] {$h'$};

 \end{scope}
\end{tikzpicture}

}
\caption{Dangling.}
\label{fig:dangle}
\end{figure}

We can identify the labeled braid group $B_n(H)$ with a subgroup of~$\spraige_{n,n}$ via $(b, \lambda)\mapsto (1_n,(b, \lambda),1_n)$. In particular, for any $n,m\in\BN$ there is a right
action of the labeled braid group $B_m(H)$ on $\spraige_{n,m}$, by right multiplication. 
We can quotient out this action and we refer to this quotient as \emph{dangling}.  See
Figure~\ref{fig:dangle} for an example of the dangling action of $B_2(H)$
on~$\spraige_{6,2}$.

For $\sigma\in\spraige_{n,m}$, denote by $[\sigma]$ the orbit of
$\sigma$ under this action, and call $[\sigma]$ a \emph{dangling
$(n,m)$-spraige}.  We can also refer to a dangling $n$-spraige or
dangling spraige.  Note that the action of $B_m(H)$ preserves the property of being
a braige or elementary braige, so the notions of dangling
braiges and dangling elementary braiges are well-defined.

Let $\Poset$ denote the set of all dangling spraiges,
$\Poset_{n,m}$ the dangling $(n,m)$-spraiges and $\Poset_n$ the set of dangling $n$-spraiges.   Note that if $\sigma\in\spraige_{n,m}$ and~$\tau_1, \tau_2
\in \spraige_{m,k}$ with $[\sigma \ast \tau_1] =
[\sigma \ast \tau_2]$, then $[\tau_1]=[\tau_2]$. We will refer to this fact as
\emph{left cancellation}.

There is also a poset structure on $\Poset$.  For $x,y\in\Poset$, with
$x=[\sigma_x]$, say that~$x\le y$ if there exists a forest~$F$
with~$m$ leaves such that~$y=[\sigma_x \ast (F,(\id, \iota),1_m)]$.  In other
words,~$x\le y$ if~$y$ is obtained from~$x$ via splitting.  It is easy
to see that this is a partial ordering.  Also, if~$x\in\Poset_n$
and~$y\in\Poset$ with~$x\le y$ or~$y\le x$, then~$y\in\Poset_n$.
In other words, two elements are comparable only if they have the same number of heads. We further define a
relation~$\preceq$ on~$\Poset$ as follows.  If~$x = [\sigma_x] \in
\Poset$ and $y\in\Poset$ such that~$y = [\sigma_x \ast
\lambda^{(n)}_J]$ for some $n\in\BN$
and~$J\subseteq\{1,\dots,n\}$, write~$x\preceq y$.  That is,~$x\preceq
y$ if $y$ is obtained from~$x$ via elementary splitting, and this is a
well-defined relation with respect to dangling.  If~$x\preceq y$
and~$x\neq y$, then write $x\prec y$.  Note that~$\preceq$ and~$\prec$
are not transitive, though it is true that if~$x\preceq z$ and~$x\le
y\le z$, then $x\preceq y$ and~$y\preceq z$.  This is all somewhat similar to the corresponding
situation for~$F$ and~$V$ discussed for example in
~\cite[Section 4]{Br92}.
We remark that a totally analogous construction yields the notion of a
\emph{dangling pure spraige} and \emph{dangling cyclic spraige}, where the dangling is now via the action
of the pure labeled braid group or cyclic labeled braid group.  We also have dangling pure/cyclic braiges and
dangling elementary pure/cyclic braiges.  All of the essential results above
still hold.

\subsection{The Stein space}\label{sec:def_stein_space}
In this subsection, we construct a space $X$ on which $bV_{d,r}(H)$ acts and which
we call the \emph{Stein space} for~$bV_{d,r}(H)$. A similar space can also
be constructed using pure braids and cyclic braids to get spaces $X(bF_{d,r}(H))$ and $X(bT_{d,r}(H))$ on which~$bF_{d,r}(H)$ and $bT_{d,r}(H)$ act, respectively, and we will say more about this at the end of the section.

Once we have the Stein space, we will apply Brown's criterion to the action on $X$ to deduce the positive finiteness properties of the labeled braided Higman--Thompson groups. First we recall Brown's criterion \cite[Theorem 2.2, 3.2]{Br87}.  Recall that a filtration $(X_j)_{j\ge
1}$ of $X$ is called \emph{essentially $n$-connected} if for every $i\geq 1$, there exists $i'\geq i$ such that $\pi_l(X_i\to X_{i'})$ is trivial for all $l\leq n$.

\begin{thm*}[Brown's criterion]
\label{brownscriterion}
Let $n\in \BN$ and assume a group $G$ acts on an $(n-1)$-connected CW-complex $X$. Assume that the stabilizer of every $k$-cell of $X$ is of type $F_{n-k}$. Let $\{X_j\}_{j\ge
1}$ be a filtration of $X$ such that each~$X_j$ is finite
$\textnormal{mod}~G$. Then $G$ is of type $F_n$ if and only if $\{X_j\}_j$ is essentially $(n-1)$-connected.
\end{thm*}
 
For several of the results in this section, we direct the reader to \cite[Section 2]{BFM+16}. Although, that paper only directly addresses the case where $d=2$, $r=1$, and $H$ is the trivial group, their proofs often generalize directly to higher $d$, $r$ and arbitrary $H$. We first consider only the groups $bV_{d,r}(H)$ and then remark on $bF_{d,r}(H)$ and $bT_{d,r}(H)$ at the end.

Our starting point is the poset $\Poset_r$ of dangling $r$-spraiges,
i.e., dangling spraiges with $r$ heads.  Consider the geometric realization $|\Poset_r|$, i.e., the
simplicial complex with a~$k$-simplex for every chain $x_0<\cdots<x_k$
in $\Poset_r$.  We will refer to $x_k$ as the \emph{top} of the
simplex and $x_0$ as the \emph{bottom}.  Call such a simplex
\emph{elementary} if~$x_0\preceq x_k$. 

\begin{defn}\label{def:stein_space}
Define the \emph{Stein space} $X$ for $bV_{d,r}(H)$ to be the subcomplex
of~$|\Poset_r|$ consisting of all elementary simplices. 
\end{defn}

Since faces of elementary simplices are elementary, this is indeed a
subcomplex. 

There is also a coarser cell decomposition of $X$, as
a cubical complex, which we now describe.  First we define the cubes as well as their top and bottom. 

\begin{defn}\label{defn:cubecomplex}
For $x\le y$ define,
the closed interval $[x,y]:= \{z\mid x\le z\le y\}$.  Similarly,
define the open and half-open intervals $(x,y)$, $(x,y]$ and $[x,y)$.
Note that if $x\preceq y$, then the closed interval $[x,y]$ is a
Boolean lattice, and so the simplices in its geometric realization
fit together into a cube.  The \emph{top} of the cube is $y$ and the
\emph{bottom} is~$x$.  
\end{defn}

Now observe that every elementary simplex is contained in such a
cube, and the face of any cube is clearly another cube.  Also, the
intersection of cubes is either empty or is itself a cube; this is
clear since if $[x,y]\cap[z,w]\neq\emptyset$, then $y$ and $w$ have a
lower bound, and we get that $[x,y]\cap[z,w]=[\sup(x,z),\inf(y,w)]$. Note that as in \cite[Proposition 2.1]{BFM+16}, any two elements in $\Poset_r$ have a least upper bound and when two elements have a lower bound, they have a greatest lower bound.
Therefore,~$X$ has the structure of a cubical complex, in
the sense of~\cite[p.~112, Definition~7.32]{BH99}.  

Recall that a poset $(Y, \ll)$ is called \emph{conically contractible} if there is a $y_0$ in $Y$ and a map $g: Y \rightarrow Y$ such that $z \gg g(z) \ll y_0$ for all $z$ in  $Y$. A consequence of a poset being conically contractible is that its geometric realization is contractible. See the discussion in \cite[Section 1.5]{Qui78} for more details.

\begin{lem}\label{lem:contractibleinterval}
For $x< y$ with $x\not\prec y$, $|(x,y)|$ is contractible.
\end{lem}

\Proof
We will prove that $(x,y)$ is conically contractible and that then implies the lemma. For $x, y \in \Poset_r$, we will declare $x \ll y$ if and only if $y\leq x$. Now given any $z \in (x,y]$, define $g(z)$ to be the largest element of $[x,z]$ such that $x\preceq g(z)$. By our hypothesis, $g(z)$ is in  $[x,y)$ and also clearly in $(x,y]$, so therefore $g(z)\in (x,y)$. Let $y_0=g(y).$ Note that for any $z\in (x,y)$, we have $g(z) \leq y_0$. Whence $z \gg g(z) \ll y_0$ and $(x,y)$ is conically contractible.
\qed

\begin{cor}\label{cor-steincontract}
The space $X$ is contractible.
\end{cor}
\Proof
We first see that $\Poset_r$ is directed since any two elements have a least upper bound, as discussed after Definition \ref{defn:cubecomplex}. Therefore, $|\Poset_r|$ is contractible.

Now, as in \cite{BFM+16}, we will build up from $X$ to $|\Poset_r|$ by attaching new subcomplexes in such a way as to not change the homotopy type. Given a closed interval $[x,y]$, define $r([x,y]):=f(y)-f(x)$. We attach the contractible subcomplexes $|[x,y]|$ for $x\not\preceq y$ to $X$ in increasing order of $r$ value, attaching $|[x,y]|$ along $|[x,y)\cup (x,y]|$. This is the suspension of $|(x,y)|$ and hence is contractible by Lemma~\ref{lem:contractibleinterval}.  Therefore, attaching $|[x,y]|$ does not change the homotopy type and we conclude that $X$ is contractible.
\qed

There is a natural action of $bV_{d,r}(H)$ on the vertices of $X$.  Namely,
for $g\in bV_{d,r}(H)$ and~$\sigma\in\spraige_r$ with $x=[\sigma]$, define $gx
:= [g\ast\sigma]$.  This action preserves the
relations~$\le$ and~$\preceq$, and thus extends to an action on the
whole space.

For each $m\in\BN$, define $X^{\le m}$ to be the full subcomplex of~$X$
spanned by vertices~$x$ with $f(x)\le m$.  Note that the $X^{\le m}$
is invariant under the action of~$bV_{d,r}$. Now the same proof in \cite[Lemma 2.5]{BFM+16} works to show the following.

\begin{prop}
\label{prop:stein_space_cible}
For each $m\geq 1$, the sublevel set $X^{\leq m}$ is finite modulo $bV_{d,r}(H)$.
\end{prop}

We now consider the vertex and cell stabilizers.

\begin{defn}
Let $J\subseteq\{1,\dots,m\}$.  Let $b\in B_m$ and let $\rho_b$ be the
corresponding permutation in $S_m$.  If $\rho_b$ stabilizes $J$
set-wise, call $b$ a $J$-\emph{stabilizing braid}. Let $B_m^J\le
B_m$ be the subgroup of~$J$-stabilizing braids and $B_m^J(H) \cong H^m\rtimes B_m^J$.
\end{defn}

\begin{prop}\label{cor:cell_stabs}
Let $x$ be a vertex in $X$, with $f(x)=n$ and $x=[\sigma]$, and let
$F^{(m)}_J$ be an elementary forest.  If
$y=[\sigma\ast\lambda^{(m)}_J]$, then the stabilizer in~$bV_{d,r}(H)$ of the
cube $[x,y]$ is isomorphic to~$B_m^J(H)$.  In particular, if $H$ is of type $F_n$, then so are the cell stabilizers.
\end{prop}
\Proof
The first part of the statement follows directly from the proofs of Lemma 2.6 and Corollary 2.8 in \cite{BFM+16}. For the second statement, observe that $B_m^J(H)$ has finite index in $B_m(H)\cong H^m \rtimes B_m$, that the braid groups are of type $F_\infty$ and that finiteness properties are preserved under extensions by a group of type $F_\infty$ \cite[Theorem 7.2.21]{Ge08}.
\qed

The complex $X$ and the filtration $\{X^{\le m}\}_{m}$ has so far been shown to
satisfy all the conditions of Brown's criterion save one, namely that
the filtration $\{X^{\le m}\}_m$ is essentially $(n-1)$-connected.  We will prove this in
Corollary~\ref{cor:desc_link_conn}, using the Morse Lemma. 

Note that every cell of~$X$ has a unique vertex maximizing~$f$, so $f$ is a height function. Hence we can inspect the connectivity of $\{X^{\le m}\}_m$
by looking at descending links with respect to~$f$.  In the rest of
this section, we describe a convenient model for the descending links,
and then analyze their connectivity in the following sections.

Recall that we identify $\Poset_r$ with the vertex set of $X$, and
cubes in $X$ are (geometric realizations of) intervals $[y,x]$
with $x,y\in\Poset_r$ and $y\preceq x$.  For $x\in\Poset_r$, the
descending star $\dst(x)$ of $x$ in $X$ is the set of cubes $[y,x]$
with top $x$.  For such a cube $C=[y,x]$ let $\Bot(C):= y$ be the
map giving the bottom vertex.  This is a bijection from the
set of such cubes to the set $D(x) := \{y\in\Poset_r\mid y\preceq
x\}$.  The cube $[y',x]$ is a face of~$[y,x]$ if and only if
$y'\in[y,x]$, if and only if $y'\ge y$.  Hence $C'$ is a face of $C$
if and only if $\Bot(C')\ge \Bot(C)$, so $\Bot$ is an order-reversing
poset map.  By considering cubes~$[y,x]$ with $y\neq x$ and
restricting to $D(x) \setminus \{x\}$, we obtain a description of~$\dlk(x)$.  Namely, a simplex in~$\dlk(x)$ is a dangling spraige~$y$
with~$y\prec x$, the rank of the simplex is the number of elementary
splits needed to get from $y$ to $x$ (so the number of elementary
merges to get from~$x$ to~$y$) and the face relation is the reverse of
the relation~$<$ on $D(x)\setminus\{x\}$.  Since $X$ is a cubical
complex, $\dlk(x)$ is a simplicial complex.

We proceed to describe a convenient model for the descending link. If $f(x)=m$, then thanks to left cancellation, $\dlk(x)$ is isomorphic to the simplicial complex $\elbraigecpx_d^m$ of dangling
elementary $m$-braiges $[(1_m,(b, \lambda),F_J^{(m-(d-1)|J|)})]$ for
$J\neq\emptyset$, with the face relation given by the reverse of the
ordering~$\le$ in $\Poset_r$.  See Figure~\ref{fig:desc_lk_to_EB} for
an idea of the correspondence between~$\dlk(x)$ and $\elbraigecpx_d^m$.  We will
usually draw braiges as emerging from a horizontal line, as a visual
reminder of this correspondence. We will prove that $\elbraigecpx_d^m$ is highly connected in
Corollary~\ref{cor:desc_link_conn}. 

\begin{figure}[h]
\centering
\begin{tikzpicture}

  \filldraw[lightgray]
   (0,0) -- (1,1) -- (2,0) -- (1,-1) -- (0,0);
  \draw
   (2,0) -- (1,-1) -- (0,0);
  \draw[darkgray, line width=1.5pt]
   (0,0) -- (1,1) -- (2,0);
  \filldraw
   (0,0) circle (1.5pt)   (1,1) circle (1.5pt)   (2,0) circle (1.5pt)   (1,-1)
circle (1.5pt);

  \draw
   (0.6,1.4) to [out=90, in=90, looseness=2] (1.4,1.4) -- (0.6,1.4)
   (0.75,1.4) -- (0.75,1.2) (0.85,1.4) -- (0.85,1.2)   (0.95,1.4) -- (0.95,1.2)   (1.05,1.4) --
(1.05,1.2)   (1.15,1.4) -- (1.15,1.2) (1.25,1.4) -- (1.25,1.2);
  \node at (1,1.52) {$x$};

 \begin{scope}[xshift=-1.2cm,yshift=-1cm]
  \draw
   (0.6,1.4) to [out=90, in=90, looseness=2] (1.4,1.4) -- (0.6,1.4)
   (0.75,1.4) -- (0.75,1.2)-- (0.85,1.1) (0.85,1.4) -- (0.85,1.1)-- (0.95,1.2)    (0.95,1.4) -- (0.95,1.2)   (1.05,1.4) --
(1.05,1.2)   (1.15,1.4) -- (1.15,1.2) (1.25,1.4) -- (1.25,1.2);
  \node at (1,1.52) {$x$};
 \end{scope}

 \begin{scope}[xshift=1.2cm,yshift=-1cm]
  \draw
   (0.6,1.4) to [out=90, in=90, looseness=2] (1.4,1.4) -- (0.6,1.4)
   (0.75,1.4) -- (0.75,1.2) (0.85,1.4) -- (0.85,1.2)  (0.95,1.4) -- (0.95,1.2)   (1.05,1.4) --
(1.05,1.2) -- (1.15,1.1)  (1.15,1.4) -- (1.15,1.1)-- (1.25,1.2) -- (1.25,1.4);
  \node at (1,1.52) {$x$};
 \end{scope}

 \begin{scope}[yshift=-3cm]
  \draw
  (0.6,1.4) to [out=90, in=90, looseness=2] (1.4,1.4) -- (0.6,1.4)
   (0.75,1.4) -- (0.75,1.2)-- (0.85,1.1) (0.85,1.4) -- (0.85,1.1)-- (0.95,1.2)    (0.95,1.4) -- (0.95,1.2)   (1.05,1.4) --
(1.05,1.2) -- (1.15,1.1)  (1.15,1.4) -- (1.15,1.1)-- (1.25,1.2) -- (1.25,1.4);
  \node at (1,1.52) {$x$};
 \end{scope}

  \node at (3.5,0) {$\longleftrightarrow$};

  \draw[line width=1.5pt, gray]
   (4.8,0) -- (6.8,0);
  \filldraw[darkgray]
   (4.8,0) circle (1.5pt)   (6.8,0) circle (1.5pt);

 \begin{scope}[xshift=4.7cm, yshift=-0.2cm, scale=2]
  \draw
   (-0.25,0.5) -- (-0.25,0.3) -- (-0.15,0.2) -- (-0.15,0.3) -- (-0.15,0.5)  (-0.15,0.2) -- (-0.05,0.3) -- (-0.05,0.5) 
(0.05,0.5) -- (0.05,0.3)   (0.15,0.5) -- (0.15,0.3) (0.25,0.5) -- (0.25,0.3);

  \draw
   (-0.3,0.5) -- (0.3,0.5);

(0.5pt)   (-0.05,0.3) circle (0.5pt)   (-0.05,0.5) circle (0.5pt)   (0.05,0.5)
circle (0.5pt)   (0.05,0.3) circle (0.5pt)   (0.15,0.5) circle (0.5pt)  
(0.15,0.3) circle (0.5pt);
 \end{scope}

 \begin{scope}[xshift=6.9cm, yshift=-0.2cm, scale=2]
   \draw
   (0.25,0.5) -- (0.25,0.3) -- (0.15,0.2) -- (0.15,0.3) -- (0.15,0.5)  (0.15,0.2) -- (0.05,0.3) -- (0.05,0.5) 
(-0.05,0.5) -- (-0.05,0.3)   (-0.15,0.5) -- (-0.15,0.3) (-0.25,0.5) -- (-0.25,0.3);

  \draw
   (-0.3,0.5) -- (0.3,0.5);

(0.5pt)   (-0.05,0.3) circle (0.5pt)   (-0.05,0.5) circle (0.5pt)   (0.05,0.5)
circle (0.5pt)   (0.05,0.3) circle (0.5pt)   (0.15,0.5) circle (0.5pt)  
(0.15,0.3) circle (0.5pt);
 \end{scope}

 \begin{scope}[xshift=5.8cm, yshift=-1.3cm, scale=2]
  \draw
   (0.25,0.5) -- (0.25,0.3) -- (0.15,0.2) -- (0.15,0.3) -- (0.15,0.5)  (0.15,0.2) -- (0.05,0.3) -- (0.05,0.5) 
(-0.25,0.5) -- (-0.25,0.3) -- (-0.15,0.2) -- (-0.15,0.3) -- (-0.15,0.5)  (-0.15,0.2) -- (-0.05,0.3) -- (-0.05,0.5) ;

  \draw
   (-0.3,0.5) -- (0.3,0.5);

(0.5pt)   (-0.05,0.3) circle (0.5pt)   (-0.05,0.5) circle (0.5pt)   (0.05,0.5)
circle (0.5pt)   (0.05,0.3) circle (0.5pt)   (0.15,0.5) circle (0.5pt)  
(0.15,0.3) circle (0.5pt);
 \end{scope}
\end{tikzpicture}
\caption{The correspondence between $\dlk(x)$ and $\elbraigecpx_d^m$.}
\label{fig:desc_lk_to_EB}
\end{figure}

We end this section with some remarks on $bF_{d,r}(H)$ and $bT_{d,r}(H)$.  Restricting to
pure labeled braids or cyclic labeled braids everywhere in this section does not affect any of the
proofs, so we can simply say that $X(bF_{d,r}(H))$ and $X(bT_{d,r}(H))$ are the contractible
cubical complexes of dangling pure $r$-spraiges, understood in the same
way as $X$ (though now dangling is only via pure labeled braids or cyclic labeled braids).  We will
also denote by $f$ the height function ``number of feet'' on
$X(bF_{d,r}(H))$ and $X(bT_{d,r}(H))$.  The filtration is still cocompact and the stabilizers are
still of type $F_n$ whenever $H$ is, being finite index subgroups of the corresponding labeled braid
groups.  As for descending links, the descending link of a dangling
pure $(r,m)$-spraige in $X(bF_{d,r}(H))$ is isomorphic to the simplicial
complex~$\elpbraigecpx_d^m$ of dangling elementary labeled pure~$m$-braiges and in $X(bT_{d,r}(H))$ to the simplicial complex $\elcbraigecpx_d^m$ of dangling elementary labeled cyclic $m$-braiges.

\subsection{The complex related to the descending link for $bV_{d,r}(H)$} \label{Sct:conn-dcpx}
Let $S_{b,m}^g$ be a compact oriented surface of genus $g$ with $b$ boundary components and $m$ marked points or punctures such that the marked points are in the interior of the surface.

 A \emph{$d$-arc} on the surface is an embedded path in $S^g_{b,m}\backslash \partial S^g_{b,m} $ that begins and ends at marked points and passes through a total of precisely $d$ marked points. We call a collection of $d$-arcs $\{\alpha_0, \alpha_1, \dots, \alpha_k\}$ a \emph{$d$-arc system} if for all $i\neq j$, the $d$-arcs $\alpha_i$ and $\alpha_j$ are disjoint up to isotopy. Note that the isotopies here are required to fix the marked points.

\begin{defn} The \emph{$d$-arc matching complex} $\marc_d(S^g_{b,m})$ on $S^g_{b,m}$ is the simplicial complex with a $k$-simplex for each isotopy class of a $d$-arc system $\{\alpha_0, \alpha_1 \dots, \alpha_k\}$ and the face relation given by the subset relation. 
\end{defn}

When $d=2$, our complex is just the matching complex  $\marc(\Gamma_m)$ over the surface $S_{b,m}^g$ in \cite[Section 3]{BFM+16}.

\begin{lem}\label{lem:arcsoverclasses}
Assume $m\geq 3$, given finitely many homotopy classes of $d$-arcs $[\alpha_0], [\alpha_1], \dots, [\alpha_k]$ there exist representatives $\alpha_0, \alpha_1, \dots, \alpha_k$ such that $|\alpha_i \cap \alpha_j|$ is minimal among all representatives of $[\alpha_i]$ and $[\alpha_j]$ for $0\leq i \leq j \leq k$.  In particular, any simplex is represented by disjoint $d$-arcs.
\end{lem}

\Proof Notice that each $d$-arc corresponds to a collection of $(d-1)$ $2$-arcs and so  it suffices to show it is true for $2$-arcs. But this was proven  in \cite[Lemma 3.2]{BFM+16}.  Basically, one puts a hyperbolic metric on the interior of $S^g_{b,m}$ (viewing marked points as punctures) and replaces each $2$-arc in $\alpha_i$ by a geodesic connecting the two punctures.
\qed

The lemma allows us to consider actual arcs instead of homotopy classes of arcs when $m\geq 3$ which  we will do in the rest of this section. 

We now proceed to give a connectivity bound for the complex $\marc_d(S^g_{b,m})$ closely following the strategy in \cite[Section 3.3]{BFM+16}. Let us label all the marked points in $S_{b,m}^g$ as $\{1,2,\cdots, m\}$. We put a weight on all the marked points in $S_{b,m}^g$ via the following rule: if $p\leq d$, we assign its weight to be $2^{p-1}$; if $p>d$, we assign its weight to be $0$. With this we can define a weight function $q$ on any vertex $\alpha$ in  $\marc_d(S^g_{b,m})$  by assigning $q(\alpha)$ to be the total weight of the marked points that  $\alpha$ passes through.  Note that  the zero set of $q$, which we denote as $\marc_d(S^g_{b,m}) ^{q=0}$,  can be identified with the complex $\marc_d(S_{b+d,m-d}^g)$. Here the surface $S_{b+d,m-d}^g$ is obtained from $S_{b,m}^g$ by deleting a small open disk around those marked points with positive weight. Now $q$ defines a height function on the relative complex $(\marc_d(S_{b,m}^g),\marc_d(S_{b,m}^g)^{q=0})$. We will use the $q$ to analyze the connectivity of $\marc_d(S_{b,m}^g)$.

\begin{thm}\label{thm:conn-i-d-arc-complex}
For any $d\geq2$, the complex $\marc_{d}(S_{b,m}^g)$ is $(\lfloor\frac{m+1}{2d-1}\rfloor-2)$-connected.
\end{thm}

\Proof We prove the theorem by induction on $m$. Note that when $m\geq d$, the complex $\marc_{ d}(S_{b,m}^g)$ is nonempty, hence the theorem is valid for $m\leq 4d-4$. Now assume $m>4d-4$. Given any vertex $\alpha$ in $\marc_d(S_{b,m}^g)$ such that $q(\alpha)\neq 0$, the descending link of $\alpha$ is the full subcomplex of $\marc_d(S_{b,m}^g)$ with vertices $\alpha'$ where $q(\alpha') < q(\alpha)$ and $\alpha'$ is disjoint from $\alpha$. The point here is that the descending link is again a $d$-arc matching  complex  over some surface, and the new surface now has at least $m-2d+1$ marked points. In fact, in the worst case, $q(\alpha) =1$ and $\alpha$ contains the marked point of weight $1$ and $d-1$ marked points of weight $0$. Thus any vertex $\alpha'\in \dlk (\alpha)$ must have weight $0$. This means $\alpha'$ only passes those marked points of weight zero outside $\alpha$.  Therefore, $\dlk (\alpha)$ in this case can be identified with $\marc_d(S_{b+d,m-2d+1}^g)$, where the surface $S_{b+d,m-2d+1}^g$ is obtained from $S_{b,m}^g$ by deleting a small open neighborhood of $\alpha$ and the marked points $2,\cdots,d$. By induction, the descending link is at least $(\lfloor\frac{m+1}{2d-1}\rfloor-3)$-connected.

By the Morse lemma (cf. Lemma \ref{lemm-Morse}), the pair $(\marc_d(S_{b,m}^g),\marc_d(S_{b,m}^g)^{q=0})$ is $(\lfloor\frac{m+1}{2d-1}\rfloor-2)$-connected, that is, the inclusion ~$\iota\colon
\marc_d(S_{b,m}^g)^{q=0}\hookrightarrow \marc_d(S_{b,m}^g)$ induces an isomorphism on $\pi_n$ for $n \le \lfloor\frac{m+1}{2d-1}\rfloor-3$ and an epimorphism for $n = \lfloor\frac{m+1}{2d-1}\rfloor - 2$. We could now invoke induction and use that $\marc_d(S_{b,m}^g)^{q=0}$ is at least $(\lfloor\frac{m+1}{2d-1}\rfloor-3)$-connected to conclude that $\marc_d(S_{b,m}^g)$ is $(\lfloor\frac{m+1}{2d-1}\rfloor-3)$-connected as well. However, since we want $\marc_d(S_{b,m}^g)$ to be $(\lfloor\frac{m+1}{2d-1}\rfloor-2)$-connected, we need a different argument and we may as well apply this for all $n$. It suffices to show that $\pi_n(\marc_d(S_{b,m}^g)^{q=0}\hookrightarrow \marc_d(S_{b,m}^g))$ is trivial
for $n\leq \lfloor\frac{m+1}{2d-1}\rfloor -2$. In other words, for any $n\leq\lfloor\frac{m+1}{2d-1}\rfloor -2$, every map  $\overline{\psi}\colon S^n\to\marc_d(S_{b,m}^g)^{q=0}$  can be homotoped to a constant map in $\marc_d(S_{b,m}^g)$.

First we check the hypothesis on $\marc_d(S_{b,m}^g)$ allows us to
apply Lemma~\ref{lem:injectifying}, namely that the link of a
$k$-simplex should be $(n-k-2)$-connected.  A $k$-simplex~$\sigma$ is
determined by $k+1$ disjoint d-arcs.  Hence, the link of $\sigma$ is
isomorphic to $\marc_d(S_{b+k+1,m-(k+1)d}^g)$.  By induction, this is $(\lfloor\frac{m-(k+1)d+1}{2d-1}\rfloor-2)$-connected, which is at least $(n-k-2)$-connected.

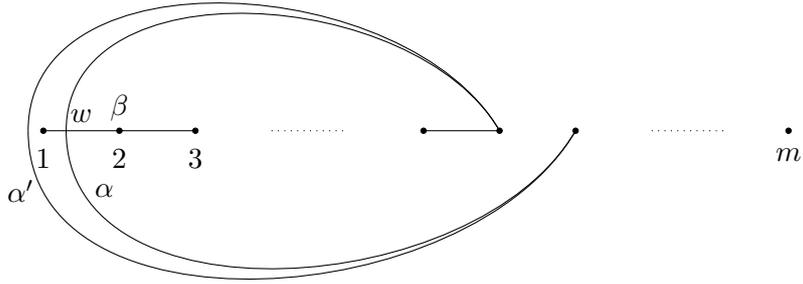
\begin{figure}
\centering
\begin{tikzpicture}
\coordinate (a) at (-5.5,0);
\coordinate (b) at (-4.5,0);
\coordinate (c) at (-3.5,0);

\draw (a) -- (b) node [anchor=south] {$\beta$};

\draw (1.5,0) to [out=-120, in=-90] (-5.2,0);
\draw (0.5,0) to [out=120, in=90] (-5.2,0);

\draw (-5.2,0) + (0.5,-0.8) node {$\alpha$};

\draw (-5,0.2) node  {$w$};

\draw (1.5,0) to [out=-120, in=-90] (-5.7,0);
\draw (0.5,0) to [out=120, in=90] (-5.7,0);

\draw (-5.7,0) + (-0.1,-0.8) node {$\alpha'$};

\draw[dotted] (-2.5,0)--(-1.5,0);

\draw[dotted] (2.5,0)--(3.5,0);

\draw (-4.5,0)--(-3.5,0);
\draw (-0.5,0)--(.5,0);
\filldraw

	(a) circle (1pt) node[below=3pt] {$1$}
	(b) circle (1pt) node[below=3pt] {$2$}
	(c) circle (1pt) node[below=3pt] {$3$}

	(-0.5,0) circle (1pt)
	(0.5,0) circle (1pt)
	(1.5,0) circle (1pt)

	(4.3,0) circle (1pt) node[below=3pt] {$m$};
\end{tikzpicture}
\caption{Pushing part of the $d$-arc $\alpha$ over the marked point  $1$ to obtain the
$d$-arc~$\alpha'$, as described in the proof of
Theorem~\ref{thm:conn-i-d-arc-complex}. }
\label{fig:matching_flow}
\end{figure}

Let $S^n$ be a combinatorial $n$-sphere.  Let
$\overline{\psi}\colon S^n\to\marc_d(S_{b,m}^g)^{q=0}$ be a simplicial map and
let $\psi := \iota \circ \overline{\psi}$.  It suffices by simplicial
approximation~\cite[Theorem~3.4.8]{Span66} to homotope~$\psi$ to a constant
map. By Lemma~\ref{lem:injectifying}, we may assume~$\psi$ is simplexwise
injective.  Fix $\beta$ to be a $d$-arc passing through the marked points $\{1,2,\cdots, d\}$ according to the order. Then $\beta$ is the concatenation of $d-1$ arcs $\beta_1,\cdots,\beta_{d-1}$, where $\beta_i$ is an arc connecting the marked point $i$ and $i+1$. We claim that~$\psi$ can be homotoped in $\marc_d(S_{b,m}^g)$ to land in the star of $\beta$, which
will finish the proof. We will proceed in a similar way to the Hatcher flow \cite{Hat91}. Note first that none of the $d$-arcs in the image of $\psi$ pass through any marked points with positive weight, but among the finitely many such d-arcs, some might intersect nontrivially with ~$\beta$. Pick one, say $\alpha$, intersecting $\beta$ at a point, say $w$, closest along $\beta$ to the marked point $1$, and let~$x$ be a vertex of~$S^n$ mapping to $\alpha$. Without loss of generality, we can assume further that the intersection point lies in $\beta_{1}$. By simplexwise injectivity, none of the
vertices in $\Lk_{S^n}(x)$ map to $\alpha$. We will replace the arc component of $\alpha$ which contains $w$ by another arc $\alpha'$ with the same endpoints but on the other side of the marked point $1$. In fact, the new arc component together with the part of $\alpha$ that contains $w$ bound a disk $D$ whose interior contains no boundary components or marked points
other than the marked point $1$.  See Figure~\ref{fig:matching_flow} for an example.   Note that
there is no edge from $\alpha$ to $\alpha'$, so none of the vertices
in~$\Lk_{S^n}(x)$ map to $\alpha'$. Note also that
$\psi(\Lk_{S^n}(x)) \subseteq \Lk(\alpha')$ by our choice of $\alpha$.

We now want to apply the mutual link trick (cf. Lemma \ref{lemma-replace-trick}) to homotope the map $\psi'$ to a new simplexwise injective map $\psi'\colon S^n\to\marc_d(S_{b,m}^g)$ that sends
the vertex $x$ to $\alpha'$ and sends all other vertices $y$ to
$\psi(y)$. For that we only need to further check that the mutual link $\Lk_{X} (\alpha) \cap \Lk_{X}(\alpha')$ is $(n-1)$-connected. But $\Lk(\alpha)\cap\Lk(\alpha')$ is isomorphic to
$\marc_d(S_{b+1,m-d-1}^g)$, where the surface  $S_{b+1,m-{(d+1)}}^g$  is 
obtained from $S_{b,m}^g$ by removing an open neighborhood of $D\cup \alpha\cup\alpha'$.  Hence by induction $\Lk(\alpha)\cap\Lk(\alpha')$ is
$(\lfloor\frac{m-d}{2d-1}\rfloor-2)$-connected, and in particular~$(n-1)$-connected. In this way, after finitely many steps, we can homotope $\psi$ such that its image is disjoint from $\beta$. In particular, for any $i\geq 2$, we get rid of the intersection of $\psi(S^n)$ with $\beta$ in $i$ steps: first push each intersection with $\beta_i$ to $\beta_{i-1}$,  then to $\beta_{i-2}$, etc. In the last step, we push the intersections off  $\beta_1$. This has the benefit that the disk $D$ between each $\alpha$ and $\alpha'$ will only bound one marked point at a time. At the end, we can assume the image of $\psi$ is disjoint from $\beta$, hence it lies in the star of $\beta$. Therefore, $\psi$ can be homotoped to a constant map. 
\qed

As a by-product of Theorem \ref{thm:conn-i-d-arc-complex}, we also get connectivity bounds of certain disk complexes which might have independent interest. Let us introduce them now.  Given a surface $S_{b,m}^g$,  a $k$-simplex in the \emph{$d$-marked-point-disk complex}  $\BD_d(S_{b,m}^g)$ is an isotopy class of a system of disjointly embedded disks $\langle D_0,D_1,\cdots, D_k\rangle$ such that each disk $D_i$ encloses precisely $d$ marked points in its interior. Here again  the face relation given by the subset relation.  Note that except some singular cases, $\BD_d(S_{b,m}^g)$ can be viewed as a full subcomplex of the curve complex first defined by Harvey in  \cite{Hav81}. In fact, given a disk enclosing $d$ marked points,  we can take their boundary curve which gives a vertex in the curve complex unless the boundary curve bounds a disk, a punctured sphere, or an annulus on the other side.  It might also happen that two disks are disjoint up to isotopy but their boundary curves are isotopic. This case occurs when the surface is a sphere with $2d$ marked points, in which case its $d$-marked-point-disk complex is not a subcomplex of the curve complex.

There is also a canonical map 
\begin{equation}\label{eq:mapfromarctodisk}
N: \marc_d(S_{b,m}^g)\to \BD_{ d}(S_{b,m}^g)
\end{equation}
mapping each $d$-arc to  a small disk tubular neighborhood of it. We have the following.

\begin{cor}\label{cor:conn-i-dcomplex}
The map $N$ is a complete join. In particular, for any $d\geq2$, the complex $\BD_{ d}(S_{b,m}^g) $ is $(\lfloor\frac{m+1}{2d-1}\rfloor-2)$-connected.
\end{cor}

\Proof 
If two systems of $d$-arcs are isotopic, then their disk tubular neighborhoods are isotopic. Hence the map is well-defined on vertices. If a system of $d$-arcs is disjoint, we can choose their disk tubular neighborhoods to be disjoint. This shows $N$ is well-defined and simplexwise injective. To show it is surjective, given any $k$-simplex  $\sigma =\langle D_0,D_1,\cdots,D_k\rangle$, where each $D_i$ is a disk enclosing $d$ marked points, we can choose a $d$-arc in the interior of each disk which passes through the $d$ marked points inside $D$. The fact that any such $d$-arc system lies in the preimage of $\sigma$ says $N^{-1}(\sigma) = N^{-1}(D_0) \ast \cdots \ast N^{-1}(D_k)$. Thus $N$ is a complete join. The lemma now follows from Remark \ref{rem-cjoin}.
\qed

\subsection{The complex related to the descending links for $bF_{d,r}(H)$ and $bT_{d,r}(H)$}  \label{Sct:conn-i-ldcpx}

In this subsection, we introduce and calculate the connectivity for the complex related to the descending links for $bF_{d,r}(H)$ and $bT_{d,r}(H)$. 

Let us first set the stage. As before, list all the marked points in $S_{b,m}^g$ as $\{ 1,2,\cdots, m \}$. We call a $d$-arc \emph{linear} (resp. \emph{cyclic}) if the marked points it passes through, in order, are given by $p,p+1,\cdots, p+d-1$ (resp. $p,p+1,\cdots, p+d-1 \mod m$) for some $p$. We will call $p$ the \emph{initial marked point} of the linear (resp. cyclic) $d$-arc. We define the \emph{linear $d$-arc matching complex} $\lmarc_d(S_{b,m}^g)$ (resp. \emph{cyclic $d$-arc matching complex} $\cmarc_d(S_{b,m}^g)$) to be the full subcomplex of the $d$-arc matching  complex $\marc_d(S_{b,m}^g)$ such that each vertex is a linear $d$-arc (resp. cyclic $d$-arc). Furthermore, for any subset $Z$ of $\{1,2,\cdots, m-d+1\}$, let $\lmarc_d(S_{b,m}^g,Z)$ be the full subcomplex of $\lmarc_d(S_{b,m}^g)$ spanned by those vertices whose initial point lies in $Z$. Similarly for any subset $Z$ of $\{1,2,\cdots,m\}$, let $\cmarc_d(S_{b,m}^g,Z)$ be the full subcomplex of $\cmarc_d(S_{b,m}^g)$ spanned by those vertices whose initial point lies in $Z$.

The proofs of the connectivity properties of $\lmarc_d(S_{b,m}^g)$ and $\cmarc_d(S_{b,m}^g)$ now follow closely to that of Theorem \ref{thm:conn-i-d-arc-complex} or \cite[Section 3.3]{BFM+16}. Let us focus on $\lmarc_d(S_{b,m}^g)$  first. Let $Z$ be any subset of  $\{1,2,\cdots, m-d+1\}$ with maximum $p_0$. We put a weight on all the marked points in $S_{b,m}^g$ via the following rule: if $p_0 \leq p\leq p_0+d-1$,  we assign its weight to be $2^{p-p_0}$; otherwise, we assign its weight to be $0$. With this we can define a height function $q$ on any linear $d$-arc  by assigning $q(\alpha)$ to be the total weight of the marked points $\alpha$ passes through.  Note that the zero set of $q$, which we denote by $\lmarc_d(S_{b+d,m-d}^g, Z)^{q=0}$,  can be identified with the complex $\lmarc_d(S_{b+d,m-d}^g,Z\setminus (\{p_0-d+1,\cdots, p_0\}\cap Z))$ as we have chosen $p_0$ to be the greatest in $Z$. Here the surface $S_{b+d,m-d}^g$ is obtained from $S_{b,m}^g$ by deleting a small open disk around those marked points with positive weight. Now $q$ defines a height function on the relative complex $(\lmarc_d(S_{b,m}^g,Z),\lmarc_d(S_{b,m}^g,Z)^{q=0})$. We will use $q$ to analyze the connectivity of $\lmarc_d(S_{b,m}^g, Z)$.

\begin{thm}\label{thm:conn-i-ldcomplex}
 For any $d\geq 2$, the complex $\lmarc_d(S_{b,m}^g,Z)$ is $(\lfloor\frac{|Z|-1}{3d-2}\rfloor-1)$-connected.
\end{thm}

\begin{rem}
In the theorem, the values of $b$ and $g$ do not play a role in our connectivity bound of $\lmarc_d(S_{b,m}^g,Z)$ whereas the value of $m$ only serves to give an upper bound on $|Z|$. Recall by definition of $Z$, if $z\in Z$, then  $z,z+1,\cdots,z+d-1$ are legitimate marked points, in particular $m\geq z+d-1$. 
\end{rem}

\Proof We prove the theorem by induction on $|Z|$. Note that as long as $|Z|>0$,  the complex  $\lmarc_d(S_{b,m}^g,Z) $ is nonempty and hence the theorem is valid when $|Z| < 3d-1$.  Now assume $|Z|\ge 3d-1$. Recall that the definition of our height function $q$ is based on the greatest $p_0\in Z$. Given any linear $d$-arc $\alpha$ in $\lmarc_d(S_{b,m}^g,Z)$ such that $q(\alpha)\neq 0$, the descending link of $\alpha$ is the full subcomplex of $\lmarc_d(S_{b,m}^g,Z)$ such that any vertex $\alpha'$ in it has the property that $q(\alpha') < q(\alpha)$ and $\alpha'$ is disjoint from $\alpha$. This complex can be identified with the linear disk complex with $\lmarc_d(S_{b+1,m-d}^g,Z')$ for some $Z'$, where $S_{b+1,m-d}^g$ is obtained from $S_{b,m}^g$ by cutting out a small open disk around $\alpha$. In the worst case, $q(\alpha)=1$ and $\alpha$ has an initial marked point $p_0-d+1$. In this case $Z'= Z \setminus (Z\cap \{p_0-2d+2,p_0-2d+2,\cdots, p_0\})$. Thus $|Z'|\geq |Z|-2d+1$. By induction,  $\lmarc_d(S_{b+1,m-d}^g,Z')$ is at least $(\lfloor\frac{|Z|-1}{3d-2}\rfloor-2)$-connected.

Now, as before, by the Morse lemma, the pair $(\lmarc_d(S_{b,m}^g,Z),\lmarc_d^{q=0}(S_{b,m}^g,Z))$ is $(\lfloor\frac{|Z|-1}{3d-2}\rfloor-1)$-connected, i.e. the inclusion ~$\iota\colon\lmarc_d^{q=0}(S_{b,m}^g,Z)\hookrightarrow \lmarc_d(S_{b,m}^g)$ induces an isomorphism in $\pi_n$ for $n \le \lfloor\frac{|Z|-1}{3d-2}\rfloor-2$ and an epimorphism for $n = \lfloor\frac{|Z|-1}{3d-2}\rfloor - 1$. On the other hand, by induction $\lmarc_d^{q=0}(S_{b,m}^g,Z)$ is at least $(\lfloor\frac{|Z|-1}{3d-2}\rfloor-2)$-connected. Hence $\marc_d(S_{b,m}^g,Z)$ is $(\lfloor\frac{|Z|-1}{3d-2}\rfloor-2)$-connected. But this is not enough as we want $\lmarc_d(S_{b,m}^g,Z)$ to be $(\lfloor\frac{|Z|-1}{3d-2}\rfloor-1)$-connected. Just as in the proof of Theorem \ref{thm:conn-i-d-arc-complex}, it is sufficient to show that $\pi_n(\lmarc_d^{q=0}(S_{b,m}^g,Z)\rightarrow\lmarc_d(S_{b,m}^g,Z))$ is trivial for $n\leq \lfloor\frac{|Z|-1}{3d-2}\rfloor-1$. In other words, we will show that when $n\leq \lfloor\frac{|Z|-1}{3d-2}\rfloor-1$, every map $\bar{\psi}: S^n\to \lmarc_d^{q=0}(S_{b,m}^g,Z)$ can be homotoped to a point in $\lmarc_d(S_{b,m}^g,Z)$.

We sketch how to proceed as in the proof of Theorem \ref{thm:conn-i-d-arc-complex}. First, we can apply Lemma~\ref{lem:injectifying} to $\psi = \iota\circ \bar{\psi}$ and assume $\psi$ is simplexwise injective. Now fix $\beta$ to be linear $d$-arc passing through the marked points $\{p_0,p_0+1,\cdots, p_0+d-1\}$. We claim that~$\psi$
can be homotoped in $\marc_d(S_{b,m}^g)$ to land in the star of $\beta$, which
will finish the proof. By assumption, none of the
$d$-arcs in the image of~$\psi$ will pass through positive valued marked points, but among the finitely
many such $d$-arcs, some might intersect nontrivially with ~$\beta$. Pick one, say $\alpha$,
intersecting $\beta$ at a point closest along $\beta$ to the marked point $p_0+d-1$, and let~$x$ be a
vertex of~$S^n$ mapping to $\alpha$. We now use the  mutual link trick (cf. Lemma \ref{lemma-replace-trick}) to push the intersection with $\beta$ towards the $p_0+d-1$ direction step by step. Note that our pushing direction is different than in Theorem \ref{thm:conn-i-d-arc-complex}.  In each step, we replace $\alpha$ by $\alpha'$ by pushing the intersection point closest to $p_0+d-1$ along $\beta$ towards $p_0+d-1$ across precisely one marked point. At the end, the image of $\psi$ will be disjoint from $\beta$. The only thing we need to worry about in order to do this is the connectivity of the mutual link. Let $D$ be the disk bounded by $\alpha$ and $\alpha'$ which contains one extra marked point $p'\in \{p_0+1,\cdots, p_0+d-1\}$ in its interior.  The mutual link again can be identified with $\lmarc_d(S_{b+1,m-d-1}^g,Z')$ for some subset $Z'$ of $Z$ where the surface $S_{b+1,m-d-1}^g$ is obtained from $S_{b,m}^g$ by cutting out a small open neighbourhood of $D\cup \alpha \cup \alpha'$. To obtain the subset $Z'$, we must remove any point from $Z$ which is a marked point that $\alpha$ crosses. Note that the initial marked point of $\alpha$ can be any point in $Z$. In the worst case, this results in removing the initial marked points in $Z\cap \{p_\alpha-d+1,p_\alpha-d+2,\cdots, p_\alpha+d-1\}, $ where $p_\alpha$ is the initial marked point of $\alpha$.  We also cannot have any $d$-arcs passing through $p'$. The worst case is when $p'=p_0+1$ which excludes linear $d$-arcs with an initial marked point $\{p_0-d+2,\cdots, p_0\}$. In total, we are throwing away at most $3d-2$  points in $Z$, thus by induction the mutual link is $( \lfloor\frac{|Z|-1}{3d-2}\rfloor-2)$-connected.
\qed

Taking $Z = \{1,2,\cdots, m-d+1\}$, we have the following.
\begin{cor}
For any $d\geq 2$, the complex $\lmarc_d(S_{b,m}^g) $ is $(\lfloor\frac{m-d}{3d-2}\rfloor-1)$-connected.
\end{cor}

Similarly, we have the following theorem.

\begin{thm}\label{thm:conn-i-cdcomplex}
 The complex $\cmarc_d(S_{b,m}^g,Z) $ is $(\lfloor\frac{|Z|-1}{3d-1}\rfloor-1)$-connected. In particular, the complex  $\cmarc_d(S_{b,m}^g)$ is $(\lfloor\frac{m-1}{3d-1}\rfloor-1)$-connected.
\end{thm}
\sProof The proof runs parallel to that of Theorem \ref{thm:conn-i-ldcomplex}. We can define a height function $q$ exactly as before except now there is no largest number $p_0\in Z$  as $Z$ is  cyclically ordered. So instead, we just pick an arbitrary $p_0$. This will affect the following calculations.
\begin{enumerate}
    \item The calculation of the descending link changes. Given any vertex $\alpha$ in the complex $\cmarc_d(S_{b,m}^g,Z)$ such that $q(\alpha)> 0$, the descending link of $
    \alpha$ can be identified with the complex $\cmarc_d(S_{b,m}^g,Z')$ for some $Z'$. In the worst case, $q(\alpha)=1$ and $\alpha$ has an initial marked point $p_0-d+1$. In this case $Z' = Z \setminus (Z\cap \{p_0-2d+2,p_0-2d+3,\cdots, p_0+d-1\})$. Thus $|Z'|\geq |Z|-3d+2$. By induction,  $\cmarc_d(S_{b,m}^g,Z')$ is at least $(\lfloor\frac{|Z|-1}{3d-1}\rfloor-2)$-connected.
    
    \item The calculation of the mutual link changes. Suppose for some point $x\in S^n$, its image $\psi(x) =\alpha$ intersects with $\beta$ nontrivially. We will replace $\alpha$ by $\alpha'$, where $\alpha'$ is obtained from $\alpha$ by pushing the intersection part along $\beta$ across one marked point. Let $p_\alpha$ be the initial point of $\alpha$. In the worst case, we have to remove from $Z$ any vertices with initial points in  $Z\cap \{p_\alpha-d+1,p_\alpha-d+2,\cdots, p_\alpha+d-1\}$. We also cannot allow the $d$-arcs which touch the marked point $p'\in \{p_0+1,\cdots,p_0+d-1\}$ in the disk bounded by $\alpha$ and $\alpha'$, i.e. vertices with initial marked points $p'-d+1,\cdots, p'$. In total, we are throwing away $3d-1$ elements in $Z$. Hence the mutual link is at least $(\lfloor\frac{|Z|-1}{3d-1}\rfloor-2)$-connected by induction.
    \end{enumerate}
Taking $Z$ to be the set of all marked points, we get the second part of the statement.
\qed

Similarly to how one defines the disk complex $\BD_d(S_{b,m}^g)$, one can define linear (resp. circular) disk complexes, $\BL\BD_d(S_{b,m}^g)$ (resp. $\BC\BD_d(S_{b,m}^g)$) by requiring the disks to enclose $d$ adjacent vertices ordered linearly (resp. circularly). In this case, one can consider the map $N$ defined in \ref{eq:mapfromarctodisk} but with the restricted domain of $\lmarc_d(S_{b,m}^g)$ (resp. $\cmarc_d(S_{b,m}^g)$). In either situation, the proof of Corollary~\ref{cor:conn-i-dcomplex} extends identically to give the following corollary.

\begin{cor}
The maps 
\[N|_{\lmarc_d}: \lmarc_d(S_{b,m}^g)\to \BL\BD_{ d}(S_{b,m}^g)\]
and
\[N|_{\cmarc_d}:\cmarc_d(S_{b,m}^g)\to \BC\BD_{ d}(S_{b,m}^g)\]
are complete joins. In particular, for any $d\geq 2$, the complex $\BL\BD_{ d}(S_{b,m}^g)$ is $(\lfloor\frac{m-1}{3d-2}\rfloor-1)$-connected and the complex $\BC\BD_{ d}(S_{b,m}^g)$ is $(\lfloor\frac{m-1}{3d-1}\rfloor-1)$-connected.
\end{cor}

\subsection{Finiteness of $H$ implies finiteness of braided Higman-Thompson groups} In this subsection, we  prove the ``if part" of Theorem \ref{thm-fin-bthomp} by studying the connectivity properties of the descending links in the Stein space $X$ with respect to the height function $f$. Recall that the descending link of a vertex~$x$ with
$f(x)=m$ is isomorphic to the complex $\elbraigecpx_d^m$ of dangling elementary $(d,m)$-braiges
$[(1_m,(b, \lambda),F_J^{(m-(d-1)|J|)})]$ with~$J\neq\emptyset$.  We will now construct a
projection from $\elbraigecpx_d^m$ to the $d$-marked-point-disk complex $\BD_d(S_{b,m}^g)$ and show it is a complete join. Since we have calculated the connectivity of $\BD_d(S_{b,m}^g)$ already, we can then apply our connectivity tools from Section \ref{subsec-com-join} to obtain the necessary connectivity of $\elbraigecpx_d^m$. We will wait until the end of the section to mention the ``pure''and ``cyclic'' cases.

Let $L_{m-1}$ be the \emph{linear graph
with}~$m$ vertices, that is the graph with $m$ vertices labeled~$1$ through~$m$, and $m-1$ edges, one connecting~$i$ to~$i+1$ for each~$1\le i<m$. Call a subgraph of $L_{m-1}$ a \emph{$d$-matching} on $L_{m-1}$ if each connected component of it is a subgraph of length $d-1$. Clearly, the set of $d$-matchings form a simplicial complex called the $d$ matching complex, denoted by $\CM_d(L_{m-1})$, where a matching forms a $k$-simplex whenever it consists of $k+1$ disjoint paths and the face relation is given by inclusion.

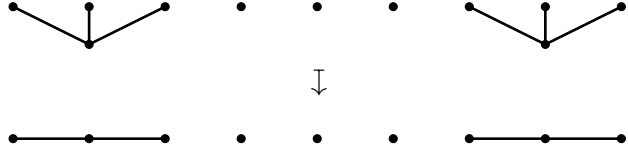
\begin{figure}[h]
\centering
\vskip2ex
\begin{tikzpicture}
  \filldraw
   (0,0) circle (1.5pt)
   (1,0) circle (1.5pt)
   (2,0) circle (1.5pt)
   (3,0) circle (1.5pt)
   (4,0) circle (1.5pt)
   (5,0) circle (1.5pt)
   (6,0) circle (1.5pt)
   (7,0) circle (1.5pt)
   (8,0) circle (1.5pt)

   (1,-0.5) circle (1.5pt)
   (7,-0.5) circle (1.5pt);

  \draw[line width=1pt]
  (0,0) -- (1, -0.5)
  (1,0) -- (1, -0.5)
  (2,0) -- (1, -0.5)
  (6,0) -- (7, -0.5)
  (7,0) -- (7, -0.5)
  (8,0) -- (7, -0.5);
  
  \node[rotate=-90] at (4,-1) {$\mapsto$};

 \begin{scope}[yshift=-1.75cm]
  \filldraw 
   (0,0) circle (1.5pt)
   (1,0) circle (1.5pt)
   (2,0) circle (1.5pt)
   (3,0) circle (1.5pt)
   (4,0) circle (1.5pt)
   (5,0) circle (1.5pt)
   (6,0) circle (1.5pt)
   (7,0) circle (1.5pt)
   (8,0) circle (1.5pt);
  
  \draw[line width=1pt]
   (0,0) -- (1,0)  -- (2,0) 
   (6,0) --  (7,0) -- (8,0);
 \end{scope}
\end{tikzpicture}
\vskip2ex
\caption{An example of the bijective correspondence between $3$-ary elementary forests
with~$9$ leaves and simplices of $\CM_3(L_8)$.}
\label{fig:matchings_to_forestsrep}
\end{figure}

We now observe that there is a bijection between the set of elementary $d$-ary forests with $m$ leaves and the set of $d$ matchings on $L_{m-1}.$ Under the identification, carets correspond to paths of length $d-1$. See Figure~\ref{fig:matchings_to_forestsrep} for an example. 

In light of the observation, we can denote an elementary~$(d,m)$-braige by
$((b, \lambda),\Gamma)$, where~$b\in B_m$, $\lambda$ is a labeling, and~$\Gamma$ is a $d$-matching on $L_{m-1}$. As usual, the equivalence class under dangling will be denoted $[(b, \lambda),\Gamma]$.

Let~$S=S_{1,m}^0$ be the unit disk with $m$ marked points given by fixing an embedding $L_{m-1} \hookrightarrow S$ of the linear graph with~$m-1$ edges into~$S_{1,0}^0$. With these data in place we can consider~$\marc_d(S)$, the $d$-arc matching complex on $S$,
and we have an induced embedding of simplicial complexes $\CM_d(L_{m-1})
\hookrightarrow \marc_d(S)\rightarrow \BD_d(S)$ where the second map is the map $N$ given in equation~\ref{eq:mapfromarctodisk}.  The braid group $B_m$ on $m$ strands is
isomorphic to the mapping class group of the disk with $m$ marked points~\cite{Bir74}, so we have an action of $B_m$
on~$\BD_d(S)$.  We will consider this as
a right action (in the same way as dangling is a right action on braiges), so for~$b\in
B_m$ and $\sigma\in\BD_d(S)$ we will write~$(\sigma)b$ to denote the image
of~$\sigma$ under~$b$.

Define a map $\pi$ from $\elbraigecpx_d^m$ to $\BD_d(S)$ as follows. We
view~$\CM_d(L_{m-1})$ as a subcomplex of~$\marc_d(S)$, so we can associate
to any $d$ matching $\Gamma$ a simplex in $\BD_d(S)$ denoted $\Gamma_0$. Thus we can map any
elementary~$(d,m)$-braige~$((b, \lambda),\Gamma)$ to the simplex $(\Gamma_0)b^{-1}$ in
$\BD_d(S)$, forgetting the label $\lambda$, taking a tubular neighborhood of $\Gamma$ to get a set of disks $\Gamma_0$, and then applying $b^{-1}$ to $\Gamma_0$. 

\begin{lem}
The map
\begin{align*}
\pi \colon \elbraigecpx^m_d & \to \BD_d(S) \\
[((b,\lambda),\Gamma)] & \mapsto (\Gamma_0)b^{-1} \text{ .}
\end{align*}
is well defined on equivalence classes.
\end{lem}

\proof
Suppose $\Gamma$ corresponds to a $q-1$ simplex of $\CM_d(L_m)$, i.e. $\Gamma$ is a elementary forest with $q$ carets. Let $(b_1, \lambda_1)$ be in $B_q(H)$ so that $[((b,\lambda),\Gamma)]=[((b,\lambda),\Gamma)(1_q, (b_1, \lambda_1), 1_q)]$. Let $\Gamma'$ and $(b_1', \lambda_1')$ be such that
\[((b,\lambda),\Gamma)(1_q, (b_1, \lambda_1), 1_q)=((b, \lambda)(b_1', \lambda_1'), \Gamma')=(bb_1', \lambda^{b_1'}\lambda_1', \Gamma').\]
Under $\pi$, this maps to $(\Gamma'_0)(bb_1')^{-1}$. We need to show this is the same as $(\Gamma_0)b^{-1}$ or equivalently that $(\Gamma_0')(b_1')^{-1}=\Gamma_0$ or $(\Gamma_0)b_1'=\Gamma_0.$ Note that $b_1'$ is obtained from $b_1$ by turning each strand that corresponds to a root of $\Gamma$ into $d$ parallel strands and then braiding them according to the appropriate label. Since these local braids are supported on the interior of individual disks, they do not change the disk system and we conclude $(\Gamma_0)b_1'=\Gamma_0$. 
\qed 

By construction, the map $((b, \lambda),\Gamma)\mapsto (\Gamma_0)b^{-1}$ is
well defined on equivalence classes under dangling, so we obtain a simplicial
map \begin{align*}
\pi \colon \elbraigecpx^m_d & \to \BD_d(S) \\
[((b,\lambda),\Gamma)] & \mapsto (\Gamma_0)b^{-1} \text{ .}
\end{align*}
Note that~$\pi$ is surjective, but not injective at all.

\begin{figure}[h]
\centering
\begin{tikzpicture}
    \clip(-.30,0.4) rectangle (2.75,-2.95);

    \draw[line width=1pt] (0,0) to (0,-2);
    \draw[line width=1pt] (.5,0) to (.5,-2);
    
     \draw[line width=3pt, white] (1,0) to [out=270, in=90] (1.5,-2);
    \draw[line width=1pt] (1,0) to [out=270, in=90] (1.5,-2);
    
    \draw[line width=1pt] (2,0) to (2,-2);
    \draw[line width=1pt] (2.5,0) to (2.5,-2);
 
   \draw[line width= 1pt] (0,-2)-- (.5, -2.5)--(.5, -2);
   \draw[line width=1pt] (.5, -2.5) --(1,-2);
 
    \draw[line width= 1pt] (1.5,-2)-- (2, -2.5)--(2, -2);
   \draw[line width=1pt] (2, -2.5) --(2.5,-2);

        \draw[line width=3pt, white] (1.5,0) to [out=270, in=90] (1,-2);
    \draw[line width=1pt] (1.5,0) to [out=270, in=90] (1,-2);

         \draw[white,fill=white] (0,-1.6) circle (5 pt);
\draw  (0,-1.6) circle (5 pt) node[text=black, scale=.5] {$f_1$}; 

      \draw[white,fill=white] (0.5,-1.6) circle (5 pt);
\draw  (0.5,-1.6) circle (5 pt) node[text=black, scale=.5] {$f_2$}; 

       \draw[white,fill=white] (1.02,-1.6) circle (5 pt);
\draw  (1.02,-1.6) circle (5 pt) node[text=black, scale=.5] {$f_3$};  
   
        \draw[white,fill=white] (1.45,-1.6) circle (5 pt);
\draw  (1.45,-1.6) circle (5 pt) node[text=black, scale=.5] {$f_4$};  

       \draw[white,fill=white] (2,-1.6) circle (5 pt);
\draw  (2,-1.6) circle (5 pt) node[text=black, scale=.5] {$f_5$};  

       \draw[white,fill=white] (2.5,-1.6) circle (5 pt);
\draw  (2.5,-1.6) circle (5 pt) node[text=black, scale=.5] {$f_6$};

   \filldraw[]
   (0, -2) circle (1pt)
   (.5, -2) circle (1pt)
   (1, -2) circle (1pt)
   (1.5, -2) circle (1pt)
   (2, -2) circle (1pt)
   (2.5, -2) circle (1pt);

\draw[line width=1pt] (-0.25,0) -- (2.75,0);
    
\end{tikzpicture}
\quad
\begin{tikzpicture}
    \clip(-.30,0.4) rectangle (2.75,-2.95);

        \draw[line width=1pt] (0,0) to (0,-2.5);
    \draw[line width=1pt] (.5,0) to (.5,-2.5);

     \draw[line width=3pt, white] (1,0) to [out=270, in=90] (1.5,-2.5);
    \draw[line width=1pt] (1,0) to [out=270, in=90] (1.5,-2.5);
    
     \draw[line width=3pt, white] (1.5,0) to [out=270, in=90] (1,-2.5);
    \draw[line width=1pt] (1.5,0) to [out=270, in=90] (1,-2.5);
    
    \draw[line width=1pt] (2,0) to (2,-2.5);
    \draw[line width=1pt] (2.5,0) to (2.5,-2.5);
 
   \filldraw[]
   (0, -2.5) circle (1pt)
   (.5, -2.5) circle (1pt)
   (1, -2.5) circle (1pt)
   (1.5, -2.5) circle (1pt)
   (2, -2.5) circle (1pt)
   (2.5, -2.5) circle (1pt);
   
   \draw[line width=1pt, ] (0,-2.5) -- (1,-2.5);
   \draw[line width=1pt, ] (1.5,-2.5) -- (2.5,-2.5);

\draw[line width=1pt] (-0.25,0) -- (2.75,0);
\end{tikzpicture}

\quad
\begin{tikzpicture}
    \clip(-.30,0.4) rectangle (2.75,-2.95);
       \draw[line width=1pt] (-.125,-2.5) to [out=270, in=270] (1.125,-2.5);
       \draw[line width=1pt] (1.125,-2.5) to [out=90, in=90] (-.125,-2.5);

       \draw[line width=1pt] (1.375,-2.5) to [out=270, in=270] (2.625,-2.5);
       \draw[line width=1pt] (2.625,-2.5) to [out=90, in=90] (1.375,-2.5);
    
        \draw[line width=3pt, white] (0,0) to (0,-2.5);
        \draw[line width=1pt] (0,0) to (0,-2.5);
    \draw[line width=3pt, white] (.5,0) to (.5,-2.5);
    \draw[line width=1pt] (.5,0) to (.5,-2.5);

     \draw[line width=3pt, white] (1,0) to [out=270, in=90] (1.5,-2.5);
    \draw[line width=1pt] (1,0) to [out=270, in=90] (1.5,-2.5);
    
     \draw[line width=3pt, white] (1.5,0) to [out=270, in=90] (1,-2.5);
    \draw[line width=1pt] (1.5,0) to [out=270, in=90] (1,-2.5);
        
        \draw[line width=3pt, white] (2,0) to (2,-2.5);
    \draw[line width=1pt] (2,0) to (2,-2.5);
    
        \draw[line width=3pt, white] (2.5,0) to (2.5,-2.5);
    \draw[line width=1pt] (2.5,0) to (2.5,-2.5);
 
   \filldraw[]
   (0, -2.5) circle (1pt)
   (.5, -2.5) circle (1pt)
   (1, -2.5) circle (1pt)
   (1.5, -2.5) circle (1pt)
   (2, -2.5) circle (1pt)
   (2.5, -2.5) circle (1pt);
   
   \draw[line width=1pt] (-0.25,0) -- (2.75,0);

\end{tikzpicture}

\quad
\begin{tikzpicture}
 \clip(-.30,0.4) rectangle (2.75,-3.95);
 
\draw[line width=1pt] (-.125,-2.5) to [out=90, in=180] (.6,-2.35) to [out=0, in=90] (.75,-2.75)to [out=270, in=180] (1.125,-2.9) to [out=0, in=270] (1.35, -2.375) to [out=90, in=90] (1.65,-2.375);
\draw[line width=1pt] (-.125,-2.5) to [out=270, in=180] (.75,-3.1)to [out=0, in=180] (1.325,-3.1) to [out=0, in=270] (1.65,-2.375);

\begin{scope}[xscale=-1, xshift=-2.5cm, yscale=-1, yshift=5cm]
 \draw[line width=1pt] (-.125,-2.5) to [out=90, in=180] (.6,-2.35) to [out=0, in=90] (.75,-2.75)to [out=270, in=180] (1.125,-2.9) to [out=0, in=270] (1.35, -2.375) to [out=90, in=90] (1.65,-2.375);
\draw[line width=1pt] (-.125,-2.5) to [out=270, in=180] (.75,-3.1)to [out=0, in=180] (1.325,-3.1) to [out=0, in=270] (1.65,-2.375);
\end{scope}

    \draw[line width=3pt, white] (0,0) to (0,-2.5);
    \draw[line width=1pt] (0,0) to (0,-2.5);
    \draw[line width=3pt, white] (.5,0) to (.5,-2.5);
    \draw[line width=1pt] (.5,0) to (.5,-2.5);
    \draw[line width=3pt, white] (1,0) to (1,-2.5);
    \draw[line width=1pt] (1,0) to (1,-2.5);
    \draw[line width=3pt,white] (1.5,0) to (1.5,-2.5);
    \draw[line width=1pt] (1.5,0) to (1.5,-2.5);
    \draw[line width=3pt, white] (2,0) to (2,-2.5);
    \draw[line width=1pt] (2,0) to (2,-2.5);
    \draw[line width=3pt, white] (2.5,0) to (2.5,-2.5);
    \draw[line width=1pt] (2.5,0) to (2.5,-2.5);

  \filldraw[]
   (0, -2.5) circle (1pt)
   (.5, -2.5) circle (1pt)
   (1, -2.5) circle (1pt)
   (1.5, -2.5) circle (1pt)
   (2, -2.5) circle (1pt)
   (2.5, -2.5) circle (1pt);

\draw[line width=1pt] (-0.25,0) -- (2.75,0);

\end{tikzpicture}
\caption{From left to right the pictures show the process of flattening to a $3$-matching, forgetting the labels, expanding the matching to disks, and ``combing straight'' the braid.}
\label{fig:BMD-to-arc-cplx}
\end{figure}
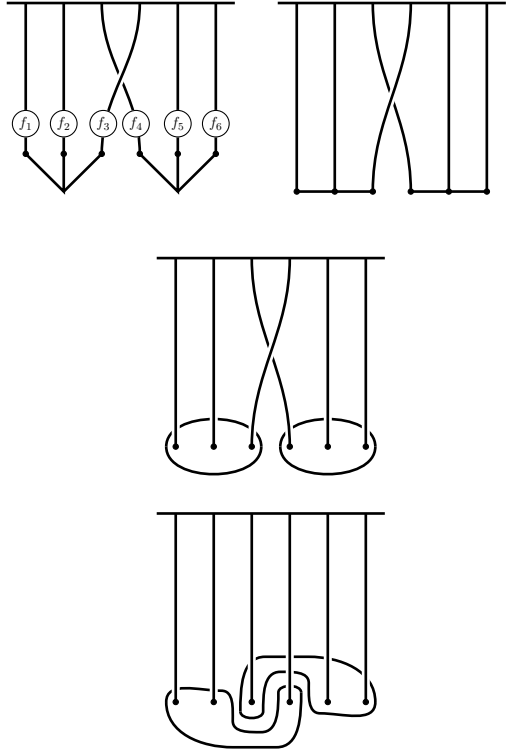

One can visualize this map by first forgetting the labels and considering the merges as $d$-disks, then ``combing
straight'' the braid and seeing where the d-disks are taken, as in
Figure~\ref{fig:BMD-to-arc-cplx}. Note that the resulting simplex $(\Gamma_0)b^{-1}$
of $\BD_d(S)$ has the same dimension as the simplex~$[((b, \lambda),\Gamma)]$ of
$\elbraigecpx_d^m$.

The next lemma and proposition are concerned with the fibers of~$\pi$.

\begin{lem}\label{lem:spx_fiber}
 Let~$E$ and~$\Gamma$ be simplices in~$\CM_d(L_{m-1})$, such that~$E$ is a $0$-simplex and~$\Gamma$ is an~$e(\Gamma)$-simplex. Let $[((b_1, \lambda_1),E)]$ and $[((b_2,\lambda_2), \Gamma)]$ be dangling elementary~$(d,m)$-braiges. Suppose that their images under the map~$\pi$ are contained in a simplex of $\BD_d(S)$.  Then there exists a simplex in $\elbraigecpx_d^m$ that contains~$[((b_1, \lambda_1),E)]$ and~$[((b_2,\lambda_2),\Gamma)]$.
\end{lem}

\Proof
We may assume that $[((b_1,\lambda_1), E)]$ is not contained in $[((b_2,\lambda_2),\Gamma)]$. There is an action of $B_m(H)$ on $\elbraigecpx_d^m$ (``from above''), given by
\[(b,\lambda)[((b',\lambda'),\Gamma')]=[((b,\lambda)(b',\lambda'),\Gamma')].\]
One can check that for each $k\ge 0$, this
action is transitive on the~$k$-simplices of~$\elbraigecpx_d^m$. We can therefore assume
without loss of generality that $(b_2, \lambda_2)=(\id, \iota)$ where $ \iota$ is the trivial labeling, and $\Gamma$ is the $d$-matching
of~$L_{m-1}$ whose components are precisely those subgraphs of length $d-1$ with starting points
$j \in \{1,d+1,\dots,de(\Gamma)+1\}$.

There is a $d$-disk $\alpha$ representing~$\pi([((b_1, \lambda_1),E)])$ that is disjoint from $\Gamma$. This
disjointness ensures that, after dangling, we can assume the following condition
on~$(b_1,\lambda_1)$:  for each component of $\Gamma$ say with endpoints $j$ and $j + d-1$, where $j \in \{1,d+1,\dots,de(\Gamma)+1\}$, $b$ can be represented as a
braid in such a way that from the $j$th to $(j+d-1)$st strands of $b$ run straight down, parallel
to each other,  no strands cross between them, and the labels on them are all trivial.  In particular, $[((b_1, \lambda_1),\Gamma)]=[((\id,\iota),\Gamma)]$,
so $[((b,\lambda_1),\Gamma\cup E)]$ is a simplex in $\elbraigecpx_d^m$ with $[(b_1,\lambda_1, E)]$
and~$[((\id,\iota), \Gamma)]$ as faces.
\qed

\begin{prop}\label{prop:fibers}
The map  $\pi \colon \elbraigecpx^m_d  \to \BD_d(S)$ is a complete join.
\end{prop}
\Proof
We have already seen that the map $\pi$ is surjective and injective on individual simplices. Let $\sigma$ be a $k$-simplex in $\BD_d(S)$ with vertices $v_0, \ldots, v_k$. To prove $\pi$ is a complete join, it just remains to show 
\begin{equation*}
\pi^{-1}(\sigma)=\Bigast_{j=0}^k\pi^{-1}(v_j).
\end{equation*}

``$\subseteq$'':\quad This inclusion just says that vertices in $\pi^{-1}(\sigma)$ that are connected by an edge map to distinct vertices under $\pi$ which is clear.

``$\supseteq$'':\quad We prove this by induction on $k$. Suppose $\pi^{-1}(\sigma) \supseteq \Bigast_{j=0}^k\pi^{-1}(v_j)$ for $k=r$. Now given an $(r+1)$-simplex $\sigma=\langle v_0,\cdots, v_{r+1}\rangle$ which is a join of $\tau =\langle v_0,\cdots, v_{r}\rangle$ and $v_{r+1}$. We just need to show for any simplex $\bar{\tau},\bar{v}_{r+1} \in \elbraigecpx_d^m$ such that $\pi(\bar{\tau})=\tau, \pi\bar{v}_{r+1}=v_{r+1}$, we have a $(r+1)$-simplex contains both $\bar{\tau}$ and $v_{r+1}$. But this is exactly Lemma \ref{lem:spx_fiber}.
\qed

\begin{cor}
\label{cor:desc_link_conn}
The complex $\elbraigecpx_d^m$ is $(\lfloor\frac{m+1}{2d-1}\rfloor-2)$-connected.  Hence for any vertex~$x$ in $X$ with~$f(x)=m$, $\dlk(x)$ is
$(\lfloor\frac{m+1}{2d-1}\rfloor-2)$-connected.
\end{cor}

\Proof
We know that $\BD_{d}(S)$ is $(\lfloor\frac{m+1}{2d-1}\rfloor-2)$-connected according to Corollary~\ref{cor:conn-i-dcomplex}.  For any
$k$-simplex $\sigma$ in $\BD_{d}(S)$,  $\Lk(\sigma)$ is isomorphic to $\BD_d(S_{k+2,m-d(k+1)}^{0})$, which is $(\lfloor\frac{m-d(k+1)+1}{2d-1}\rfloor-2)$-connected, hence at least $(\lfloor\frac{m+1}{2d-1}\rfloor-2 -(k+1))$-connected. Thus  $\BD_{d}(S)$ is wCM of dimension $\lfloor\frac{m+1}{2d-1}\rfloor-1$. Since $\pi:\elbraigecpx_d^m \to \BD_{d}(S)$ is a complete join, by Proposition \ref{prop-join-conn},  $\elbraigecpx_d^m $ is wCM of dimension $\lfloor\frac{m+1}{2d-1}\rfloor-1$. In particular, it is $(\lfloor\frac{m+1}{2d-1}\rfloor-2)$-connected.
\qed

In the other cases, we consider the descending links of vertices in
$X(bF_{d,r}(H))$ and $X(bT_{d,r}(H))$.  For a vertex~$x$ with~$m$ feet, $\dlk(x)$ is isomorphic
to $\elpbraigecpx_d^{m}$ or to $\elcbraigecpx_d^{m}$, respectively.  These project onto the complexes $\BL\BD_d(S)$ and $\BC\BD_d(S)$. Using the same argument, we have the following.

\begin{cor}\label{cor:desc_link_conn-pc}
The complex  $\elpbraigecpx_d^{m}$ is $(\lfloor\frac{m-d}{3d-2}\rfloor-2)$-connected.  Hence for any vertex~$x$ in $X(bF_{d,r}(H))$ with~$f(x)=m$, the descending link $\dlk(x)$ is
$(\lfloor\frac{m-d}{3d-2}\rfloor-1)$-connected. The complex $\elcbraigecpx_d^{m}$ is $(\lfloor\frac{m-1}{3d-2}\rfloor-2)$-connected.  Hence for any vertex~$x$ in $X(bT_{d,r}(H))$ with~$f(x)=m$, the descending link $\dlk(x)$ is
$(\lfloor\frac{m-1}{3d-1}\rfloor-2)$-connected.
\end{cor}

Combining these with the Morse lemma, we obtain the following.

\begin{cor}\label{cor:connectofpair}
For any $k\geq 0$, the filtration $\{X^{\leq m}\}_m$ is essentially $k$-connected. The same is also true for the filtration $\{X(bF_{d,r}(H))^{\leq m}\}_m$ and $\{X(bT_{d,r}(H))^{\leq m}\}_m$.
\end{cor}
\Proof 
By the Morse lemma (Lemma \ref{lemm-Morse} (2)) and Corollary \ref{cor:desc_link_conn}, we have for $m \geq 1$ the pair $(X, X^{\leq m-1})$ is $(\lfloor\frac{m+1}{2d-1}\rfloor-1)$-connected.   On the other hand, by Corollary \ref{cor-steincontract}, $X$ is contractible. This means  for any $m$ such that $\lfloor\frac{m+1}{2d-1}\rfloor-2\geq k$, we have $\pi_k(X^{\leq m})$ is trivial. Therefore, for any $k\geq 1$, the filtration $\{X^{\leq m}\}_m$ is essentially $k$-connected. The same argument implies $\{X(bF_{d,r}(H))^{\leq m}\}_m$ and $\{X(bT_{d,r}(H))^{\leq m}\}_m$ are also essentially $k$-connected for any $k\geq 0$.
\qed

We are now ready to prove the ``if part" of Theorem \ref{thm-fin-bthomp}.

\begin{thm}\label{thm-fin-if}
If $H$ is of type $F_n$, then the groups $bV_{d,r}(H)$, $bF_{d,r}(H)$, and $bT_{d,r}(H)$ are also of type $F_n$. 
\end{thm}

\Proof
Suppose that $H$ is of type $F_n$. Consider the actions of $bV_{d,r}(H)$, $bF_{d,r}(H)$, and $bT_{d,r}(H)$ on the corresponding Stein spaces which are connected by Corollary \ref{cor-steincontract}. Abusing notation, we will denote all of the Stein spaces by $X$. By Corollary \ref{cor:cell_stabs}, all of the cell stabilizers are of type $F_n$ and by Proposition \ref{prop:stein_space_cible}, each $X^{\leq m}$ is finite modulo the corresponding group. Finally, by Corollary \ref{cor:connectofpair}, the filtration $\{X^{\leq m}\}_m$ is essentially $k$-connected for any $k\geq 0$. We conclude, by Brown's criterion (Theorem \ref{brownscriterion}), that if $H$ is of type $F_n$ then so are each of groups $bV_{d,r}(H)$, $bF_{d,r}(H)$, and $bT_{d,r}(H)$.
\qed

\subsection{Quasi-retracts and finiteness properties}
The purpose of this subsection is to show that the group $bV_{d,r}(H)$ (resp. $bF_{d,r}(H)$, $bT_{d,r}(H)$) is not of type $F_n$ if $H$ is not. The proof is inspired by \cite[Section 4]{BZ20}. Basically, we will prove that $H$ is a quasi-retract of $bV_{d,r}(H)$ and $bT_{d,r}(H)$ and a retract of $bF_{d,r}(H)$.

Recall first that a group $Q$ is called a retract of a group $G$ if there is a pair of group homomorphisms
\[ Q \overset{i}{\hookrightarrow} G \overset{r}{\twoheadrightarrow} Q\]
such that $r\circ i$ is the identity on $Q$. Suppose $Q$ is a retract of $G$. Then if $G$ is of type $F_n$, so is $Q$, see for example \cite[Proposition 4.1]{Bux04}.  The same holds if one replaces retract by quasi-retract. Let us make this precise. Recall a function $f:X\to Y$ is said to be coarse Lipschitz if there exists constants $C,D>0$ so that 
$$d(f(x),f(x')) \leq Cd(x,x')+D \text{ for all } x,x'\in X$$
For example, any homomorphism between finitely generated groups is coarse Lipschitz with respect to the word metrics.
A function $\rho:X\to Y$ is said to be a quasi-retraction if it is coarse Lipschitz and there exists a coarse Lipschitz function $\iota:Y\to X$ and a constant $E>0$ so that $d(\rho\circ \iota(y),y)\leq E$ for all $y\in Y$. If such a function exists, $Y$ is said to be a quasi-retract of $X$. 

\begin{thm}\cite[Theorem 8]{Alo94}\label{thm-quasi-re-fin}
Let $G$ and $Q$ be finitely generated groups such that $Q$ is a quasi-retract of $G$ with respect to word metrics corresponding to some finite generating sets. Then if $G$ is of type $F_n$, so is $Q$. 
\end{thm}

Now let us define a map $\iota_F: H \rightarrow bF_{d,r}(H)$ via $h\mapsto [1_r, (\id, \lambda_h), 1_r]$ where $1_r$ is the trivial forest and $\lambda_h$ labels all the strings by $h$. Since $bF_{d,r}(H)\leq bT_{d,r}(H)\leq bV_{d,r}(H)$, we also have maps $\iota_V: H \rightarrow bV_{d,r}(H)$ and $\iota_T: H \rightarrow bT_{d,r}(H)$. We define another map $r_V:bV_{d,r}(H)\rightarrow H$ given by $[F_-,(b, \lambda), F_+]\mapsto \lambda(1)$. Restricting $r_V$ to $bF_{d,r}(H)$ and $bT_{d,r}(H)$, we get the maps $r_F$ and $r_T$. Note that only $r_F$ is a group homomorphism. One easily checks that $r_F\circ \iota_F =\id$. Thus we have the following.
\begin{lem}
The group $H$ is a retract of $bF_{d,r}(H)$.
\end{lem}

We do also have $r_V\circ \iota_V =\id$ and  $r_T\circ \iota_T =\id$. But since $r_V$ and $r_T$ are not group homeomorphisms now, the best we can hope for is that they are coarse Lipschitz. To prove this, we first need an understanding of the generating set. Let $T_1$ be a $(d,r)$-forest such that the first tree is a single caret and all other trees are trivial. Let $\iota': H\to bF_{d,r}(H)$ be the inclusion sending $ h $ to $[T_1,(id,\lambda_h'),T_1]$, where $\lambda_h'$ labels the first string by $h$ and all other strings by $1\in H$. Note that $\iota'(H)$ naturally sits in $bT_{d,r}(H)$ and $bV_{d,r}(H)$. On the other hand, we have $bT_{d,r}\leq bT_{d,r}(H)$ and $bV_{d,r}\leq bV_{d,r}(H)$ using the trivial labels on all strings.

\begin{prop}\label{prop-fg}
The group $bV_{d,r}(H)$ is generated by $\iota'(H)$ and $bV_{d,r}$. Similarly, the group $bT_{d,r}(H)$ is generated by $\iota'(H)$ and $bT_{d,r}$.
\end{prop}
\Proof 
We prove the proposition for $bV_{d,r}(H)$. The other case is similar. Let $G$ be the subgroup generated by $\iota'(H)$ and $bV_{d,r}$, we prove $G = bV_{d,r}(H)$ in four steps.

\begin{enumerate}[label= Step \arabic*.]
    \item Let $F$ be a $(d,r)$-forest such that the first leaf has distance $1$ to the root of the tree it is part of. Then for any $k\geq 1$, elements of the form $[F,(id,\lambda_h^k),F] $ lie in $G$, where $\lambda_h^k$ labels the $k$-th string  by $h$ and all others by $1$. In fact, let $b$ be any braid whose corresponding element in the symmetric group permutes $1$ and $k$, then  $[F,(id,\lambda_h^k),F]  = [F,(b,\lambda_0),F]^{-1}\iota'(h) [F,(b,\lambda_0),F] \in G$, where $\lambda_0$ here is the trivial labeling.
    
    \item Let $F$ by any $(d,r)$-forest, and $\lambda_h$ be a labeling of the strings such that only one string is labeled nontrivially and it is labeled by $h$, then $[F,(id,\lambda_h),F] $ lies in $ G$. If the initial leaf of the string labeled nontrivially does not lie below the leftmost vertex that has distance $1$ to the root, it is already covered by step 1. If not, we can choose any element as in step 1, and conjugate it to $[F,(id,\lambda_h),F] $ by an element in $bV_{d,r}$ using the same strategy.
    
    \item Let $F$ by any $(d,r)$-forest, and $\lambda$ be any labeling, then $[F,(id,\lambda),F] \in G$. In fact, let $\lambda^k_h$ be the labeling of the strings such that  the $k$-th string is labeled by $h$ and all other string are labeled by $1$. Then $[F,(id,\lambda),F] \in G$ is a product of $[F,(id,\lambda_h^k),F]$.
    
    \item Finally, let $[F,(b,\lambda),F']$ be any element of $bV_{d,r}(H)$, then $$[F,(b,\lambda),F'] = [F,(b,\lambda_0),F'][F',(id,\lambda),F'],$$ 
    where again $\lambda_0$ is the trivial labeling. Since $[F,(b,\lambda_0),F']\in bV_{d,r} \leq G$ and  $[F,(id,\lambda),F'] \in G$, we have  $[F,(b,\lambda),F'] \in G$.
\end{enumerate}
\qed

\begin{thm}\label{thm-quasi-re}
The group $H$ is a quasi-retract of $bV_{d,r}(H)$ and $bT_{d,r}(H)$.
\end{thm}
\proof We prove the theorem for  $bV_{d,r}(H)$.
Fix  finite generating sets $S_H$ for $H$, and $S_V$ for $bV_{d,r}$. By Proposition \ref{prop-fg}, $\iota'(S_H)\cup S_V$ is a finite generating set of $bV_{d,r}(H)$. We will show that the map  $r_V:  bV_{d,r}(H)  \rightarrow H$ is coarse Lipschitz with respect to the word metric on $bV_{d,r}(H)$ and $H$. Now
\begin{enumerate}
    \item  $r_V( g\iota'(s))\in \{r_V(g), r_V(g)s\}$ for all $s\in \iota'(S_H)$ and $g\in bV_{d,r}(H)$, and 
    \item  $r_V(gg') = r_V(g)$ for any $g'\in bV_{d,r}$ and $g\in bV_{d,r}(H)$.
\end{enumerate}
It follows that $r_V$ is nonexpanding and hence coarse Lipschitz. Since $\iota_V$ is a group homomorphism, it must be coarse Lipschitz as well. As $r_V\circ \iota_V = \id_H$, we conclude that $r_V$ is a quasi-retraction. The proof for $bT_{d,r}(H)$ is exactly the same.
\qed

\begin{thm} \label{thm-fin-bthomp}
For any $d\geq 2$ and $r\geq 1$  and any subgroup $H$ of the braid group $B_d$ (resp. of the pure braid group $PB_d$), the group $bV_{d,r}(H)$ (resp.  $bT_{d,r}(H)$ or $bF_{d,r}(H)$) is of type $F_n$ if and only if $H$ is.
\end{thm}
\Proof  For $n \geq 2$, the theorem is immediate from Theorems \ref{thm-fin-if}, \ref{thm-quasi-re-fin} and \ref{thm-quasi-re}.

For $n=1$, the only thing we need to prove is that the group $bV_{d,r}(H)$ (resp. $bF_{d,r}(H)$, $bT_{d,r}(H)$), is finitely generated, then $H$ is also finitely generated. Suppose $H$ is not finitely generated, then we have a sequence of  proper subgroups $H_1\lneq\cdots H_i\lneq H_{i+1}\lneq \cdots$ of $H$ such that $\cup_iH_i=H$. Then we have a sequence of proper subgroups $bV_{d,r}(H_1)\lneq\cdots bV_{d,r}(H_i)\lneq bV_{d,r}(H_{i+1})\lneq \cdots$ of  $bV_{d,r}(H)$ such that $\cup_i bV_{d,r}(H_i)=bV_{d,r}(H)$. This shows $bV_{d,r}(H)$ is not finitely generated.
\qed

Note that if $H$ is the trivial group, then the groups $bV_{d,r}(H)$, $bF_{d,r}(H)$, and $bT_{d,r}(H)$ are the braided Higman--Thompson groups $bV_{d,r}$, $bF_{d,r},$ and $bT_{d,r}$. Hence, we have the following immediate corollary.

\begin{cor}\label{cor-BVBTBFtypeFinfty}
The braided Higman--Thompson groups  $bV_{d,r}$, $bF_{d,r},$ and $bT_{d,r}$ are of type $F_{\infty}$.
\end{cor}

Similarly, taking $C$ to be the subgroup of $B_d$ generated by the half Dehn twist (resp. a full Dehn twist) around the boundary, we see that following Proposition \ref{lem-idf-rb-br}, the same is true for the ribbon Higman--Thompson groups.

\begin{cor}\label{cor-rthompFinfty}
The ribbon Higman--Thompson group $RV_{d,r}$ is of type $F_\infty$. Likewise, the oriented ribbon Higman--Thompson groups $RV^+_{d,r},RF^+_{d,r},$ and $RT^+_{d,r}$ are of type $F_\infty$.
\end{cor}

\bibliographystyle{alpha}
\bibliography{references.bib}

\end{document}